\documentclass{article}
\usepackage{amsmath}
\usepackage{graphicx}
\usepackage{xcolor}

\begin{document}
\title{Rodrigues, Olinde: \lq\lq Des lois g\'{e}om\'{e}triques
qui r\'{e}gissent les d\'{e}placements d'un syst\`{e}me
solide..", translation and commentary}

\author{Richard Friedberg}

\maketitle

\abstract{I provide for the first time in English a line-by-line 
translation of the entire text of the monumental 1840 memoir 
of Olinde Rodrigues, \lq\lq On the geometrical laws 
governing the motions of a solid system..." published in French in 
the  Journal de Math\'{e}matiques Pures et Appliqu\'{e}es.
I accompany the translation with copious notes in {\it italics},
 in which I explicate some passages whose
meaning is obscure in the direct translation, supply detailed
proofs, where lacking, of assertions in the original, and clarify
the overall organization of the memoir and the relation of its
sections to one another.  (In my notes, Rodrigues himself is
consistently called \lq\lq the author".)  I often supply a rendering
in modern vector notation, for equations and formulas in the
original text in which vectors are expressed laboriously in terms 
of their projections on fixed axes.}

{\it Benjamin Olinde Rodrigues (1795-1851) was a successful banker 
and ardent social critic descended from a Sephardic family 
long resident in France.  He associated himself for many years
with the utopian philosophy of Henri de Saint-Simon, supporting
the latter's movement both financially and administratively. Physics
students beyond the introductory level are likely to know of his
formula for the Legendre polynomial of any order, but unlikely to
know anything about its author.   This formula was contained in
his 1815 thesis as Docteur \`{e}s-Sciences at the University of Paris,
a thesis containing much material of far greater depth. Among his
subsequent writings for over twenty years are articles favoring the
equality of women, articles on Saint-Simonism, articles on the
reform of the banking laws, and other social issues, but none on
science or mathematics.  Then in the last eight years of his life
he published a number of mathematical papers, including the
magnum opus which is the subject of this translation.

The paper is a double tour-de-force in which the same material
is expounded first geometrically in Sections 1-14 and then
algebraically in Sections 15-22.  In reading the second part,
one will necessarily experience some d\'{e}ja-vu, but one should
keep in mind that although the same assertions or formulas
reappear, the logic binding them is not necessarily the same,
as the starting assumptions may be different.  

Then Sections 23-33 contain a kind of coda in which new issues
are presented; these sections are well worth study for themselves alone,
especially Section 33 which gives a meaning to Lagrange parameters
surely at odds with that intended by Lagrange himself.

The paper now follows in English translation.}

\section{General idea of translation and rotation of a solid system.}

I understand by a solid system any assemblage of points, either continuous
or discontinuous, that are mutually bound in a fixed way, such that if any three
of these points are located at positions not in a straight line, and all their
distances from other points of the system are given, the placement of the system 
will be completely determined for any placement of the triangle formed by these
three points.

Such a system can actually exist, since on a given triangular base, with given 
lengths of its sides, one can construct only a single pyramid identical to - that is, 
superposable on - another given pyramid.  {\it A pyramid obtained by reflecting
the first pyramid through its base would not be considered superposable.}
Accordingly, three noncollinear points of a solid being fixed, no displacement 
of the solid is possible.

But if only two points of the system are required to stay unmoved,
by fixing the distances of all other points from these two, one assures,
to begin with, that all points in a line with these two are unmovable. 
This line becomes a {\bf fixed axis}, and any other point of the system can 
only move on the circumference of a circle concentric with and normal to
the axis.  Since all the points of the system are bound unchangeably to any
one of them and to the fixed axis, the {\bf rotation} of one implies the 
rotation of all, and the {\bf amplitude} of this rotation is the same for all the points 
of the system.

Any displacement of a {\bf solid} about two fixed points reduces,
therefore, to a {\bf rotation}, of equal amplitude and {\bf in the same sense}
for all the points of the system, about the axis formed by the two fixed points.

Here it must be remarked that any given rotation can be replaced by a rotation 
of the opposite sense, of an amplitude complementary
to that of the first rotation with respect to $360^{\circ}$.

Different rotations about the same axis result in a rotation equal to their sum; here
one must take care to assign contrary signs to the amplitudes of rotation that are effected
in opposite directions, but otherwise the order of succession of the rotations remains
arbitrary.

If the amplitude of the rotations is infinitely small, the arcs described by the points of the
solid located at a finite distance from the axis become indistinguishable from their chords,
while the latter are variously inclined according to the angles of the rays  drawn from the
axis to the points of the system.

But suppose that the axis, while still firmly bound to the system under consideration, is
infinitely far removed from it, and that the system undergoes an infinitely small rotation
about the axis, of an infinitesimal order {\bf reciprocal} to that of the distance of the
axis of the system.  The effect will be to make all those points describe {\bf equal} and
{\bf parallel} straight lines, so that the system will have simply undergone a
{\bf translation}, that is to say a displacement resulting in an equal transport of all
its points in a certain direction.

Thus any {\bf translation} of a system can be considered rigorously as a rotation of 
infinitesimal amplitude about an axis infinitely far off and normal to the direction of that translation.

It is no surprise, then, to find moreover that all the properties of translations are implied by those
of rotations, just as those of a straight line are implied by those of a large circle to which
the line is tangent.  We need not linger over this.

We shall complete this general exposition of the displacement about a fixed axis by the 
following theorem, which is evident from the figure and whose consequences will be
of use in what is to follow.

\subsection{Movement of the axis of rotation without changing its direction.}

The rotation of a solid system about a fixed axis can be replaced by an equal rotation 
about another parallel axis, followed by a translation of the system equal and parallel to the 
chord of the arc described by a point of the second axis about the first.  Or, what
comes to the same thing, to the chord of the arc that would be described by a point of the 
first axis about the second, except that the direction of the translation must be reversed.

{\it The author apparently considers this theorem as having been proved by the preceding
discussion, along with the figure referred to.  In fact, the theorem is an immediate corollary to
the theorem on couples to be discussed in {\bf 10.1}. }

\section{Displacements by translation.}

If two situations of the same solid are such that all the lines joining a point of the solid
in one situation to the corresponding point in the other are equal and parallel, the solid
may clearly arrive from the first to the second situation by sliding parallel to itself
along one of these lines.  The length and direction of the line will measure those of the 
translation of the system.

\subsection{General law of composition of successive translations.}

If  the system undergoes several consecutive translations, differing both in direction and in extent,
it is evident that all these translations sum up to a single unique translation, equal and parallel
to the line that joins a point of the first situation to its corresponding point in the second.
This line would close the polygon traced by the successive translations of the point in question, and its
length and direction, as is well known, depend only on those of the various other sides of the polygon
and are {\bf independent of their order of succession}.

By means of this law, one can reciprocally decompose any given translation into a succession of diverse 
translations, provided only that the sum of their projections on three perpendicular axes, or 
more generally on {\bf any} arbitrary axis, is equal to the projection of the given translation on the
same axis, or to the sum of its projections on the three axes.  {\it The projections are here being \lq\lq added"
in what we would call a \lq\lq vectorial" sense.}

One may call this law of composition the {\bf law of the polygon of translations}.

\section {On the displacement of a system about a fixed point.}

{\it What follows is a remarkably concise proof of a celebrated theorem due to Euler.\cite{Eulfix}}

Let two situations of the same system be given, sharing a point $O$ which remains fixed 
in passing from one situation to the other.  {\it It should be understood that the displacement 
can be expressed in some way as a succession of rotations, with no reflection.}  Consider
two arbitrary points $A,B$ in the first situation, different from $O$ and not collinear with it.
In the second situation the corresponding points are $A',B'$.  {\it Here and elsewhere I have
used letters to designate points, lines, etc. as an aid to the reader's comprehension.  In
the French original no letters are used, all geometrical entities being described purely
by words.}

The isosceles triangles $AOA'$, lying in plane $P_a$, and $BOB'$, lying in plane $P_b$, share a common vertex $O$.  Through $O$
pass a plane $N_a$ which is normal to $P_a$ and bisects the vertex angle $\theta_a$. {\it See Fig. 1A$_i$ and 1A$_{ii}$.  In Fig. 1A$_{ii}$
the points $A,A'$ are superimposed, as the line $AA'$  is normal to the plane of the diagram.  The line marked $x_A$ is
actually the common perpendicular to $AA'$ and $L$, a line defined below.}

Likewise a plane $N_b$ normal to $P_b$ and bisecting $\theta_b$. {\it See Fig. 1B$_i$ and 1B$_{ii}$. Again, $x_B$ denotes
the common perpendicular to $BB'$ and $L$.}

The intersection of $N_a$ and $N_b$ 
will be a line $L$ {\it (See Fig. 1C)} passing through $O$ and normal to both $P_a$ and $P_b$.  
{\it In Fig. 1C, it must
be understood, if we regard the line $L$ as \lq vertical', that the triangle $OAA'$ is not \lq\lq horizontal", although its
base $AA'$ is horizontal.  Likewise $OBB'$ is not horizontal, although its base $BB'$ is.  Also, although $AA'$ and $BB'$
are both horizontal, they lie in general in different horizontal planes and consequently do not intersect.}

Any point $S$ on $L$ {\it (being equidistant from $A$ and $A'$, as well as from $B$ and $B'$)} can be considered 
as the common summit of two {\bf identical} or {\bf superposable} pyramids $OSAB$ and $OSA'B'$ 
{\it see the second paragraph of Section 1}, having as bases the triangles $OAB$ and $OA'B'$; so that
the line $L$, being invariably bound to the displaced system, remains unchanged by the displacement.
Hence this displacement reduces to a rotation around the fixed axis $L$.  {\it ($L$ is prevented from
sliding along its length by the immobility of $O$. This short passage is the whole proof, different
from that given by Euler.)}

{\it Two pyramids are congruent if they agree in the lengths of all six edges.  The edges
$OA, OB, AB$ are equal respectively to $OA', OB', A'B'$  because the displacement
is rigid.  The equality $SA = SA'$ follows from the construction of the plane $N_a$,
and $SB = SB'$ from that of $N_b$.  The edge $OS$ is the same in both pyramids
because $S$ is defined as a point on $L$ a certain
distance from $O$.  The line $L$, however, has been located by a construction drawing on
both the initial and the final configurations.  But now the pyramid $OSAB$ makes it possible
to determine the position of $L$ in terms of the initial configuration alone; and alternatively
$OSA'B'$ by the final configuration alone.  The congruence of the two pyramids 
then ensures that the location of $S$ is the same (both in the body and in space)
in both configurations, so that $L$ is unmoved by the
displacement. (The author's primary emphasis
on pyramids was already evident in the second sentence of {\bf 1}.)}

In the singular case (it could be avoided if we wished) in which the planes $N_a$ and $N_b$ coincide in a single plane $N$, we see that the axis is simply the line
of intersection of the planes $Q,Q'$ containing the original triangles $OAB$ and $OA'B'$.
Any other line, should it lie in the plane $N$, forms with the sides of these triangles two
trihedral angles symmetrically related but not superposable. 

To sum up, any displacement of a system {\it (achievable by rotations)} about a fixed 
point reduces to a rotation about a fixed axis passing through this point.  Or, more
generally, it reduces to an equal rotation about a different fixed axis parallel to
the first one, but otherwise located where one will, provided that the rotation be
followed by a translation of the same extent and direction as the chord of
the arc described by the original fixed point $O$ under the rotation about the new axis;
but the sense  of the translation must be opposite to that of the chord.  (This generalization
follows from the theorem in Section 1 relating to the parallel transport of an axis of rotation.)              

\section {On the most general displacement of a solid system in space.}

Now let us consider any two situations whatsoever of the same solid, and seek the
simplest mode of displacement that can bring the solid from one situation to the other.
Select any point $O_1$ of the solid in the first situation, and imagine that for each point
$A_1$,  a straight line is drawn from $O_1$ that is equal and parallel to the line $O_2A_2$ in the second situation.  
Denote by $A_{12}$ the termination of the line thus drawn. {\it Thus $O_1A_{12}$ and $O_2A_2$ are opposite
sides of a parallelogram, whence the same is true of $A_{12}A_2$ and $O_1O_2$.} We have thus 
constructed an intermediate assemblage of points, $O_{12},A_{12},...$, forming
a solid entirely identical to the one under consideration, but lying in a situation intermediate
between the two given ones.  Since $O_{12}= O$ {\it (by considering the case $A_1=O_1$)}, 
the intermediate situation can be derived from the first situation (in view of the theorem 
of the preceding section) by means of a certain rotation through an angle $\theta$ 
about a fixed axis $L$ passing through $O_1$ {\it (this is indeed the fixed point theorem of Euler, 
proved by the author in {\bf 3})}, while the passage from the intermediate to the 
second situation requires only a translation whose extent and direction are the same 
as those of the line $O_1O_2$  {\it (on account of the parallelogram $O_1O_2A_2A_1$)}.

Moreover, we observe that in view of the theorem on parallel transport in Section 1, nothing
prevents us from supposing that this intermediate situation of the solid is reached 
from the first situation by a rotation through the same angle $\theta$ about an 
arbitrary line $L'$ parallel to $L$, followed by a translation equal to the chord of the arc
that would have been described by a point on $L'$ in making the rotation about $L$
that would take the first to the intermediate situation.  {\it The whole displacement 
taking each $A_1$ to $A_2$ has already been decomposed as $R+T$ where $R$
is the rotation through $\theta$ about $L$ and $T$ is the translation taking $O_1$ to $O_2$.  
By the theorem on parallel transport, $R$ can be decomposed as $R'+T_{LL'}$ where
$R'$ is the rotation through $\theta$ about $L'$ and $T_{LL'}$ is a translation normal
to $L$ and $L'$.}

But this translation and the following one taking the intermediate to the second
situation combine to make a single translation equal and parallel to the line joining any
point on $L'$ to its corresponding point in the second situation. {\it Thus the whole displacement is
$R'+T_{LL'}+T = R'+T'$ where
$T' = T_{LL'}+T$ is the translation that would take $O_1'$ to $O_2'$.} Now, for any origin
chosen instead of $O$, there is only one axis of rotation possible; therefore we have
completely demonstrated the following theorem, indisputably one of the most beautiful
in geometry, which deserves to be considered the fundamental basis of the geometric
laws of the movement {\it [of a rigid body in 3 dimensions]}.

\subsection{Fundamental theorem.}

However a solid has been transported from one place to another, the displacement can
always be considered as resulting from two consecutive displacements, a rotation 
and a translation. The rotation takes place about a fixed axis passing through an
arbitrarily chosen point in the initial situation and parallel to a certain direction.
This direction, as well as the amplitude and sense of the rotation, is invariably
determined by the initial and final situations.  The translation takes place parallel
to the line joining a point on the said axis to its corresponding point in the second
situation, and its length is the length of that line.

The order of these displacements can be reversed: the translation can precede the rotation, 
but the latter then takes place about an axis passing through the point in
the final situation that corresponds to the point that was taken to be the origin 
in the initial situation.  In addition, the direction of the axis of rotation and its amplitude
and sense are the same for all the points of the system, whether before or after the translation.

In this succession of displacements, let us observe that the line joining any point of the 
initial situation to its corresponding point in the final situation, that is the resultant line
{\bf really} traced by this point, forms the third side of a triangle, of which the first
side {\it (representing the effect of the rotation)} is variable for different points of the system 
but always normal to the axis of rotation, and the second side, constant for all the points 
of the system, measures the translation of the system.

The projection of this resultant line on the axis of rotation is therefore constant for all
the points of the solid.  This constancy can be achieved only relative to the direction
of the axis of rotation, for the projection on an arbitrary direction is the sum of two
projections, that of the chord of the arc of the rotation and that of the line traversed
in translation.  The first of these is variable, the second is constant.  Their sum cannot,
therefore by constant for any direction other than that which causes the first of these
projections to vanish, that, is, the direction of the axis of rotation.

All the points of a solid system displaced in an arbitrary manner are therefore {\bf equally} transported {\bf relatively} to the {\it (invariably determined)}
direction of the axis of rotation.

If this constant projection vanishes, the displacement reduces to a rotation
about a certain {\it (preferred)} axis without any translation.  The transverse 
location of this axis is easily found by considering a plane normal to the 
common direction of all the possible axes of rotation.  The straight line from
any initial point in this plane to its corresponding final point will be the base of
an isosceles triangle lying in the plane, of which the angle at the vertex will
equal the amplitude of the rotation, and the vertex itself will lie on the {\it (preferred)}
axis to be found.                                                                                                                  

In the more general case {\it in which the displacement is not a pure rotation}, 
this constant projection is the measure of the {\bf absolute} 
translation of the system, by which is meant the minimum translation  among
all those associated with the possible axes of rotation, variously located 
 {\it but all having a common direction}; it is none other than the translation
of those points of the system that are displaced parallel to {\it [this common direction]}.  
If any of these points is chosen as the origin $O$, the axis of rotation, which we shall
distinguish by calling it the {\bf central axis of the displacement}, is also the axis
of translation.   

Thus the displacement reduces, with respect to its central axis, to turning about this 
axis while {\bf sliding} parallel to its direction: a kind of movement that has been 
compared to that of a screw turning in its nut.  This is the simplest expression of the
fundamental theorem, in which the two displacements of rotation and translation
take place {\it [simultaneously and]} orthogonally.  {\it Here the author attributes the theorem 
to M. Chasles\cite{screw}.}

{\it This theorem may be viewed as the natural generalization of Euler's fixed point
theorem - a generalization in which the displacement is assumed to result from an 
arbitrary series of rotations (excluding reflections) about axes that are {\bf not} 
assumed to have a common point of intersection.}

\section {Locating the central axis.}

But we have now to find this central axis, that is to find those points of the system 
carried by the displacement along a line parallel to the common direction of the
possible axes of rotation.  Now, one may arrive at this by the following construction.

Let $A$ be an arbitrary point in the first situation, and $A'$ its corresponding point 
in the second.  Let $L$ be the line through $A$ in the common direction which is
the {\it (known)} direction of the central axis. {\it Let us think of this direction as \lq\lq vertical".}
Let $P$ be the {\it (vertical)} plane containing $L$  and $AA'$. {\it (See Fig. 2A.)} Within $P$, erect from $A$ a
{\it (horizontal)} perpendicular to $L$ which terminates
at the point $Z$ chosen so that $AZA'$ is a right angle.  {\it We may introduce Cartesian
coordinates $x,y,z$, with origin at $A$, and let $L$ be the $z$-axis and
$AZ$ the $x$-axis.  Then $P$ is the $(x,z)$ plane.}

Now let $R$ be the plane
containing $AZ$ and normal to $L$.  {\it  $R$ is the $(x,y)$ plane. See Fig. 2B.
Note that $Z$ and $A'$ are both in the $(x,z)$ plane P and that the line
$ZA'$ is perpendicular to the $x$-axis $AZ$; therefore $ZA'$ is vertical, parallel
to the $z$-axis $L$.  Hence $ZA'$ is normal to the horizontal plane $R$. In other words,
$Z$ is the foot of the perpendicular dropped from $A'$ to $R$.} 

Within $R$, construct an isosceles triangle 
having $AZ$ as base and vertex angle equal to the {\it (known)} amplitude $\theta$
of the rotation. The vertex $V$ of this triangle will lie on the {\bf central axis},
provided only that the isosceles triangle is placed, relatively to $AZ$, in the sense
of the rotation.   {\it On a given base in a given plane, there are two ways to erect an
isosceles triangle with a given vertex angle, related by reflection in the base. Only
one is right.}

For it is evident that this vertex $V$, turning through $\theta$ about $L$, will
reach a point $\bar{V}$ in the plane $R$, such that the chord from $V$ to
$\bar{V}$ is equal and parallel to the line from $Z$ to $A$  {\it because the
triangle $VA\bar{V}$ is congruent to the triangle $ZVA$}, and moreover
that the translation from  $\bar{V}$ to the image $V'$ of $V$ in the second
situation is equal and parallel to the line from $A$ to $A'$.  {\it The whole 
motion from $V$ to $V'$ is made by the rotation carrying $V$ to $\bar{V}$, 
followed by the translation from $\bar{V}$ to $V'$; whereas the motion
of $A$ to $A'$ consists entirely of the translation.  Since a translation
affects all points equally, it follows that $\bar{V}V'$ is equal and parallel
to $AA'$. See Fig. 2C.} Hence the
resultant motion from $V$ to $V'$ is equal and parallel to the resultant
motion from $Z$ to $A'$ {\it that is, $V\bar{V} + \bar{V}V' = ZA + AA'$}
which in turn is parallel to $L$ {\it as noted in previous paragraph but one}.
But this is just the condition satisfied by points on the central axis.

And reciprocally, taking $VV'$ as the axis of rotation, the point $A$ on
rotation through $\theta$ about $V$ will travel to $Z$, and then by
the translation parallel to $VV'$ will travel from $Z$ to $A'$, reaching
its given destination.

And this construction shows that when the axis of relative translation $AA'$
is normal to the axis of rotation $L$, the whole displacement reduces
to a simple rotation about the central axis, since then $Z$ coincides
with $A'$ and so $V$ and $V'$ are the same.

Also, if a displacement of the solid is such that all the points of
the solid remain in mutually parallel planes, the displacement
reduces to a rotation about some fixed axis normal to these planes. 

{\it Although the foregoing construction correctly locates the central axis,
the information required to carry it out as well as the reasoning to 
justify it are drawn from the whole of the preceding paragraphs of the
essay.  But if one studies the construction in this light, one is forced
to keep in mind much redundant material.  If, on the other hand,
we disregard the preceding paragraphs, we are in danger of drawing
false inferences.  For example, it looks as if this construction requires
only one pair $A,A'$ of corresponding points.  But that is not so: the
fixed axis theorem of Euler, proved in {\bf 3}, requires two such
pairs, and the result is a necessary part of the proof of the author's
construction to locate the central axis.  Therefore he gives another
construction, entirely self-contained and requiring only the location
of two pairs $A,A'$ and $B,B'$ as well as the direction of the central
axis.  He attributes this construction to \lq\lq mon ami M. L\'{e}vy".}

Drop from $A$ a line perpendicular to the central axis at a point $C$,
and likewise from $A'$ a perpendicular at $C'$, and consider the
quadrilateral $ACC'A'$.  Let $A^m,C^m$ be the midpoints of $AA'$, $CC'$;
then the line $A^mC^m$ will be perpendicular to both $AA'$ and $CC'$.
{\it This follows from the symmetry of the figure with respect to a $180^{\circ}$
rotation about $A^mC^m$.} This property gives the following construction of
the central axis, being given only the points $A,A'$, another 
corresponding pair $B,B'$, and the direction of the central axis:

Through each point $A_m$ (resp. $B_m$), pass a line $L_a$ (resp. $L_b$)
parallel to the central axis, as well as a line $N_a$ (resp. $N_b$)
normal to both $AA'$ and $L_a$ (resp. $BB'$ and $L_b$). Let $P_a$ (resp. $P_b$)
be the plane formed by $L_a$ and $N_a$ (resp. $L_b$ and $N_b$). Then the planes
$P_a$ and $P_b$ will intersect precisely on the central axis.

{\it To justify the construction, we must refine our notation to take into
account the use of two pairs: the quadrilateral based on $A$ must now
be called $AC_aC'_aA'$, and that based on $B$ must be called
$BC_bC'_bB'$.  The points $C_a,C_b$ are not necessarily the same,
but both lie on the central axis, as well as $C'_a,C'_b$.

Since $L_a$ is parallel to $C_aC'_a$, $N_a$ is normal to $C_aC'_a$ as
well as to $AA'$; thus $N_a$ is the very line of which $A^mC_a^m$ is
a segment.  Therefore $C_a^m$  lies on $N_a$.  But then the plane
$P_a$ can equally well be described as formed by $AA'$ and $C_aC'_a$
instead of by $AA'$ and $L_a$, so that it contains the whole line $C_aC'_a$.
that is the whole central axis.  Likewise the plane $P_b$ is formed by
$BB'$ and $C_bC'_b$ and also contains the central axis.  Hence 
$P_a$ and $P_b$ intersect on the central axis, q.e.d.

There are degenerate cases, which will be left to the reader's study.}

\section {Consequences of the Fundamental Theorem.}

It now behooves us to set forth the principal corollaries that follow from the
fundamental theorem, relating to the particular displacements of the points,
the lines, and the planes of a solid system.  

(a) The distances separating each point of the solid in the first situation
from its corresponding point in the second all have equal projections on the direction 
of the central axis; this common projection is the measure of the {\bf absolute}
translation of the system.

(b) As any line belonging to the displaced system does no more, 
{\bf with respect to its direction}, than turn about the axis of rotation, there 
results a very simple relation between the angle formed by this line with
the axis of rotation and that formed between the initial and final directions of this line,
to wit:

\lq\lq The sine of the half-angle of displacement of any line belonging 
to a displaced system is equal to the sine of the half-rotation of the system, 
multiplied by the sine of the angle between this line and the axis of rotation."

{\it The author here is concerned only with the {\bf directions} of the line in
the initial and final configurations.  Now, the set of all possible directions can be
mapped onto the points on the surface of a sphere (say of radius $r$)
by mapping each direction $D$ to the unique point $P$ for which the line
from $O$, the center of the sphere, to $P$ has the direction $D$.  The initial
direction of the line in question thus maps to a point $A$ on the surface
of the sphere, and its final direction to a point $B$.  The direction of the
central axis (the \lq\lq axis of rotation") is mapped to a point $C$, which
we shall call the North Pole so as to make use of the ideas of latitude
and longitude.  The central axis is then the diameter from North to South
Pole, and the rotation about this axis clearly does not change the angle
made with it by the line in question. Thus $A$ and $B$ lie on the same
circle of latitude, which we shall call $L$, and the arc from any point $P$
on $L$ to the North Pole is of a fixed length, making a constant angle $POC$
which we shall call $\chi$.  The \lq\lq angle of rotation of the system" is
the difference in longitude between $A$ and $B$, which we may call $\phi$,
and the \lq\lq angle of displacement" of the line in question is the angle
$\theta = AOB$, which measures the geodesic distance from $A$ to $B$.

With these definitions, the proposition rendered above in words may be
expressed in symbols as
\begin{equation}\label{geo}
\sin(\theta/2) = \sin(\phi/2)\sin\chi.
\end{equation}

\noindent We can understand the r\^{o}le of $\chi$ by noting that 
when $\chi$ becomes small, the geodesic distance from $A$ to $B$ 
becomes small because of the pinching of the base of the isosceles triangle $ABC$,
even though the difference $\phi$ in longitude is kept constant.

The equation \eqref{geo} can be derived as a restriction of the
spherical Law of Cosines to isosceles triangles; but the author
gives no indication of having such reasoning in mind.  Instead, he
says briefly that the proposition \lq\lq is made evident by observing 
the figure".  Unfortunately, the figures originally appearing in the
Journal de Math\'{e}matiques have been lost, but I believe I
have closely reconstructed the one referred to in this passage,
with the aid of my associate Dr. Familton.

The easy demonstration of \eqref{geo} depends on the construction
of planes and straight lines in the interior of the sphere (see Fig. 3A), particularly
the straight line $\bar{AB}$.  On the one hand, this line is a chord
of the great circle $G$, of radius $r$ and center $O$ (see Fig. 3B) upon which
both $A$ and $B$ are located; since this chord subtends an angle
$\theta$ at $O$, we have
\begin{equation}\label{th}
\bar{AB} = 2r\sin(\theta/2).
\end{equation}

\noindent On the other hand, this same line is a chord of the
{\bf small} circle of latitude, $L$ (see Fig. 3C), whose center $K$ lies
on the axis $OC$ and whose radius we shall call $\rho$. Since
the angle subtended at $K$ is $\phi$, we have
\begin{equation}\label{ph}
\bar{AB} = 2\rho\sin(\phi/2).
\end{equation}

\noindent Comparing \eqref{ph} with \eqref{th}, we find
\begin{equation}\label{thph}
\frac{\sin(\theta/2)}{\sin(\phi/2)} = \rho/r.
\end{equation}

\noindent Finally we determine $\rho/r$ by
considering the right triangle $OKP$ (see Fig. 3D) for an
arbitrary point $P$ on $L$. The hypotenuse $OP$
has length $r$, the angle $POK$ is $\chi$, and
the side $KP$ opposite this angle has length
$\rho$.  Therefore 
\begin{equation} \label{rhrch}
\rho/r = \sin\chi.
\end{equation}

\noindent Substituting \eqref{rhrch} into \eqref{thph},
we obtain \eqref{geo}.}

(c) Any line parallel to the axis of rotation is transported parallel to itself,
while any line normal to that axis suffers an angular displacement equal
to the amplitude of the rotation. {\it (Special cases of (b)}]]

(d) Any plane invariably bound to the displaced system, 
and normal to the axis of rotation, is therefore transported
into a plane parallel to the initial one, at a distance equal to the 
absolute translation of the system.

(e) The midpoint of the line that joins any point of the system
to its correspondant {\it i. e., to its final position after displacement}
is the point of that line that approaches most closely to the central axis
of the displacement.

(f) the midpoints of all the lines that join the various points of a plane figure
to their correspondants after an arbitrary displacement lie in a single plane,
which also contains the midpoints of the lines joining any point {\bf outside} 
the plane figure to its {\bf symmetric} correspondant.

This plane makes equal angles with the planes of the two plane figures, 
as well as with the {\bf corresponding} lines bound to the two figures,
but if not within their respective planes, then {\bf symmetrically} inclined.

{\it The meaning of \lq\lq symmetrically" in (f) is \lq\lq by reflection in
the special plane under consideration."  I leave the study of (f) to 
the most ambitious of readers.}

\section{The decomposition of any displacement into two pure rotations.}

Having presented the fundamental geometric law of the passage of a solid 
from one given situation to another, also given in an arbitrary way, we have
now to study the law of composition of successive displacements; by 
means of this law one can construct or calculate the elements of the 
composite displacement, that is the position of its central axis, the 
amplitude of its rotation, and the extent of its translation.

We have already presented the law of composition of translations; we
shall next give that of rotations about different fixed axes, and finally
that of arbitrary displacements, each resulting from a combined 
translation and rotation.

From this law of composition of rotations about different axes, we shall
deduce an important transformation of the fundamental theorem
{\it ({\bf 4.1})}, to wit:

\lq\lq Any displacement of a solid system can be represented, in an 
infinite number of ways, as the composition of two successive rotations 
of the system about two nonintersecting fixed axes.  The product
of the sines of half these rotations, multiplied by the sine of the
angle between the two axes and by the minimum distance between
them, is equal - for each of these {\bf conjugate pairs} of axes -
to the product of the sine of half the angle of rotation of the
system about the {\bf central} axis of the {\it (total)} displacement with
half the length of its absolute translation." 

{\it Let the arbitrary (total) displacement be characterized by a central
axis $\cal C$, a rotation angle $\theta_C$, and a translation distance
$T$.  It is asserted that there are infinitely many \lq\lq conjugate pairs",
each of which consists of a rotation $\theta_A$ about an axis
$\cal A$ followed by a rotation $\theta_B$ about an axis
$\cal B$, such that the composition of these two rotations
is equivalent to the total displacement under consideration. 
Each of these conjugate pairs is related to the total 
displacement by the equation
\begin{equation}\label{comp}
D\sin\nu\sin(\theta_A/2) \sin(\theta_B/2) = (1/2)T\sin(\theta_C/2),
\end{equation}

\noindent where $D$ is the minimal distance between the two axes
and $\nu$ is the angle between their directions.

The author only states this remarkable theorem here, deferring its proof
to {\bf 11}.}

To put it another way, the volume of the tetrahedron of which two 
opposite edges lie anywhere along the respective conjugate axes,
provided that the length of each of these edges is proportional to
the sine of the corresponding half-angle of rotation, is the same for
all conjugate pairs of rotations whose composition is equivalent to a 
given displacement.

{\it The equivalence of this second statement of the theorem to the first
is based on a theorem of geometry, that the volume of any
tetrahedron is one-sixth the product of the lengths of any two opposite edges,
times the minimal distance between them, times the sine of the angle 
between them.  This formula can be established by a variety of methods;
further study is left to the reader.}

Thus, any displacement of a solid system reduces to a rotation about 
one or two fixed axes.

In the case where one of these axes is parallel to the central axis 
{\it (of the total displacement)}, it follows from the law of composition
of rotations that its conjugate is situated at infinity and that the
rotation corresponding to it becomes infinitely small and therefore
amounts to a simple translation.  This leads to the original version
{\it ({\bf 4.1})} of the fundamental theorem, so that that version is
no more than a particular case of the theorem we have just given.

\section{The composition of two given rotations.}

We have now to present the law of composition of successive 
rotations of a solid about different axes.  {\it The author now
begins the succession of {\bf composition} theorems that will 
lead ultimately to the {\bf decomposition} theorem stated
in the previous section.}

Let us begin by considering only two intersecting axes, and
let us seek to determine the {\bf resultant} axis of these 
two rotations - the one about which the given solid will be
finally found to have turned, in order to arrive from the
initial to the final situation.  {\it (A problem closely related
to this was studied in \cite{Eulcomp}.)}

This resultant axis must be placed in such a way that in
being subjected to the two rotations indicated about the 
supposed intersecting axes, it comes back to its initial
position.  If, therefore, through each of the given axes
one passes a plane that makes an angle with the plane
of the two axes that is equal to half the rotation about
that same axis, the intersection of these two planes
will be the resultant we seek, as it arrives by virtue
of the first rotation at the position symmetric by reflection
in the plane of the two axes, and returns by the second
rotation to its original position.

{\it Call the two intersecting axes ${\cal A}$ and ${\cal B}$; they
determine a plane that we shall call $\cal AB$.  The resultant axis
we call ${\cal C}$, and let $\cal AC, BC$ denote the planes formed
respectively by ${\cal A}$ and ${\cal C}$, and by ${\cal B}$ and ${\cal C}$. The
angles of rotation about ${\cal A}$ and ${\cal B}$ are $\theta_A$
and $\theta_B$. 

It is assumed that the axes ${\cal A}$ and ${\cal B}$, intersecting
at the origin $0$, as well as the amplitudes of rotation $\theta_A$ 
and $\theta_B$, are given, and the problem posed is to determine 
the axis  ${\cal C}$ and rotation amplitude $\theta_C$ such that the 
rotation $\theta_A$ about ${\cal A}$
followed by the rotation $\theta_B$ about ${\cal B}$ will produce as resultant
the rotation $\theta_C$ about ${\cal C}$.  By Euler's theorem ({\bf 3}) 
such a ${\cal C}$ exists and passes through the origin $0$.

In understanding the solution proposed by the author, it is
necessary to distinguish the angle between two {\bf axes}
(e. g. between ${\cal C}$ and ${\cal B}$) from
the angle between two {\bf planes} (e. g. between
$\cal AC$ and $\cal AB$).  

The solution proposed is that
$\cal AC$ should make an angle $\theta_A/2$ with $\cal AB$,
and $\cal BC$ should make an angle $\theta_B/2$ with $\cal AB$.
In this way the planes $\cal AC$ and $\cal BC$ are determined,
and ${\cal C}$ is their intersection.  

The argument is that if ${\cal C}^*$ is defined as the reflection
of ${\cal C}$ in $\cal AB$, then the angle between $\cal AC$ and $\cal AC^*$
will be twice that between $\cal AC$ and $\cal AB$, that is twice $\theta_A/2$,
so that the first rotation (of $\theta_A$ about ${\cal A}$) will 
bring $\cal AC$ to $\cal AC^*$ and hence ${\cal C}$ to ${\cal C}^*$
(since the angle between ${\cal C}$ and ${\cal A}$ is unchanged
by the rotation).  Then by a similar argument, the second
rotation (of $\theta_B$ about ${\cal B}$) will take ${\cal C}^*$
back to ${\cal C}$.  So the combined effect of the two partial
rotations will be to leave ${\cal C}$ unaltered.  But this is what
is required in order that ${\cal C}$ be the resultant axis.

(It may well be objected that it has not been shown that the axis
${\cal C}$ is restored to its original position {\bf in the same sense}
in which it began.  But if not, the situation may be remedied by taking
the appropriate sign of $\theta_C$.)}

 At the same time one sees that the angle between the two planes
{\it ($\cal AC$ and $\cal BC$)} will be half the angle of the resultant rotation {\it ($\theta_C$)}.  For the first axis, 
which does not move under the first rotation, is displaced only
by the second, and it describes about the resultant axis, 
determined as shown above, an angle twice that between 
the two planes.

{\it The phrase \lq\lq determined as shown above" refers to
the determination of $\cal C$ as the intersection of two planes, of
which in particular the plane $(\cal BC)$ forms an angle $\theta_B/2$ with
the plane $(\cal AB)$. It follows that a rotation about $\cal B$ through
an angle $\theta_B$ brings the plane $(\cal AB)$ to its reflection in $(\cal BC)$.
and in particular it brings ${\cal A}$ to its reflection ${\cal A}^*$ in $(\cal BC)$.

But since  ${\cal A}$ was unchanged by the first rotation about itself,
this second rotation about $\cal B$ must have the same effect on $\cal A$
as does the composite rotation through the angle $\theta_C$ (to be found)
about $\cal C.$ That is, this composite rotation must take $\cal A$ 
to ${\cal A}^*$. But this requires that $\theta_C$ be twice the angle
between $(\cal AC)$ and $(\cal BC)$, as asserted by the author.}

%

Here let us note that the half-rotation {\it (of a plane)} about each axis can
be measured equally well by the interior or the exterior angle 
of the two planes passing through that axis, only the sense
of the rotation depending on which measure one adopts, 
since any rotation {\it (of a point)} about an axis in one sense is equivalent
to a rotation in the opposite sense with an amplitude
complementary to the first by 360$^{\circ}$.

{\it The conscientious
reader may wish to ascertain that if a consistent sign convention
be followed, whereby the angle of rotation about an axis is
measured either always clockwise or always counterclockwise
looking along the direction of the axis, this construction
will yield the correct sign of $\theta_C$ in relation to those of
$\theta_A$ and $\theta_B$.}

Furthermore, as to the order of succession of these rotations,
it comes about that if the two rotations are supposed to take
place in a certain order, leading to a particular resultant
axis ${\cal C}$, then by reflecting this axis in the plane of the of the
two given axes one obtains the resultant ${\cal C}^*$ of the same two
rotations in the inverse order.  From this we see that the amplitude
of the resultant rotation is independent of the order of the
two given ones, but that the position of the resultant axis
depends on this order, and that in the composition of
more than two rotations about arbitrary intersecting axes,
the order cannot be modified without altering {\bf both}
the position of the resultant axis and the amplitude of
the resultant rotation.

{\it The last statement is a bit too strong: the resultant amplitude
will be unaltered if the order is {\bf completely} inverted, as 
from 1234 to 4321, no matter how long the succession is. 
Likewise the amplitude is preserved under a cyclic permutation
as from 12345 to 34512.  A corollary of these two facts is that
the amplitude cannot be altered by any permutation of the
composing displacements unless they number at least 4.}

Such is the characteristic difference between the composition
of rotations and that of successive translations.  In fact, these
two kinds of composition are analogous in a way similar
to the properties of a plane triangle and those of a spherical triangle.
For if one compares the translations parallel to the three sides of
a planar triangle to the sines of the half-angles of rotations effected
around the three sides of a trihedral angle, the values of these
translations and those of these sines are in equal manner proportional
to the sines of the angles opposite to the respective sides of the 
planar triangle and of the trihedral angle.

{\it The author is essentially comparing the Law of Sines for a planar
triangle to that for a spherical triangle.  But he adds the complication 
of associating the sides of the triangle to the corresponding
translations or rotations of a solid system.}

\section {Composition of infinitesimal rotations.}

But these two resultant axes ${\cal C}$ and ${\cal C}^*$, corresponding to the same two rotations 
in two different orders of succession, will {\bf coincide} in the plane 
of the two axes if the rotations become infinitely small, and from this
there follow two important consequences:

First, the order of succession of infinitesimal rotations about two 
intersecting axes (and, as it follows, about as many such axes as one wishes)
is immaterial.  And second, the axis and amplitude of the infinitesimal
rotation resulting from the succession of two infinitesimal rotations 
${\cal A}$ and ${\cal B}$ about two intersecting axes are determined 
in the same way as the axis and translation length that would result from 
two successive translations proportional to the given rotations and parallel to their axes.

{\it The author is referring here to the theorem enunciated in Section 1:

\lq\lq Thus any {\bf translation} of a system can be considered rigorously as a rotation of 
infinitesimal amplitude about an axis infinitely far off and normal to the direction of that translation."

\noindent Note that this theorem concerns a {\bf single} axis, and in
the present context it applies {\bf separately} to ${\cal A}$ and to ${\cal B}$. 
Note also that the angle between ${\cal A}$ and ${\cal B}$ is finite, unlike
that between ${\cal C}$ and ${\cal C}^*$.}

Since, by removing the axes of rotation far away, one may transform 
the infinitesimal rotations into finite translations perpendicular to these
axes and inclined one to the other just as the axes are to each other,
one achieves all the generality of the law of composition of finite rotations,
which by mediation of the infinitesimal rotations includes also the law
of composition of translations.

{\it The author evidently does not mean that the two axes are removed 
far from each other, since they continue to intersect.  Rather, he is observing
the action of the infinitesimal rotations at points far from both of the two axes.
But these points are regarded, for the present purpose, as \lq\lq here", while
the two axes with their intersection are \lq\lq there", that is removed to infinity.

The author's point is that the law of composition of finite rotations is so powerful
that by suitable applications of it one can derive that of composition of finite 
translations as well.}

\section {Of the composition of rotations about two parallel axes.}

All the points of the system displaced by two consecutive rotations about
two parallel axes remain within parallel planes normal to these axes.  Therefore
the displacement reduces to a simple rotation about a certain axis 
parallel to the first two.  This being admitted, the mode of determination
and of construction of the resultant axis of two intersecting axes applies 
equally well to this case, and yields a resultant axis parallel to the first
two and a composite rotation equal to the sum or the difference of the 
given rotations, according as they act in the same or opposite senses.

\subsection{Couples of parallel rotations.}

But here we encounter a remarkable case, that in which the two rotations
are equal and of opposite sense.  The composite rotation is then null and 
the composite axis is placed at infinity, which causes the displacement
to amount to a simple translation.  To be precise, each point of the
displaced solid has traversed a line of length and direction constant 
for all points of the system; the direction of this line is normal to the
two axes but makes an angle with the normal to the plane of the two 
axes, equal to half the given angle of rotation about each axis; and its length
is the product of the  distance between the axes by twice the sine of
half the rotation angle.

{\it Using rectilinear coordinates $x,y,z$, let the two axes point in the $z$-direction
and their separation $d$ in the $x$-direction, so that they lie in the $x$-$z$ plane;
then the normal to their plane is in the $y$-direction.  Since the coordinate $z$ 
is unchanged in each rotation, we may regard the whole operation as confined
to the $x$-$y$ plane.  Take any initial point $P$ in this plane.  Under the rotation
about the first axis ${\cal A}$ it describes an arc of amplitude $\theta$ ending at a point $P^{\circ}$.  
Then the rotation about the second axis ${\cal B}$ carries the point along an arc also
of amplitude $\theta$ but in the opposite sense, from $P^{\circ}$ to its final 
position $P'$.  The straight line from $P$ to $P'$ is asserted by the author 
to have the same length and direction for all points $P$ - a fact by no means
obvious.  The following proof can best be followed by consulting Figs. 4a and 4b.

Let the axes ${\cal A}$ and ${\cal B}$ cut the $x$-$y$ plane at the points $A$ and $B$,
so that the distance $d$ between the axes is the line segment $AB$. And let the 
distances from $P^{\circ}$ to the axes $A$ and $B$ be respectively
$r_A$ and $r_B$.   Then the isosceles triangles $PAP^{\circ}$ and $P^{\circ}BP'$
are similar, having the same vertex angle $\theta$; and their scale ratio is $r_A/r_B$.  
Therefore $PP^{\circ}/P^{\circ}P' = r_A/r_B$.

The base angles of the isosceles triangles pertaining to various initial points $P$ 
are all equal to $\phi=(180^{\circ}-\theta)2$. Therefore the angles $AP^{\circ}B$
and $PP^{\circ}P'$ are equal, being both equal to $\phi+\chi$ where $\chi$ 
is the (undetermined) angle $AP^{\circ}P'$.  (This assertion merits close examination
of Fig. 4a.)  

Hence the triangles $AP^{\circ}B$
and $PP^{\circ}P'$ are similar, having two sides in the same ratio 
$AP^{\circ}/BP^{\circ} = PP^{\circ}/P'P^{\circ} = r_A/r_B$ and
agreeing in the included angle.

It follows that $PP'/AB = PP^{\circ}/AP^{\circ}$
or  $PP' = (AB)(PP^{\circ}/AP^{\circ})=AB(2\sin(\theta/2))$, 
as asserted by the author. Note that this equation does not
involve the angle $\chi$.

Furthermore, the angle between the directions $PP'$ and $AB$ is the angle
through which the triangle $PP^{\circ}P'$ must be turned about $P^{\circ}$ so as 
to make the angles $AP^{\circ}B$ and $PP^{\circ}P'$ coincide.  Clearly this 
angle is $\phi$, so that $PP'$ makes an angle 
$\phi = (180^{\circ}-\theta)/2 = 90^{\circ} - \theta/2$ with the $x$-axis.  Therefore 
it makes an angle $\theta/2$ with the $y$-axis, as claimed.

Figs. 4a and 4b show two realizations of this construction. Using the line segment
AB as the reference for length and direction, one sees that the point $P$ is placed
differently in the two diagrams, and that the ratio $r_A/r_B$ is also quite different,
as well as the size of the angle $\chi$.
Nevertheless, the dotted line from $P$ to $P'$ has the same length as well as
the same direction in both diagrams, and in each diagram the triangle $PP^{\circ}P'$ 
is similar to the triangle $AP^{\circ}B$ as found in the above proof.

It may help the understanding to distinguish between {\bf determinate} ratios, which 
for a given $\theta$ are independent of the choice of $P$, and {\bf indeterminate} ratios,
which are affected by that choice.  The two isosceles triangles in the diagram have 
determinate shape, which yields the ratio $2\sin(\theta/2)$ (base to sides), but their
relative size is indeterminate, depending on the ratio $r_A/r_B$. For the two 
symmetric triangles, however, the situation is reversed:  the shape of these triangles
involves $r_A/r_B$ and is hence indeterminate, while the ratio of their size
is determinate as the diagram exhibits it as tied to the shape of the isosceles
triangles.

Once it is established that the direction and length of $PP'$ is the same
for all $P$, the author's statements about this direction and length can be derived 
very easily by setting $P = A$. For then one has obviously $P^{\circ} = P$, and
$PP' = P^{\circ}P'$ which is the base of an isosceles triangle
with vertex angle $\theta$ at $B$ and side $AB$. The author's statements follow.}

The order of succession of the two rotations makes a difference; if the order
is reversed, the composite line of translation is reflected about the normal to the two axes.

All this follows easily from a comparison of similar triangles {\it (as shown above)}; and then
if the rotations are infinitesimal, the order of succession becomes immaterial, and
the translation acts along the normal to the plane of the two axes.

Thus any {\bf couple of parallel rotations} {\it (not necessarily infinitesimal)} is equivalent
to a simple translation, and reciprocally, any translation can be replaced in an infinite
number of ways by a couple of this kind. {\it The word \lq\lq couple" is meant to imply
that the two rotations are equal and opposite.}

These couples of parallel rotations compose and decompose, in accordance with the
law of translations, in an arbitrary order of succession, acting in all the positions
that correspond in length and direction to a particular translation; compositions
and decompositions which can be found by substituting for the couples 
the translations that they represent. {\it The order of different {\bf couples} 
is arbitrary; the order of the two rotations forming a single couple is not, if
the rotations are finite.}

Thus we have generalized to couples of finite rotations the law of composition
which M. Poinsot, I believe, was the first to state for couples of infinitesimal rotations.

\section {Proof of the general decomposition theorem} 

As any displacement of a solid system can be reduced to a rotation 
followed by a translation {\it (see {\bf 4})}, and this translation can always be
replaced by a {\bf couple of rotations} {\it (see {\bf 10.1})}, one
of whose axes intersects the given rotation axis of the system, 
and the rotations about these two intersecting axes can be
composed {\it (see {\bf 8})} into a single rotation, there results
immediately the proof of the transformation stated above
{\it (see {\bf 7})} of the fundamental theorem {\it (see {\bf 4.1})}, to wit: that
any displacement of a solid system can be accomplished
in an infinite number of ways by the succession of two
rotations about two fixed nonintersecting axes.

{\it The axes will be nonintersecting unless the translation
is null, in which case one rotation suffices. This one-sentence
proof can benefit, as usual, by some expansion of the
reasoning and naming of the geometrical entities.

Let the axis of the given rotation be called $\cal C$ and
the accompanying translation be called $T$.  As shown
in {\bf 10.1}, $T$ can be replaced by a couple of
equal but opposite rotations about two parallel axes
which we may call ${\cal B}'$ and $\cal B$, of which
${\cal B}'$ may be located so as to intersect $\cal C$.
(The axis ${\cal B}'$ is the one that was called $\cal A$
in {\bf 10.1}; here we shall define an entirely different
axis as $\cal A$.) 

Then, by the method of {\bf 8}, the
rotations about $\cal C$ and ${\cal B}'$ can be composed
to make one about an axis we shall call ${\cal A}$; the two rotations about the
nonintersecting axes $\cal A$ and $\cal B$, performed
successively, generate the same displacement  as
$\cal C$ and $T$.

This is quite clear as far as it goes, but it gives no clue as how
to derive eq. \eqref{comp} of {\bf 7}.  This will be done in {\bf 13};
it is unnecessary to carry out the demonstration also in the present
context, as it would involve the same steps sometimes done
in reverse order.}

\section {Rotations about an arbitrary number of fixed nonintersecting axes.}

Finally, there is the composition of rotations about an arbitrary number of
fixed nonintersecting axes.  Let us take a point in the space, upon which 
we shall study the effect of all these rotations in their order.  We have seen
that any rotation about a fixed axis can be replaced by another equal rotation,
accomplished about another axis parallel to the first, followed by a translation
equal to the chord of the arc described by a point of the new axis about the 
first in consequence of the rotation given at the outset.  We have also seen 
that a translation {\bf followed} by a rotation about an axis passing through
the endpoint of the axis of translation can instead be {\bf preceded} by it,
if the axis of rotation passes through the {\bf origin} of the axis of translation.

{\it This last statement deserves careful examination.  Let a translation $T$
be followed by a rotation $R$ whose axis is $L$.  Decompose $T$ into
$T_l$ parallel to $L$ and $T_t$ transverse to $L$.  Since $T_l$ commutes with $R$, 
it suffices to consider $T_t$ and to project the whole situation onto a transverse plane.

We may represent points on this plane by complex numbers, and for simplicity
let us represent the displacement $T_t$ by the number $1$.  Take $R$ to be
a rotation of the plane through an angle $\theta$ about the point $1$, and 
identify the author's phrase \lq\lq axis of translation" as the line from
$0$ to $1$.  Thus the endpoint of the axis of translation is at $1$, and its
origin is at $0$.  The claim is that the effect of $T_t$ followed by $R$ is the
same as that of $R'$ followed by $T-t$, where $R'$ is the rotation through $\theta$
about $0$.

But this is easily proven.  Start with a point $z$ and first perform $R'$; this takes
$z$ to $ze^{i\theta}$.  Then perform $T_t$; the result is $ze^{i\theta}+1$.  On the 
other hand, $T_t$ acting first on $z$ produces $z+1$, and the radius vector 
from $1$ to $z+1$ is $z$ so that the second transformation $R$ replaces 
the term $z$ by $ze^{i\theta}$, yielding the final result $ze^{i\theta}+1$
as before.}

This being given, if, through the point of origin we have chosen for study,
we pass axes parallel to each of the given nonintersecting axes, the
displacement of the system operates successively about these axes,
by means of the transport of the rotations to the intersecting axes 
respectively parallel to the original ones.  By virtue of the successive
replacement of rotations about axes passing through the endpoints
of the translation lines by rotations about axes passing through their origins,
the displacement of the system will be partitioned into a series 
of rotations respectively equal to those originally given, taking place successively 
about {\bf intersecting} axes parallel to the first series, followed
by a series of translations resulting from the chords successively
traversed by the chosen point about the original nonintersecting axes
{\bf in the order ascribed to the rotations}.

{\it One may ask whether the final composite translation would be the same
if a different point of origin \lq\lq $p$" had been chosen for study.
The answer is yes, because a change in $p$ can be simulated by keeping
$p$ unchanged and rigidly changing the positions of all the rotation axes.}

(The composition of rotations about {\bf intersecting} axes and that of translations
will take place in the manner described above; in this case the composite
displacement will reduce to a rotation and a translation whose axes both pass 
through the point of intersection.)

We see from {\it (the general)} construction that the elements of the final composite
rotation depend only on the amplitude and direction of the individual rotations, 
and are not changed by any parallel movement of their axes; while the length
and direction of the composite translation depend, as well, on the
{\bf positions} of the individual rotations, seeing that the chords successively
described by the chosen point vary in length and direction according to
the successive positions that the displaced point takes relatively to the
various given axes.

If, in the system of these axes, there are found consecutive pairs that form
{\bf couples of parallel rotations}, it is evident that these couples do
not contribute anything to the determination of the direction and amplitude
of the resultant rotation, and that they influence only the length and direction 
of the resultant translation, as the point whose successive rotations determine 
this translation will be found, upon completion of each couple, to have
described the translation equivalent to the couple.

\section {The case of only two nonintersecting fixed axes.}

{\it The author intends in this section to give a second proof of the \lq\lq two-axis" theorem
of {\bf 7}, which has already been proved in {\bf 8-11}.  But in that first proof, he started 
with a displacement described in screw form with the central axis
${\cal C}$ and rotation angle $\theta_C$ given, and {\bf constructed} the two
(usually) nonintersecting axes ${\cal A},{\cal B}$ of rotations $\theta_A,\theta_B$
whose composition gives this displacement.  In this second proof he assumes that
${\cal A},{\cal B}$ are given, along with $\theta_A,\theta_B$, and constructs
$\cal C$ and $\theta_C$.  Thus the {\bf decomposition} theorem of {\bf 7} is replaced
in the present section by a {\bf composition} theorem, a generalization of the
theorem of {\bf 8} in that the axes ${\cal A}$ and ${\cal B}$ are no longer 
required to intersect.  And this time the author carries through the derivation
of eq. \eqref{comp} in {\bf 7}, which was omitted in {\bf 11}.}


Consider two nonintersecting axes and their shortest distance $D$, and take 
as the origin of the displacement the end $A$ of that shortest distance lying 
on the first axis of rotation ${\cal A}$ - that is, the rotation to be first executed.  Then,
on passing through that origin an axis ${\cal B}'$ parallel to the second one ${\cal B}$ given, 
the two intersecting axes will be composed into a third, which will be
the {\bf axis of rotation}  of the displacement relative to that origin $A$.

{\it Let the symbol $D$ refer equally to the shortest distance between ${\cal A}$ 
and ${\cal B}$ and to the line segment of that length, pointing from its intersection $A$ with ${\cal A}$ 
to its intersection $B$ with ${\cal B}$. The author has introduced a substitute axis ${\cal B}'$, 
parallel to ${\cal B}$ but passing through $A$. Regarding $A$ as Euler's ({\bf 3}) fixed center, 
the fixed-point theorem says that the rotations $\theta_A , \theta_{B'}$ about the intersecting axes 
${\cal A}$, ${\cal B}'$ can be composed to make a rotation $\theta_{C'}$ about a third axis 
${\cal C}'$, also passing through $A$. The author calls ${\cal C}'$ the \lq\lq axis of rotation of the displacement relative to A."

At the same time (see {\bf 1.1} and {\bf 4}), the rotation about $\cal B$ can be accomplished by performing
the rotation about ${\cal B}'$ followed by a translation, which we may call $T'$.  Thus
the entire displacement $({\cal A},\theta_A)$ followed by $({\cal B},\theta_B)$ is equivalent to
the three actions $({\cal A},\theta_A)$, $({\cal B}',\theta_{B'})$, $T'$ taken in sequence,
which in turn is the same as $({\cal C}',\theta_{C'})$ followed by $T'$.  More compactly, we can
write this result as
\begin{equation}\label{ABeqCprTpr}
({\cal A},\theta_A)({\cal B},\theta_B) = ({\cal C}',\theta_{C'})T'.
\end{equation}
}

The accompanying {\bf axis of translation} $T'$ will
be given {\it ({\bf 1.1})} as the chord of the arc described by this same origin $A$ in consequence
of the rotation about the second {\it (${\cal B}$)} of the two axes given.

{\it It is desired, however, to express the whole displacement in
terms of the {\bf central axis} $\cal C$ and its associated quantities; thus we must have
\begin{equation}\label{ABeqCT}
({\cal A},\theta_A)({\cal B},\theta_B) = ({\cal C},\theta_C)T,
\end{equation}

\noindent where $T$, called the {\bf absolute translation} of the displacement,
is directed along the axis $\cal C$. 

Comparing \eqref{ABeqCT} to \eqref{ABeqCprTpr}, we find that the rotation $({\cal C},\theta_C)$
differs from the rotation $({\cal C}',\theta_{C'})$ only by a translation, which may be written $T'-T$
since translations compose by addition.  It follows immediately that $\cal C$
and ${\cal C}'$ are parallel and that $\theta_C = \theta_{C'}$.  Furthermore,
both these rotations move points only within planes perpendicular to $\cal C$ and to ${\cal C}'$.
Therefore the difference $T'-T$ is perpendicular to $\cal C$, hence perpendicular to $T$
which lies along $\cal C$.  It follows that $T$ is the projection of $T'$ onto $T$,
as stated in the author's next remark:}

%

The projection of this chord onto the composed axis of rotation, determined 
as above, measures the absolute translation of the {\it (composite)} displacement.
It is equal to the sum of the projections of the two sides of the isosceles triangle
of which it is the base.  The two sides are equal {\it (in length)} to the shortest 
distance between the two given nonintersecting axes.

{\it Let the whole displacement under study carry the point $A$ to its final
position $A_f$.  Since the first rotation about $\cal A$ left $A$ unmoved,
the movement from $A$ to $A_f$ is accomplished entirely by the second
rotation, about $\cal B$.  But the point $B$, lying on $\cal B$, is unmoved
by this rotation; since the displacement is rigid, the distances $AB$ and
$A_fB$ are equal. That is, the triangle $ABA_f$ is isosceles, with vertex
at $B$ and sides $AB$, $A_fB$ both equal to $D$ as defined above.

The base of this triangle is the chord $AA_f$, which has previously been identified
(author's comment after eq.\eqref{ABeqCprTpr}) as giving the translation $T'$.
Its projection on $\cal C$ is equal to the sum of the projections of
$AB$ and $BA_f$.  (This is readily understood by thinking of $AB, BA_f, AA_f$ 
as vectors.)}

It is easily
demonstrated in addition that the two sides are equally inclined to the 
composed axis. 

{\it This is an important claim. We see at once that the projection of $AB$
on $\cal C$ is $D\cos S$ where $S$ is the angle between $\cal C$ and $D$.
The author wishes to establish that $BA_f$ makes the same angle with $\cal C$
and therefore has the same projection.  This will establish the important
equation $T = 2D\cos S$.}

%


In fact, this shortest distance {\it ($D$)} is normal to the plane ${\cal AB}'$ of the two 
intersecting axes.   Now, in 
considering the angle formed by this normal $D$ with the {\bf reflection} {\it ${\cal C}^*$} of
    the resultant axis $\cal C$ in this same plane ${\cal AB}'$, one sees that this angle 
does not change when one supposes it mobile and displaced rigidly by the second
rotation {\it (about $\cal B$)}, which brings the reflected line $\cal C^*$of which we speak into coincidence
with the resultant axis $\cal C$.

But, in this rotation, the normal $D$ is rotated through a plane
perpendicular to the second axis $\cal B$ {\it (sweeping out a cone with vertex at $B$)} 
so as to become parallel to the second
side {\it $A_fB$} of the isosceles triangle we are considering, and as it is evident
that the angle of the normal $D$ with the resultant axis $\cal C$ is supplementary to
that which it forms with the reflection $\cal C^*$ of that axis, one sees that the
resultant axis $\cal C$ is, as we have just said, equally inclined with respect to the 
two sides $AB$, $BA_f$ of this isosceles triangle. {\it The projection of $\cal C^*$ 
onto the initial position of $D = AB$ is the negative of the projection of $\cal C$;
that is, it equals the projection of $\cal C$ onto $BA$.  Therefore the $\cal B$
rotation sweeps the angle under discussion to the angle between $\cal C$
and $BA_f$ as asserted.

The preceding argument is best understood by comparing it with the
construction in {\bf 8}, where $\cal A$ and $\cal B$ intersect;
$\cal B$ and ${\cal B}'$ coincide, $D$ vanishes but its direction
is still defined as the normal to the plane $AB$; $A_f$ is identical
to $A$ since both rotations leave $A$ fixed; and there is no
translation $T'$ or $T$.  The key specification is that $\theta_A$
is twice the angle between the planes $\cal AC$ and $\cal AB$, and 
$\theta_B$ is twice the angle between the planes $\cal BC$ and $\cal AB$, 
so that $\cal C$ is reflected in $\cal AB$ by the $A$-rotation and
reverse-reflected by the $B$-rotation. The author's determination
of $\theta_C$ depends on computing the final position of a 
line whose initial position was $\cal A$, deduced on the one hand from the effect
of the $B$-rotation and on the other hand from that of the $C$-rotation.

In the present case the plane of reflection is taken as ${\cal AB}'$, but
it could be any plane parallel to both $\cal A$ and $\cal B$, without
changing the angles.  Instead of $\cal A$, the initial position of
the moving line is taken as $D$, and again a comparison is
made between two ways of finding its final position.  The required
agreement between the two ways, as in  {\bf 8}, yields the result,
which in this case is the equality of the angles made by $AB$ and
$BA_f$ with the resultant axis $\cal C$.}

From which one finally concludes that the absolute translation of
a solid system arising from the succession of two rotations about
two fixed nonintersecting axes is equal to double the distance 
between the two axes, projected on the direction of the composite,
or resultant, axis. {\it That is, if the side $AB$ of the isosceles triangle
$ABA_f$ is extended to twice its length, the projection of this doubled
side onto the direction of $\cal C$ will fall on $A_f$. The resulting
equation is
\begin{equation}\label{TDS}
T = 2D \cos S
\end{equation}

\noindent as anticipated above.

%
%
%

 }

But it is evident that the cosine of the angle $S$ of this distance $D$ with
the composite axis $\cal C$ is equal to the sine of the angle made by
$\cal C$ with the plane of the two composing axes $\cal A,B'$ 
{\it (since $D$ itself is normal to this plane)}. which
is found to be equal to the product
of the sines of the given half-rotations by the sine of the angle
between the two axes, divided by that of the composite
half-rotation.  {\it That is, 
\begin{equation}\label{SAB/C}
\cos S = \sin h = \sin (\theta_A/2) \sin(\theta_B/2) \sin\nu /\sin (\theta_C/2),
\end{equation}

\noindent where $h = \angle(\cal C, \cal AB)$ and $\nu = \angle(\cal A, \cal B)$.}

Equation \eqref{SAB/C} follows from the law of proportion of the 
sines of the half-rotations to those of the angles included
between the opposing axes.

{\it Here again the author has compressed many steps into
one.  We can understand \eqref{SAB/C} more readily by noting
that $D$ and $T$ no longer appear, so that the equation 
involves only {\bf directions}.  This enables us to associate
the directions $\cal A,\cal B,\cal C$ with points $A,B,C$ on a 
sphere, forming the vertices of a spherical triangle. Since we
are dealing only with directions, we need not distinguish between
$C$ and $C'$ or between $B$ and $B'$. The angle
$\nu = \angle(\cal A, \cal B)$ is just the arc $c$, the side
of the triangle opposite to $C$. Thus \eqref{SAB/C} becomes
\begin{equation}\label{cAB/C}
\sin h = \sin(\theta_A/2)\sin(\theta_B/2)\sin c/\sin (\theta_C/2).
\end{equation}

\noindent where $h$ is what we
might call the {\bf altitude} of the triangle, that is the arc
running from $C$ to $c$ and making a right angle with the latter.

The author asserts that \eqref{cAB/C} is a consequence of the
\lq\lq law of proportion..." which can be stated as
\begin{equation}\label{Rodprop}
\frac{\sin(\theta_A/2)}{\sin a} = \frac{\sin(\theta_B/2)}{\sin b} = \frac{\sin(\theta_c/2)}{\sin c}.
\end{equation}

\noindent We must ask how this law is arrived at, and also how it leads to  \eqref{cAB/C}.

Let us adopt the custom, with spherical triangles, of allowing the letters $A,B,C,$
to stand also for the spherical angles at the respective vertices, thus
$B = \angle(\cal AB,\cal BC)$, etc. Now referring to the law of composition
in {\bf 8}, we see that $\angle(\cal AB,\cal BC)$ is set equal to half the rotation angle $\theta_B$,
and likewise for the other axes.  Therefore \eqref{Rodprop} is equivalent to

\begin{equation}\label{sinprop}
\frac{\sin A}{\sin a} = \frac{\sin B}{\sin b} = \frac{\sin C}{\sin c},
\end{equation}

\noindent which is just the Law of Sines for a spherical triangle.
(I am unable to see how the author could have arrived at \eqref{Rodprop}
without relying on this law from spherical trigonometry.)

Now we must show how \eqref{Rodprop} leads to \eqref{cAB/C}.
If $h$ were the altitude of a planar triangle we would obviously have 
$h=a\sin B$. The spherical analogue is $\sin h = \sin a\sin B$, 
also an application of the Law of Sines to the special case of
a right triangle.  But since $B = \theta_B/2$ from {\bf 8}, we have
\begin{equation}\label{hgood} 
\sin h = \sin a \sin(\theta_B/2) = \sin(\theta_A/2)\sin(\theta_B/2)\frac{\sin a}{\sin(\theta_A/2)}.
\end{equation}

\noindent Now, applying \eqref{Rodprop}, we may replace $\frac{\sin a}{\sin(\theta_A/2)}$
by $\frac{\sin c}{\sin(\theta_C/2)}$, obtaining \eqref{cAB/C} as asserted by the author.}

From which we arrive at the modified Fundamental Theorem in the form ({\bf 7}) 
in which we have already pronounced it, namely:

\lq\lq any displacement of a solid system can always arise,
in an infinite number of ways, from the succession of two 
rotations about two nonintersecting fixed axes, {\bf provided
that} the product of the sines of the successive 
half-rotations by the distance between the two conjugate axes
and by the sine of the angle of these axes is equal to the product
of the absolute half-translation of the displaced system by 
the sine of the resultant half-rotation."

{\it That is, $\sin(\theta_A/2)\sin(\theta_B/2) D \sin\nu = (1/2) T \sin(\theta_C/2)$,
where $D$ is the distance and $\nu$ the angle between the two axes.
(Having established \eqref{cAB/C}, we can discard the spherical
triangle representation, writing once more $\nu$ in place of $c$,
and substitute \eqref{SAB/C} into \eqref{TDS}.
This gives indeed
\begin{equation}\label{Rodbig}
2D \sin(\theta_A/2) \sin(\theta_B/2) \sin\nu = T\sin(\theta_C/2)
\end{equation}

\noindent in agreement with eq. (6) of {\bf 7}.)}

\section {Composition of general displacements.}

We are now in a position to resolve completely the following general problem,
in which one considers the succession of {\it (an arbitrary number of)}
displacements of the same solid.

Being given the axes of rotation and translation as well as the amplitude
of the rotations and extent of the translations for each successive 
displacement of a system, it is required that we construct the axes
[[and amplitudes]] of rotation and translation of this system relative
to a given origin.

The solution of this problem is evidently the same as that of the
previous one, where it was only a matter of rotations about fixed axes,
since the translations can be replaced by couples of rotations
about fixed axes.  We have briefly indicated the solution in the
last paragraph of  {\bf 12}.  We therefore need not
linger over it further.

\subsection{The particular case of infinitely small displacements.}

The solution is considerably simplified when one considers only
infinitely small displacements.  First of all, the order of the 
rotations is indifferent, and their composition by whatever
number around intersecting axes operates like that of
translations proportional to these rotations and parallel 
to these axes. Second, the order of the rotations and 
translations successively accomplished by the origin
of the displacement is equally indifferent, and each of
these rotations and translations can be established
directly and separately, as though the point to be 
displaced were displaced only alternately and not
successively, which follows from the fact that the space
traversed by each of these displacements is infinitely small.
The composition of these partial translations resulting
from withdrawing from the given axes, or from the translations 
themselves that are joined to the rotations, acts 
in accordance with the same law as that of the rotations.

\section {}

We have now to apply {\bf calculation} to the geometric laws that we have just presented 
concerning the general displacements of a solid system.  We shall start by deriving the
formulas for change of coordinates of points in the solid system, which hold such a 
large place in analytical mechanics.

Let $x,y,z$ and $x+\Delta x,y+\Delta y,z+\Delta z$ be the coordinates of two points 
of which the first is moved to the second by the displacement of the system,
and let $\xi,\eta,\zeta$ be the coordinates of the midpoint of the line joining the two, so that
\begin{equation}\label{xietazetaorig}
\xi = x +(1/2)\Delta x, \;\;\;\;\; \eta = y +(1/2)\Delta y,\;\;\;\;\;\zeta = z +(1/2)\Delta z.
\end{equation}

\noindent Furthermore, let $g,h,l$ be the angles formed by the direction of the axis of rotation 
with the three coordinate axes, $\theta$ the amplitude of the rotation, $t$ the 
absolute size of the translation, and $X,Y,Z$ the coordinates of an arbitrary point 
on the central axis of the displacement.  {\it In much of what follows, the {\bf origin of coordinates}
may be assumed to lie on the central axis; that is, we can take $X = Y = Z = 0$, or in the vector
notation to be introduced, $\vec{W}=0$.} 

Consider the right triangle whose hypotenuse is formed by the line joining the initial and final point 
and whose sides are given, one by the arc of the chord described by the initial point 
under the rotation $\theta$, and the other by the line traversed in a translation by this same point 
after undergoing the rotation. Clearly, the changes $\Delta x,\Delta y,\Delta z$ are respectively 
equal to the projections of this hypotenuse, that is to the sum of the projections 
of the other two sides of this triangle on the respective coordinate axes. 

Now, the side equal and parallel to the absolute translation $t$ gives the three projections 
$\cos g, \cos h, \cos l$; the other side is equal to $2u\tan(\theta/2)$, $u$ denoting the distance 
from the central axis to this same side {\it (the one formed by the chord)}.   Let us call $G,H,L$ 
the angles between this side and the coordinate axes. Then we have immediately
\begin{eqnarray}\label{Deldeforig}
\Delta x = t\cos g + 2u\tan(\theta/2) \cos G   \nonumber \\
\Delta y = t\cos h + 2u\tan(\theta/2) \cos H  \nonumber \\
\Delta z = t\cos l + 2u\tan(\theta/2) \cos L.
\end{eqnarray}

{\it This and the following equations can be better understood if translated into 
modern vector notation. Let us define a right-handed orthonormal system $\hat{t},\hat{u},\hat{v}$
where $\hat{t}$ points along the central axis; $\hat{v}$ along the chord; and 
$\hat{u}$, perpendicular to both, points to the midpoint of the chord from the base of the perpendicular 
dropped from that midpoint to the central axis.  The whole displacement 
$(\Delta x, \Delta y, \Delta z)$ may be designated as $\vec{\Delta}$; then the above equations say that
\begin{equation}\label{Deldef}
\vec{\Delta} = \vec{t}+ 2\hat{v}u\tan(\theta/2),
\end{equation}

\noindent where $\vec{t} = \hat{t} t$ is the translation vector.}

Since one has necessarily
\begin{equation}\label{tvorth}
\cos g \cos G + \cos h \cos H + \cos l \cos L = 0,
\end{equation}
\noindent {\it (this says that $\hat{t}\cdot\hat{v} = 0$)}
one deduces
\begin{equation}\label{Deldtt}
\Delta x\cos g + \Delta y\cos h + \Delta z\cos l = t
\end{equation} 

\begin{equation}\label{Delsq}
(\Delta x)^2 + (\Delta y)^2+ (\Delta z)^2 = t^2 + 4u^2\tan^2(\theta/2).
\end{equation} 

\noindent {\it That is, $\vec{\Delta}\cdot\hat{t} = t$ and $\vec{\Delta}\cdot\vec{\Delta} = t^2 + (2u\tan(\theta/2))^2$.}

The first terms {\it (of \eqref{Deldeforig})} $t\cos g, t\cos h, t\cos l$ represent the part of
the changes that arise from the absolute translation displacement; the second 
{\it (set of three)} terms, the part due to the rotation performed by the displacement. 
In comparing these first terms to the second, one finds that the first, which measure
the effect or {\bf moment} of the translation of the system, have for value the 
projections of this translation on each axis of the coordinates, while the second,
which represent for each point the effect or {\bf moment} of the rotation of
the system, have for value the projection upon each coordinate axis of {\it (the area of
an isosceles)} triangle whose vertex is the midpoint of the line finally traversed
by the point considered, and whose base is a line {\it [segment]} directed along
the central axis and of length $4\tan(\theta/2)$.

In the case of an infinitely small displacement, this midpoint {\it (the vertex of the triangle)}
coincides with the initial point, and consequently the moment of the rotation, relative to any given direction,
is equal to double the projection on that direction of {\it (the area of)} a triangle whose
vertex is the point {\it (under consideration)} and whose base taken on the central axis
is equal to the rotation of the system.

This explains how the theory of projections applies to the laws of translation through
linear projections, and to those of rotation through the projection of areas. {\it A possible
influence of Grassmann here? or independent? Let us
continue.

\noindent (Wordy explanation here omitted: terms in $t$ represent the effect of translation, 
those in $u$ the effect of rotation.)}

The chord $2u\hat{v}\tan(\theta/2)$ being normal both to the central axis
and to the perpendicular dropped from the point $\xi,\eta,\zeta$ to this axis, which has length $u$, 
we have
\begin{eqnarray}\label{vomWtorigL}
u\cos G = (\eta - Y)\cos l - (\zeta - Z)\cos h \nonumber \\
u\cos H = (\zeta - Z)\cos g - (\xi - X)\cos l \nonumber \\
u\cos L = (\xi - X)\cos h - (\eta - Y)\cos g. 
\end{eqnarray}

{\it Here we need to define more vectors.  Let $\vec{W}= (X,Y,Z)$ represent the point 
on the central axis that has been identified as locating that axis 
with respect to the immovable space within which the solid exists.  
Let $\vec{w}$ stand for the base of the perpendicular dropped from the midpoint of the chord of 
rotation to the central axis.  ($\vec{W}$, although chosen arbitrarily along the axis, is fixed
for a particular displacement of the solid, whereas $\vec{w}$ slides along the axis as we consider 
the trajectories of different points of the solid.) And let $\vec{\omega} = (\xi,\eta,\zeta)$ 
represent the position of the midpoint of the chord.  Then \eqref{vomWtorigL} tells us that
\begin{equation}\label{vomWtL}
\hat{v}u = (\vec{\omega} - \vec{W})\times\hat{t}.
\end{equation}

I note here that all the author's equations starting with \eqref{vomWtorigL} are 
consistent with $(\hat{t},\hat{u},\hat{v})$ being a left-handed system, whereas the 
motion described in the Chasles theorem at the end of {\bf 4}
(\lq\lq mouvement ... de la vis dans son \'{e}crou", movement of a screw in its nut)
makes it definitely right-handed.  I choose to write vector equations in the original
right-handed notation.  When equation arrays are written out in the original
in terms of the components, I shall reproduce them without change; but in writing
these arrays in vector notation I shall reverse the order of all cross-products.  
Thus, I shall interpret \eqref{vomWtorigL} as
\begin{equation}\label{vtomWR}
\hat{v}u = \hat{t}\times(\vec{\omega} - \vec{W})
\end{equation}

\noindent instead of as \eqref{vomWtL}.  I shall do this consistently
without further comment.

Moreover, $\hat{t}\times(\vec{w}-\vec{W}) = 0$ 
since $\vec{w}-\vec{W}$ lies along the central axis.  Therefore
$\vec{\omega} - \vec{W}$ may be replaced by $\vec{\omega} - \vec{w}$
which is just $\vec{u} = \hat{u}u$.  So \eqref{vomWtorigL} becomes
\begin{equation}\label{vtu}
\hat{v}u = \hat{t}\times\vec{u}
\end{equation}

\noindent which need not surprise us.}

And in consequence,
\begin{eqnarray}\label{DelGam+omLqorig}
\Delta x = A + p\eta - n\zeta, \nonumber \\
\Delta y = B + m\zeta - p\xi, \nonumber \\
\Delta z = C + n\xi - m\eta,
\end{eqnarray}

\noindent $A,B,C,m,n,p$ being six constants that depend on the
position of the central axis, the length of the translation, and
the amplitude of the rotation, as follows:
\begin{eqnarray}\label{Gamdeforig}
A = t\cos g + 2\tan(\theta/2)(Z\cos h - Y\cos l), \nonumber \\
B = t\cos h + 2\tan(\theta/2)(X\cos l - Z\cos g), \nonumber \\
C = t\cos l + 2\tan(\theta/2)(Y\cos g - X\cos h),
\end{eqnarray}

\noindent and
\begin{eqnarray}\label{qdeforig}
m = 2\tan(\theta/2) \cos g \nonumber \\
n = 2\tan(\theta/2) \cos h \nonumber \\
p = 2\tan(\theta/2) \cos l .
\end{eqnarray}

{\it By including \lq\lq the position of the central axis" as a variable,
the author signals that the following calculations do not assume
that this axis passes through the origin.  Indeed, much of what
follows in this Section becomes trivial if that assumption (the 
\lq\lq null central axis" assumption or NCA) is made.  For example,
\eqref{Gamdeforig} becomes $A = t\cos g$, etc.

\noindent Let us put $\vec{\Gamma} = (A,B,C)$, $\vec{q} = (m,n,p)$;
then these definitions become
\begin{equation}\label{Gamdef}
\vec{\Gamma} = \vec{t} + 2\vec{W}\tan(\theta/2)\times\hat{t}
\end{equation}

\noindent and
\begin{equation}\label{qdef}
\vec{q} = 2\hat{t}\tan(\theta/2);
\end{equation}

\noindent under NCA, \eqref{Gamdef} becomes $\vec{\Gamma} = \vec{t}$.

Going back to \eqref{DelGam+omLqorig}, in vector notation it reduces to
\begin{equation}\label{DelGam_qRom}
\vec{\Delta} =\vec{\Gamma} + \vec{q}\times\vec{\omega}
\end{equation}

\noindent ($\vec{\Delta} = \vec{\Gamma}$ under NCA).

Applying \eqref{Gamdef} and \eqref{qdef}, this becomes
\begin{eqnarray}\label{Delt+2tRomWq}
\vec{\Delta} & = & \vec{t} - 2\hat{t}\times \vec{W}\tan(\theta/2) + 2\hat{t}\tan(\theta/2)\times\vec{\omega} \nonumber \\
& = & \vec{t} + 2\hat{t}\times(\vec{\omega} - \vec{W})\tan(\theta/2),
\end{eqnarray}

\noindent agreeing with \eqref{Deldef} in view of \eqref{vtomWR}.}

If we denote by $\alpha,\beta,\gamma$ \lq\lq les variations des coordonn\'{e}es de l'origine des axes coordonn\'{e}s"
{\it (the coordinates of the point to which the origin of coordinates is carried by the displacement)}, we have
the following relation:
\begin{eqnarray}\label{deldeforig}
\alpha = A + (1/2)(p\beta - n\gamma)\nonumber \\
\beta = B + (1/2)(m\gamma - p\alpha)\nonumber \\
\gamma = C + (1/2)(n\alpha -m\beta).
\end{eqnarray} 

{\it Let us introduce the vector $\vec{\delta} = (\alpha,\beta,\gamma)$.Then \eqref{deldeforig} becomes
\begin{equation}\label{deldef}
\vec{\delta} =\vec{\Gamma} +(1/2)\vec{q}\times\vec{\delta}.
\end{equation}

\noindent To arrive at this equation, consider \eqref{DelGam_qRom} and recall that $\omega = \vec{r} + (1/2)\vec{\Delta}$,
where the displacement carries $\vec{r}$ into $\vec{r} + \vec{\Delta}$. Now take the special case $\vec{r} = 0$ (the
origin of coordinates).  In this case $\vec{\Delta}$ takes the value of $\vec{\delta}$, by definition of the latter.  Thus
\eqref{DelGam_qRom} reduces to $\vec{\delta} = \vec{\Gamma} +(1/2)\vec{q}\times\vec{\delta}$ which is exactly
\eqref{deldef}.  (Under NCA, one has simply $\vec{\delta} = \vec{\Gamma}$.)

Without NCA, equation \eqref{deldef} looks indeterminate since $\vec{\delta}$ is defined in terms of itself.  
But actually the system \eqref{deldeforig} provides three linear equations in the three unknowns 
$(\alpha,\beta,\gamma)$ which are perfectly determinate.

The solution is obtained transparently by vector algebra.  Let $\vec{U}$ be the 
position vector of the point on the central axis closest to the origin, and 
$\hat{V} = \hat{t}\times\hat{U}$, then $\hat{t},\hat{U},\hat{V}$ are orthonormal.
Moreover, \eqref{Gamdef} reduces to
\begin{equation}\label{Gamt+URt}
\vec{\Gamma} = \vec{t} + 2\vec{U}\times\hat{t}\tan(\theta/2) 
\end{equation}

\noindent since $\vec{W}$ and $\vec{U}$ are both on the central axis.  This can be written
\begin{equation}\label{GamtV}
\vec{\Gamma} = \vec{t} - 2U\hat{V}\tan(\theta/2) = \vec{t} - Uq\hat{V}
\end{equation} 

\noindent in view of \eqref{qdef}.  From this we have
\begin{equation}\label{GamUorth}
\vec{\Gamma}\cdot\vec{U} = 0
\end{equation}

\noindent and also
\begin{equation}\label{Gamt}
\vec{\Gamma}\cdot\hat{t} = t.
\end{equation}

\noindent Then from \eqref{deldef} we find
\begin{equation}\label{delt}
\vec{\delta}\cdot\hat{t}= t,
\end{equation} 

\noindent and also, in view of \eqref{GamUorth},
\begin{equation}\label{delUdelUqV}
\vec{\delta}\cdot\vec{U}= -(1/2)\vec{\delta}\times\vec{q}\cdot\vec{U}    
= -(1/2)\vec{\delta}\cdot \vec{q}\times\vec{U} =-(1/2)\vec{\delta}\cdot Uq\hat{V} . 
\end{equation}

\noindent which can be rewritten as
\begin{equation}\label{delUdelV}
\vec{\delta}\cdot\hat{U } = -\vec{\delta}\cdot\hat{V}\tan(\theta/2).
\end{equation}

\noindent On the other hand, \eqref{GamtV} also gives
\begin{equation}\label{GamUq}
\vec{\Gamma}\cdot\hat{V}= -Uq
\end{equation}

\noindent whence
\begin{equation}\label{delVUdelU}
\vec{\delta}\cdot\hat{V} =-Uq+ \vec{\delta}\cdot\hat{U}\tan(\theta/2).
\end{equation}

Substituting \eqref{delUdelV} into \eqref{delVUdelU} gives
\begin{equation}\label{delVUdelV}
\vec{\delta}\cdot\hat{V} = -Uq - \vec{\delta}\cdot\hat{V}\tan^2(\theta/2),
\end{equation}

\noindent or
\begin{equation}\label{delVUsin}
\vec{\delta}\cdot\hat{V} = -Uq\cos^2(\theta/2) = -2U\cos^2(\theta/2)\tan(\theta/2) = - U\sin\theta
\end{equation}

\noindent whence by \eqref{delUdelV}
\begin{equation}\label{delUUcos}
\vec{\delta}\cdot\hat{U} = +2U\sin(\theta/2)\cos(\theta/2)\tan(\theta/2) = U(1-\cos\theta).
\end{equation}

Combining \eqref{delt},\eqref{delVUsin}, and \eqref{delUUcos},we obtain
the  formula
\begin{equation}\label{deltUV}
\vec{\delta} = \vec{t} + \vec{U} - U[\hat{U}\cos\theta + \hat{V}\sin\theta], 
\end{equation}

\noindent which tells us that the endpoint of $\vec{\delta}$ can be located
by passing through the origin $0$ a circle in the $\hat{U},\hat{V}$ plane with
center at $\vec{U}$ on the central axis, moving counterclockwise on this circle
from $0$ through an arc subtending an angle $\theta$ at the center,
and erecting on the endpoint of this arc the translation vector $\vec{t}$. Indeed, this is 
just the operation (rotation $\theta$ about the central axis followed by translation $\vec{t}$)
that takes the origin of coordinates into its image under the displacement considered.}

By means of \eqref{deldeforig} one may replace, if one wishes, the constants 
$A,B,C$ by their values in terms of $\alpha,\beta,\gamma$.  This leads to
\begin{eqnarray}\label{Deldelorig}
\Delta x = \alpha + 2\tan(\theta/2)(\eta-\beta/2) \cos l - (\zeta - \gamma/2) \cos h)  \nonumber \\
\Delta y = \beta + 2\tan(\theta/2)(\zeta-\gamma/2) \cos g - (\xi - \alpha/2) \cos l)  \nonumber \\
\Delta z = \gamma + 2\tan(\theta/2)(\xi-\alpha/2) \cos h - (\eta - \beta/2) \cos g) 
\end{eqnarray}

\noindent where the first terms $\alpha,\beta,\gamma$ express the translation
{\bf relative to the origin of coordinates}  (the length being $\sqrt{\alpha^2+\beta^2+\gamma^2}$) 
and those pertaining to the rotation express the rotation {\bf about an
axis passing through the origin}.

{\it To understand the last remark, let us call the direction of the central axis \lq\lq vertical"
and recall the statement of the \lq\lq Fundamental Theorem" in {\bf 4}.  There it is
pointed out that {\bf any} vertical axis can be chosen as the rotation axis, and that
the accompanying translation, while constant for all points considered, is not vertical
unless the central axis is the one chosen.  In general, it is only the {\bf projection}
of this translation on the axis of rotation that is vertical. There is an additional horizontal
translation that compensates for the change of rotation axis.

The equation \eqref{Deldelorig} can be written in vector notation as
\begin{equation}\label{Deldelmodq}
\vec{\Delta} = \vec{\delta} + \vec{q}\times(\vec{\omega} - (1/2)\vec{\delta}),
\end{equation}

\noindent where $\vec{\delta}$ is given by \eqref{deltUV}.  In \eqref{deltUV} the term
$\vec{t}$ can be taken as a vertical translation and the remaining term as
the trajectory of the initial point $\vec{r} = 0$ under rotation about the central
axis. But now the author desires to consider the vertical axis through the origin
(call it the 0-axis) as the axis of rotation.  Then the initial point $\vec{r} = 0$
does not move under the rotation, and hence the entire expression \eqref{deltUV}
must be viewed as translation, containing both a vertical and a horizontal part.

For a general point $\vec{r}$, still taking the 0-axis as the axis of rotation, the
whole of $\vec{\delta}$ is still translation and therefore the second term of
\eqref{Deldelmodq} gives the rotation about the 0-axis.

Of course, under NCA there is no change of \eqref{deltUV} and one still has $\vec{\Delta} = \vec{\delta} = \vec{t}$.}

We could have established these formulas directly, as well as those which precede,
and which we have constructed on the central axis.

\section {Equations of the central axis.}

The equations of the central axis follow in their turn from the above formulas, with the greatest
simplicity; for the rotation has no effect on any point lying on that axis, and one then has
$\Delta x = t\cos g,\;\;\;\Delta y = t\cos h,\;\;\;\Delta z = t\cos l$.  {\it (That is, $\vec{\Delta} =\vec{t}$.)}
The coordinates $x,y,z$ describe some point on that axis  {\it (I shall denote this point by $\vec{r}$)}
and taking this into account one has the equations sought, expressed by means of $\alpha,\beta,\gamma$:
\begin{eqnarray}\label{taloneorig}
py-nz+A = t\cos g = \frac{m(Am+Bn +Cp)}{m^2+n^2+p^2} = \frac{m(\alpha m+\beta n +\gamma p)}{m^2+n^2+p^2}, \nonumber \\
mz-px+B = t\cos h = \frac{n(Am+Bn +Cp)}{m^2+n^2+p^2} = \frac{n(\alpha m+\beta n +\gamma p)}{m^2+n^2+p^2}, \nonumber \\
nx-my+C = t\cos l = \frac{p(Am+Bn +Cp)}{m^2+n^2+p^2} = \frac{p(\alpha m+\beta n +\gamma p)}{m^2+n^2+p^2}.
\end{eqnarray}

\noindent {\it In vector notation this array becomes
\begin{equation}\label{talone} 
\vec{\Gamma} + \vec{q}\times\vec{r} = \vec{t} = \vec{q}\,\vec{\Gamma}\cdot\vec{q}/q^2 = \vec{q}\,\vec{\delta}\cdot\vec{q}/q^2, 
\end{equation}

\noindent where the first equality is simply a rearrangement of \eqref{Gamdef} taking account
that both $\vec{r}$ and $\vec{W}$ are arbitrary points on the central axis, and also using \eqref{qdef}; 
the second equality follows from \eqref{qdef} which gives
$\vec{q}\,\vec{\Gamma}\cdot\vec{q} = q^2 \hat{t}\,\vec{\Gamma}\cdot\hat{t}$; and the third from \eqref{delt} with \eqref{Gamt}.}

This can be simplified by eliminating the constants $A,B,C$; we find {\it (using \eqref{deldeforig})}
\begin{equation}\label{alphbetgamcentralaxorig}
\frac{x-(1/2)\alpha - \frac{p\beta-n\gamma}{m^2+n^2 + p^2}}{m} = \frac{y-(1/2)\beta - \frac{m\gamma-p\alpha}{m^2+n^2 + p^2}}{n} 
 = \frac{z-(1/2)\gamma- \frac{n\alpha-m\beta}{m^2+n^2 + p^2}}{p}.
 \end{equation} 

\noindent {\it To derive this, we combine \eqref{talone} with \eqref{deldef}, and eliminating $\vec{\Gamma}$,we have
\begin{equation}\label{deircentralax}
\vec{\delta}-(1/2)\vec{q}\times\vec{\delta} + \vec{q}\times\vec{r}= \vec{q}\,\vec{\delta}\cdot\vec(q)/q^2
\end{equation}

\noindent which can be rearranged, using the identity
$(\vec{\delta}\times\vec{q})\times\vec{q} = \vec{q}\vec{\delta}\cdot\vec{q}- q^2\vec{\delta}$, to yield
\begin{equation}\label{r-r0}
(\vec{r}-\vec{r_0})\times\vec{q} = 0
\end{equation}

\noindent where
\begin{equation}\label{r0def}
\vec{r}_0 = (1/2)\vec{\delta} - \vec{\delta}\times\vec{q}/q^2;
\end{equation}

\noindent these two equations are exactly the content of \eqref{alphbetgamcentralaxorig},
which says that $\vec{r}-\vec{r}_0$ is parallel to $\vec{q}$ when $\vec{r}_0$ is given by \eqref{r0def}.}

In these equations the coordinates subtracted from $x,y,z$ {\it (that is, the components of $\vec{r}_0$)}
are precisely those of the vertex of an isosceles triangle normal to the plane of the two relative
axes of translation and of rotation passing through the origin of the coordinate axes, raised on\
a base perpendicular to the relative axis of rotation, cutting at its midpoint the relative axis
of translation and equal in length to the translation itself, the vertex angle being equal to the
amplitude of rotation $\theta$; all this in agreement with the construction given previously.

{\it I find this description difficult to visualize.  Again I resort to the equations to find the 
geometric meaning of $\vec{r}_0$.

From \eqref{delt}, \eqref{delVUsin}, and \eqref{delUUcos}, let us extract the half-angle formula for $\vec{\delta}$:
\begin{equation}\label{delhalf}
\vec{\delta} = \vec{t}+ 2\vec{U}\sin^2(\theta/2) -2U\hat{V}\sin(\theta/2)\cos(\theta/2).
\end{equation}

\noindent Then, since $\vec{q}/q^2 = \hat{t}/(2\tan(\theta/2))$,
\begin{eqnarray}\label{delq/qsq}
[\vec{\delta}\times\vec{q}/q^2 & = & \vec{\delta}\times\hat{t}/(2\tan(\theta/2)) \nonumber  \\
& = & [2\vec{U}\sin^2(\theta/2)- 2U\hat{V}\sin(\theta/2)\cos(\theta/2)]\times\hat{t}/(2\tan(\theta/2)) \nonumber\\
&=&[\vec{U}\sin(\theta/2)\cos(\theta/2) - U\hat{V}\cos^2(\theta/2)]\times\hat{t}\nonumber \\
& = & -U\hat{V}\sin(\theta/2)\cos(\theta/2) - \vec{U}\cos^2(\theta/2)
\end{eqnarray}

\noindent Now \eqref{delq/qsq} and \eqref{r0def} combine harmoniously to yield
\begin{eqnarray}\label{r0calc}
\vec{r}_0 & = & (1/2)\vec{\delta} - [- U\hat{V}\sin(\theta/2)\cos(\theta/2) - \vec{U}\cos^2(\theta/2)] \nonumber \\
& = & (1/2)[\vec{t} + 2\vec{U}\sin^2(\theta/2) -2U\hat{V}\sin(\theta/2)\cos(\theta/2)] + [U\hat{V}\sin(\theta/2)\cos(\theta/2) + \vec{U}\cos^2(\theta/2)] \nonumber \\
& = & (1/2)\vec{t} + \vec{U}[\sin^2(\theta/2) + \cos^2(\theta/2)] + U\hat{V} (\sin(\theta/2)\cos(\theta/2)- \sin(\theta/2)\cos(\theta/2)] \nonumber \\
& = & (1/2)\vec{t} +\vec{U}
\end{eqnarray}

\noindent saying that $\vec{r}_0$ is the midpoint of the line joining the point $\vec{U}$, 
the point on the central axis closest to the origin, to the point $\vec{U} + \vec{\Delta} = \vec{U}+ \vec{t}$
to which the point $\vec{U}$ is carried by the displacement.  Under NCA, $\vec{U} = 0$ and $r_0 = \vec{t}/2$.}

These equations {\it (referring to \eqref{alphbetgamcentralaxorig})} can also be written as follows,
on introducing the rotation angle and the direction of the axis of rotation:
\begin{eqnarray}\label{alphbetgamentralaxcot}
\frac{x-(1/2)\alpha - (1/2)\cot(\theta/2)(\beta\cos l-\gamma\cos h)}{\cos g} \nonumber \\
= \frac{y-(1/2)\beta- (1/2)\cot(\theta/2)(\gamma\cos g-\alpha\cos l)}{\cos h} \nonumber \\
= \frac{z-(1/2)\gamma - (1/2)\cot(\theta/2)(\alpha\cos h-\beta\cos g)}{\cos l} 
\end{eqnarray}

\noindent {\it These equations follow immediately from \eqref{alphbetgamcentralaxorig} by using \eqref{qdeforig}. 
In modern notation they become just another way of saying that $\vec{r}_0$ is given by \eqref{r0def}.}

\subsection{General equation of the central axis}.

But one may represent these three equations for the projections of the central axis {\it (it is
not clear which three equations are meant, since \eqref{alphbetgamentralaxcot} consists
of only two equations)} by a single equation with undetermined coefficients, namely
\begin{equation}\label{Deltabc}
\cos a\Delta x + \cos b\Delta y + \cos c\Delta z = t\cos(t,a,b,c)
\end{equation} 

\noindent where
\begin{equation}\label{tabc}
\cos(t,a,b,c)= \cos a\cos g + \cos b\cos h  +  \cos a\cos l,
\end{equation}

\noindent $a,b,c$ being the three angles formed by any direction whatever with the 
coordinate axes.  {\it That is, if $\hat{d} = \cos a,\cos b,\cos c $ is a unit vector in any direction whatsoever, then
$\vec{\Delta}\cdot\hat{d} = t\hat{d}\cdot\hat{t} = \vec{t}\cdot\hat{d}$.  But this is
tantamount to saying that $\vec{\Delta}= \vec{t}$, which is just the property of
initial points $\vec{r}$ that lie on the central axis.}

\section {The case of an infinitesimal displacement.}

If one considers only an infinitesimal displacements, that is where 
$\Delta x,\Delta y, \Delta z$ are infinitely small and can be replaced by
$dx, dy, dz$ and correspondingly the constants $\alpha, \beta, \gamma, \theta$ 
are small of the same order, then, neglecting terms of the second 
or higher order, the equations \eqref{Deldelorig} and \eqref{qdeforig} reduce to
\begin{equation}\label{ddelorig}
dx = \alpha +py -nz, \;\;\;\;\;dy = \beta + mz-px, \;\;\;\;\;dz = \gamma + nx - my;
\end{equation}

\begin{equation}\label{qdori}
m = \theta\cos g,\;\;\;\;\;n = \theta\cos h,\;\;\;\;\;p = \theta\cos l;
\end{equation}

\noindent and for the equations of the central axis,
\begin{equation}\label{centd}
\frac{\theta x + \gamma\cos h - \beta\cos l}{\cos g} = \frac{\theta y + \alpha\cos l - \gamma\cos g}{\cos h} = \frac{\theta z + \beta\cos g - \alpha\cos h}{\cos l}.
\end{equation}

\section {}

Let us return to the general formulas \eqref{Deldelorig}.  From them one derives
\begin{equation}\label{18mod}
(\vec{\Delta} - \vec{\delta})^2 = 4\tan^2(\theta/2)\;[(\vec{\omega}-(1/2)\vec{\delta})^2 -  ((\vec{\omega}-(1/2)\vec{\delta})\cdot\hat{t})^2].
\end{equation}
 
 \noindent {\it I shall henceforth omit some of the more complex equations given in the original text,
 supplying only my transcription into modern notation.  The derivation of 
 \eqref{18mod} may proceed as follows: \eqref{Deldelorig} can be transcribed to
 \begin{equation}\label{Deldelmod}
 \vec{\Delta} = \vec{\delta} + 2\tan(\theta/2)\hat{t}\times(\vec{\omega}-(1/2)\vec{\delta})
 \end{equation}
 
 \noindent where $\vec{\delta}$ is given by \eqref{deltUV}.  To see how \eqref{Deldelmod} 
 leads to \eqref{18mod}, let us first define $\vec{\,\omega}' = \vec{\omega} - (1/2)\vec{\delta}$.
 Then \eqref{Deldelmod} can be written
\begin{equation} \label{Deldelmodpr}
\vec{\Delta} =  \vec{\delta} + 2\tan(\theta/2)\hat{t}\times\vec{\,\omega}'
\end{equation}

\noindent which leads immediately to 
\begin{equation}\label{18modprvar}
(\vec{\Delta} - \vec{\delta})^2 = 4\tan^2(\theta/2) \; (\hat{t} \times \vec{\,\omega}')^2, 
\end{equation} 

\noindent whereas \eqref{18mod} becomes
\begin{equation}\label{18modpr}
(\vec{\Delta} - \vec{\delta})^2 = 4\tan^2(\theta/2)\;[(\vec{\,\omega}'^2 -  (\vec{\,\omega}'\cdot\hat{t})^2].
\end{equation}

\noindent But \eqref{18modprvar} and \eqref{18modpr} are identical, from the familiar identity
$(\vec{a}\times\vec{b})^2 = a^2b^2 - (\vec{a}\cdot\vec{b})^2$.  Thus \eqref{18mod} 
is a consequence of previously derived formulas.}

To the equation \eqref{18modpr} {\it (which gives the magnitude of $\vec{\Delta} - \vec{\delta}$)}
one may adjoin the two equations
\begin{equation}\label{Dmindperpt} 
(\vec{\Delta} - \vec{\delta})\cdot\hat{t} = 0
\end{equation} 

\begin{equation}\label{Dmindperpompr} 
(\vec{\Delta} - \vec{\delta})\cdot\vec{\,\omega}' = 0
\end{equation} 

\noindent {\it (which give its direction, normal to both $\hat{t}$ and $\vec{\,\omega}'$; these
equations follow directly from \eqref{Deldelmodpr}.)}

But now let us replace $\vec{\omega}$ by its value $\vec{r} + (1/2)\vec{\Delta}$.
We thus obtain from \eqref{Deldelmodpr}
\begin{equation}\label{intermed}
\vec{\Delta} - \vec{\delta} = 2\tan(\theta/2)\hat{t}\times\vec{r} + \tan(\theta/2)\hat{t}\times(\vec{\Delta} - \vec{\delta}) 
\end{equation}

\noindent {\it (I have dissected my own composite $\vec{\,\omega}'$ as well as making the author's substitution)}
and from this in turn
\begin{equation}\label{Delsinsinhalfsq} 
\vec{\Delta} = \vec{\delta} + (\sin\theta)\hat{t}\times\vec{r} + 2\sin^2(\theta/2)[\hat{t}\vec{r}\cdot\hat{t}-\vec{r}],
\end{equation}

\noindent of which only the first two terms will survive on passing from finite to infinitesimal displacements.

{\it In order to derive \eqref{Delsinsinhalfsq}, one must substitute \eqref{intermed} into itself as follows:
Evaluate
\begin{eqnarray}\label{intermedcross}
\hat{t}\times(\vec{\Delta} - \vec{\delta}) & = & \hat{t}\times[2\tan(\theta/2)\hat{t}\times\vec{r}] + \hat{t}\times [\tan(\theta/2)\hat{t}\times(\vec{\Delta} - \vec{\delta})] \nonumber \\
& = & 2\tan(\theta/2)[\hat{t}\hat{t}\cdot\vec{r} - \vec{r}] + \tan(\theta/2)[\hat{t}\hat{t}\cdot(\vec{\Delta} - \vec{\delta}) - (\vec{\Delta} - \vec{\delta})] \nonumber \\
& = & 2\tan(\theta/2)[\hat{t}\hat{t}\cdot\vec{r} - \vec{r}] + \tan(\theta/2)[0 - (\vec{\Delta} - \vec{\delta})]. 
\end{eqnarray}

\noindent Now substitute this expression into the last term of \eqref{intermed}:
\begin{equation}\label{intermedsubs}
\vec{\Delta} - \vec{\delta} = 2\tan(\theta/2)\hat{t}\times\vec{r} + 2\tan^2(\theta/2)(\hat{t}\hat{t}\cdot\vec{r} - \vec{r}) - \tan^2(\theta/2)(\vec{\Delta} - \vec{\delta}).
\end{equation}

\noindent Transposing the last term to the left side of the equation, and noting that $1+\tan^2(\theta/2) = 1/\cos^2(\theta/2)$, we find
\begin{equation}\label{intermedtrans}
(\vec{\Delta} - \vec{\delta}) / \cos^2(\theta/2) = 2\tan(\theta/2) \hat{t}\times\vec{r} + 2\tan^2(\theta/2)(\hat{t}\hat{t}\cdot\vec{r} -\vec{r})
\end{equation} 

\noindent Finally multiplying the equation by $\cos^2(\theta/2)$, we obtain\eqref{Delsinsinhalfsq}.}
   
\subsection{Expressions for finite displacements as rational functions of $\vec{\delta}$ and $\vec{q}$} 

If one retains in these formulas the primitive constants $m,n,p$, and directly extracts the values of $\Delta x, \Delta y, \Delta z$..., {\it (the author now gives formulas which
I shall abbreviate by recalling that these three $\Delta$'s form a single vector $\vec{\Delta}$, that $\alpha,\beta,\gamma$ form the vector $\vec{\delta}$, and that $m,n,p$ 
are the components of $\vec{q} = q\hat{t}$, and which I shall derive as follows:
Let us write all the trigonometric functions in \eqref{intermedtrans} in terms of $\tan(\theta/2)$, thus:
\begin{equation}\label{intermedtanhalf}
\vec{\Delta} - \vec{\delta} = [2\tan(\theta/2)\hat{t}\times\vec{r} + 2\tan^2(\theta/2)(\hat{t}\hat{t}\cdot\vec{r} -\vec{r})] / (1 + \tan^2(\theta/2)).
\end{equation}

\noindent Then recalling the definition \eqref{qdef}, we have
\begin{equation}\label{Deldelq}
\vec{\Delta} = \vec{\delta}+ [\vec{q}\times\vec{r} + (1/2)(\vec{q}\vec{q}\cdot\vec{r} -q^2\vec{r})](1+ (1/4)q^2)
\end{equation}

\noindent in which the displacement $\vec{\Delta}$ is given as a \emph{rational} function of $\vec{\delta}$ and $\vec{q}$.
That is, in the words of the author)},  

one has the following expressions {\it (for $\Delta x, \Delta y, \Delta z$)} in 
{\bf rational} functions of the six constants $\alpha,\beta,\gamma, m, n, p$: {\it (the expressions he gives
are essentially those of \eqref{Deldelq} written out in components.)}

\subsection{Important consequence relating to the formulas for transforming rectangular coordinates.}

Comparing these expressions with those that one would obtain by considering 
the transformation of rectangular coordinates, as this will be indicated in {\bf 26},
one obtains a way of reducing the nine coefficients that enter into the 
formulas for this transformation to three independent variables $m,n,p$,
entirely free of irrationals, which I believe has not been given before.

{\it Until now,the author has been occupied with calculating $\vec{\Delta}$, the 
straight-line displacement from an initial point $\vec{r}$ to its destination $\vec{\,r}'$.
He now turns his attention to the \emph{transformation} of $\vec{r}$ to $\vec{\,r}'$, 
in the case of a pure rotation.  Nowadays we think of this transformation as given
by a $3\times3$ matrix transforming $(x,y,z)$ to $(x',y',z')$.  The author
frames this in a complementary way as making a transformation of the
\emph{coordinate system} $\hat{x},\hat{y},\hat{z}$, such that
if $\vec{r} = x\hat{x} + y\hat{y} + z\hat{z}$ then $\vec{\,r}'  =  x\hat{\,x}' + y\hat{\,y}' + z\hat{,z}'$
- the coordinates, not the components, being altered.  The \lq\lq nine coefficients" he speaks of
are the nine elements of the matrix

\[
(M) =	\left(
		\begin{array}{ccc}
		\hat{x}\cdot\hat{x}' & \hat{y}\cdot\hat{x}' & \hat{z}\cdot\hat{x}'  \\
		\hat{x}\cdot\hat{y}' & \hat{y}\cdot\hat{y}' & \hat{z}\cdot\hat{y}'  \\ 
		\hat{x}\cdot\hat{z}' & \hat{y}\cdot\hat{z}' & \hat{z}\cdot\hat{z}'  
		\end{array}
	\right).
\]

\noindent The author gives formulas for these nine elements which have three notable 
features: (1) they involve no input other than the three parameters $m,n,p$;
(2) they contain no irrational expressions; and (3) the three parameters are completely 
independent.  He simply lists these nine formulas; I shall display them in
matrix form, but I shall write each matrix element exactly as it appears in the
original text except that in order to save space I shall write $mn/2$ rather than
$\frac{mn}{2}$, etc., and I shall place the common denominator $1 + (m^2+n^2 + p^2)/4$
outside the matrix as a prefix.  Here is the result:}

\[
(M)=	[1+ (m^2 + n^2 + p^2)/4]^{-1} \left(
		\begin{array}{ccc}
\		 1+(m^2 -n^2-p^2)/4 & mn/2 -p & pm/2 + n \\
		 mn/2 + p &1+  (n^2 - p^2 - m^2)/4 & np/2 - m \\ 
		pm/2 - n &  np/2 + m  & 1+ (p^2 - m^2 - n^2)/4
		\end{array}
	\right).
\]

On eliminating $m,n,p$ from the {\it (off-diagonal)} coefficients, one obtains the
 formulas of Monge, in {\bf irrational} functions, for the three {\it (diagonal)} coefficients.

{\it It is of some interest to decompose $(M) = (M)_1 + (M)_2 + (M)_3$ where
\[
(M)_1 =	[1 - (m^2 + n^2 + p^2)/4]/[1+ (m^2 + n^2 + p^2)/4] \left(
		\begin{array}{ccc}
\		1 & 0 & 0 \\
		0 &  1 & 0\\ 
		0 & 0 &  1 
		\end{array}
	\right),
\]

\[
(M)_2 =	[1+ (m^2 + n^2 + p^2)/4]^{-1} \left(
		\begin{array}{ccc}
\		 m^2/2  & mn/2  & pm/2  \\
		 mn/2 &  n^2/2  & np/2  \\ 
		pm/2 &  np/2  &   p^2/2 
		\end{array}
	\right),
\]

\[
(M)_3 =	[1+ (m^2 + n^2 + p^2)/4]^{-1} \left(
		\begin{array}{ccc}
\		0 &  -p &  n  \\
		  p & 0 &  - m \\ 
		- n & m & 0
		\end{array}
	\right).
\]

Passing over to modern vector notation, this gives us
\begin{equation}\label{rtorpr}
\vec{\,r}' = (M)\cdot\vec{r} = (1+ q^2/4)^{-1} [(1- q^2/4)\vec{r} + (1/2)\vec{q}\vec{q}\cdot\vec{r} + \vec{q}\times\vec{r}].
\end{equation}
 
\noindent Now  if we readmit the trigonometric functions via \eqref{qdef}, we have
\begin{equation}\label{trig1}
[1 - q^2/4]/(1+ q^2/4) = (1-\tan^2(\theta/2))/(1+\tan^2(\theta/2)) = \cos\theta,
\end{equation}

\begin{equation}\label{trig2}
(1/2)\vec{q}\vec{q}/(1+ q^2/4) = 2\hat{t}\hat{t}\tan^2(\theta/2)/\sec^2(\theta/2) = 2\hat{t}\hat{t}\sin^2(\theta/2) = \hat{t}\hat{t}(1-\cos\theta), 
\end{equation}

\begin{equation}\label{trig3}
\vec{q}/(1+ q^2/4) = 2\hat{t}\tan(\theta/2)\cos^2(\theta/2) = \hat{t}\sin\theta;
\end{equation}

\noindent and so we arrive at the celebrated Rodrigues rotation formula}
\begin{equation}\label{RRot}
\vec{\,r}' = \vec{r}\cos\theta + \hat{t}\hat{t}\cdot\vec{r} (1-\cos\theta) +\hat{t}\times\vec{r}\sin\theta.
\end{equation}

\section{On the composition of two displacements.}

{\it The following three paragraphs summarize the entire essay, looking backward
to the beginning as well as forward to the end.}

From a geometric consideration of the displacement of a solid system,
we began {\it (this essay)} by deducing the characteristic properties or general laws
of the displacement, which always reduces to a rotation followed
by a translation, or equivalently to a single couple of rotations
about two fixed axes. (The two axes may be either parallel or not parallel; in the first
case the displacement reduces to a simple translation, provided
that the two rotations are equal and in opposite directions.) 

From these properties we have now derived the analytic expression 
for the transformation, either finite or infinitely small, of the
coordinates of a solid system undergoing an arbitrary displacement.

{\it The author seems to be referring here to the formula for the matrix $(M)$, which
turned out in modern notation to be the rotation formula \eqref{RRot}.  But in the reasoning
to follow, he works not from \eqref{RRot} but from the earlier formula \eqref{Deldelmodq},
which gives the simple displacement $\Delta$ in terms of 
$\vec{\delta}$ and $\vec{q}$.  One readily sees that if one sets $\vec{q} = 0$ one
obtains $\Delta = \delta$, which is a pure translation since the components
$\alpha,\beta,\gamma$ of $\delta$ do not depend on the location of $\vec{r}$;
and that if $\delta = 0$ then $\Delta = \vec{\omega}\times \vec{r}$, which
describes a pure rotation since it vanishes at the origin $\vec{r} = 0$.}

We shall now deduce, from this expression, the laws of composition
of rotations and of translations that we previously exhibited synthetically.
And finally we shall establish these same formulas directly, by
three distinct analytic procedures, making use exclusively of the
invariance of the distances between points of this system.

Let us designate by $\vec{\Delta'}, \vec{\Delta''}$ the two successive changes of positions of a point in the displaced system, and by $\vec{\Delta}$ the resultant change.
Likewise by $\vec{\,\omega}', \vec{\,\omega}''$ the positions of the midpoints of the two straight-line displacements of the point, and by $\vec{\omega}$ the midpoint of the resultant displacement. Thus one has
\begin{equation}\label{ompr}
\vec{\,\omega}' = \vec{r}+(1/2)\vec{\Delta'},
\end{equation}

\begin{equation}\label{omprpr}
\vec{\,\omega}'' = \vec{r}+\vec{\Delta'}+ (1/2)\vec{\Delta'' },
\end{equation}

\noindent and

\begin{equation}\label{omtot}
\vec{\omega} = \vec{r}+(1/2)\Delta\vec{r}.
\end{equation}

\noindent Moreover,
\begin{equation}\label{adddels}
\vec{\Delta} = \vec{\Delta'} + \vec{\Delta''}
\end{equation}

\noindent {\it I am routinely putting all equations in modern form.}

Furthermore, let us designate by  $\hat{t}',\hat{t}''$ the absolute translations, and by $\theta', \theta''$ and $\hat{q}',\hat{q}''$ the rotations and the directions of the rotation axes of each of the two consecutive displacements under consideration; and by $\hat{T}, \Theta. \hat{Q}$ the analogous elements of the composite displacement.

First we shall examine separately the composition of simple translations
and that of rotations without translation.  If we suppose that the displacements
consist purely of translations, we have
\begin{equation}\label{transcomp}
\Delta' = \vec{t}', \Delta'' = \vec{\,t''}, \Delta = \vec{T} = \vec{t'} + \vec{\,t''},
\end{equation}

\noindent from which we see that the composite translation is nothing other 
than the third side of the triangle formed by the successive  passages of
a point of the system by reason of the two given translations.  {\it To compose
more than two given translations,} one may easily generalize from the triangle
 to the polygon of translations.

Now consider the composition of two rotations, without translation, about two axes that intersect at the origin of coordinates.  We shall have
\begin{eqnarray}\label{rotcomp}
\vec{\Delta'} & = & \vec{\,q}' \times \vec{\,\omega}', \nonumber \\
\vec{\Delta''} & = & \vec{\,q}'' \times \vec{\,\omega}'', \nonumber \\
\vec{\Delta} & = & \vec{\,Q} \times \vec{\,\omega},
\end{eqnarray}

\noindent where
\begin{eqnarray}\label{qmag}
\vec{q'} & = & \hat{q}' 2\tan(\theta'/2), \nonumber \\
\vec{q''} & = & \hat{q}'' 2\tan(\theta''/2), \nonumber \\
\vec{Q} & = & \hat{Q}' 2\tan(\Theta/2), \nonumber \\.
\end{eqnarray}

{\it To understand these equations, consider the third line of \eqref{rotcomp}.  
The part of $\omega$ parallel to $\vec{Q}$ does not contribute to $\vec{Q} \times \vec{\,\omega}$.
Therefore the expression $\;\hat{Q} \times \vec{\,\omega}$ simply rotates $\vec{\,\omega}$ through
a right angle in the plane perpendicular to $\vec{Q}$. This gives a vector with the right direction
(along $\vec{\Delta}$) but the wrong magnitude.  The magnitude of $\vec{\Delta}$ is obtained by
inserting a factor $2\tan(\Theta/2)$, that is by replacing $\hat{Q}$ with $\vec{Q}$ as indicated
in the third line of \eqref{qmag}.  Likewise for the first and second lines of these equations.

In all these discussions $\hat{q}$ is the direction of the rotation axis, and $\hat{t}$ of the translation.  But sometimes, as in {\bf 18}, it is assumed that one is dealing
with the central axis, so that $\hat{q} = \hat{t}$, Now, however, we have two partial rotation axes intersecting at $0$, which does not necessarily lie on the central
axis of either partial displacement; therefore  $\hat{q}$ must be distinguished from $\hat{t}$, for each partial displacement and for the resultant.  Indeed, for the
remainder of this Section the translations $\vec{t}',\vec{t}'',\vec{T}$ are taken as zero, so that $\hat{t}',\hat{t}'',\hat{T}$ are indeterminate.}

We need to determine the composite parameters $\Theta, \vec{Q}$ as functions of the partial ones $\theta',\vec{\,q}'$ and $\theta'',\vec{\,q}''$.
It will be useful to begin by eliminating the variables $\vec{\,\omega}',\vec{\,\omega}''$ from \eqref{rotcomp}, using the relations \eqref{ompr} through \eqref{adddels}.
{\it This elimination is the key to the author's solution of the problem posed.  He writes down the elimination formulas without showing a derivation; here I shall present a
geometrical derivation.

 Let $A,B,C$ denote respectively the points whose position vectors are \linebreak $\vec{r},\vec{r}+\vec{\Delta'},\vec{r} + \vec{\Delta}$. Then 
$\vec{\,\omega}',\vec{\,\omega}'',\vec{\omega}$
are the position vectors of the midpoints of $AB,BC,AC$. To eliminate $\vec{\,\omega}'$, consider the triangle formed by $A,\vec{\,\omega}',\vec{\omega}$.  It is 
similar to $ABC$ but half as large in linear dimension.  Therefore the vector $\vec{\,\omega}'-\vec{\omega}$ is just half the vector from $C$ to $B$ which is $-\vec{\Delta''}$.
So we have
\begin{equation}\label{elimompr}
\vec{\,\omega}' = \vec{\omega} - (1/2)\vec{\Delta''}.
\end{equation}

\noindent To eliminate $\vec{\,\omega}''$, consider the triangle formed by $\vec{\omega},\vec{\,\omega}'',C$.  It is also
similar to $ABC$ but half as large in linear dimension.  Therefore the vector $\vec{\,\omega}''-\vec{\omega}$ is just half the vector from $A$ to $B$ which is $+\vec{\Delta'}$.
So we have
\begin{equation}\label{elimomprpr}
\vec{\,\omega}'' = \vec{\omega} + (1/2)\vec{\Delta'}
\end{equation}

\noindent Note the opposite placement of $'$ and $''$, as well as the change of sign of the displacement term between \eqref{elimompr} amd \eqref{elimomprpr}.  These features will
have important consequences later on.}

Substituting \eqref{elimompr} and \eqref{elimomprpr} into \eqref{rotcomp}, we obtain
\begin{equation}\label{delnoompr}
\vec{\Delta'} = \vec{\,q}'\times\vec{\omega} - (1/2)\;\vec{\,q}'\times\vec{\Delta''}
\end{equation} 

\begin{equation}\label{delnoomprpr}
\vec{\Delta''} = \vec{\,q}''\times\vec{\omega} + (1/2)\;\vec{\,q}''\times\vec{\Delta'},
\end{equation} 

\noindent from which we deduce the following values for the partial displacements:
{\it The author writes down the \lq\lq following values" without showing the derivation, which is far from trivial.
It will be observed that \eqref{delnoompr} and \eqref{delnoomprpr} constitute a system of two linear equations
in the two unknowns $\vec{\Delta'}$ and $\vec{\Delta''}$. They can be solved by substituting each one into the
other.

Substituting \eqref{delnoomprpr} into \eqref{delnoompr}, we have
\begin{eqnarray}\label{DpDpp}
\vec{\Delta'} & = & \vec{\,q}'\times\vec{\omega} - (1/2)\;\vec{\,q}'\times[\vec{\,q}''\times\vec{\omega}] + (1/2)\;\vec{\,q}''\times\vec{\Delta'} \nonumber \\
& = & \vec{a}' + \vec{b}' - (1/4)\;\vec{\,q}' \times[\vec{\,q}''\times\vec{\Delta'}] 
\end{eqnarray}

\noindent where
\begin{equation}\label{abpr}
\vec{a}' = \vec{\,q}'\times\vec{\omega},\;\; \vec{b}' = - (1/2)\;\vec{\,q}'\times[\vec{\,q}''\times\vec{\omega}]
\end{equation}

\noindent and substituting in the reverse direction,
\begin{eqnarray}\label{DppDp}
\vec{\Delta''} & = & \vec{\,q}''\times\vec{\omega} + (1/2)\;\vec{\,q}''\times[\vec{\,q}'\times\vec{\omega}]  - (1/2)\;\vec{\,q}'\times\vec{\Delta''}   \nonumber \\
& = & \vec{a}'' + \vec{b}'' - (1/4)\;\vec{\,q}''\times[\vec{\,q}'\times\vec{\Delta''}] 
\end{eqnarray}

\noindent where
\begin{equation}\label{abprpr}
\vec{a}'' = \vec{\,q}''\times\vec{\omega}, \;\; \vec{b}'' = +(1/2)\;\vec{\,q}''\times[\vec{\,q}'\times\vec{\omega}].
\end{equation}
 
\noindent The exchange of order of the two partial displacements carries $a'$ into $a''$ but $b'$ into $-b''$; this change of sign is inherited from the one between \eqref{elimompr}
and  \eqref{elimomprpr}. We can express the change by rewriting $b' = +(1/2)\;\vec{\,q}'\times[\vec{\omega}\times\vec{\,q}'']$.

In view of \eqref{adddels}, we now have
\begin{eqnarray}\label{DD}
\vec{\Delta} & = & \vec{a} + \vec{b} - (1/4) [\vec{\,q}' \times[\vec{\,q}''\times\vec{\Delta'}] + \vec{\,q}''\times[\vec{\,q}'\times\vec{\Delta''} ] \nonumber \\
& = & \vec{a} + \vec{b} + (1/4) \vec{\,q}'\cdot \vec{\,q}'' [\vec{\Delta'} \ + \vec{\Delta''} ] = \vec{a} + \vec{b} + (1/4) \vec{\,q}'\cdot \vec{\,q}'' \vec{\Delta}
\end{eqnarray}

\noindent where we have used the fact that $\vec{\,q}'\cdot\Delta'\vec{r} =\vec{\,q}''\cdot\Delta''\vec{r} = 0$, and where
\begin{equation}\label{avec}
\vec{a} = \vec{a}' + \vec{a}''= (\,\vec{\,q}' + \,\vec{\,q}'')\times \vec{\omega}
\end{equation}  

\noindent and
\begin{equation}\label{bvec}
\vec{b} = \vec{b}' + \vec{b}'' = (1/2)\;[\vec{\,q}'\times(\vec{\omega}\times\vec{\,q}'') + \vec{\,q}''\times(\vec{\,q}'\times\vec{\omega})] = -(1/2)\;\vec{\omega}\times(\vec{\,q}''\times\vec{\,q}').
\end{equation} 

\noindent In the last step I have used the Jacobi identity for the cyclic sum of a double cross product.

The last term in \eqref{DD} can be transposed to the left side, giving \linebreak $\vec{\Delta}(1-(1/4)\; \vec{\,q}'\cdot \vec{\,q}'') = \vec{a} + \vec{b}$ or
\begin{equation}\label{Delqq}
\vec{\Delta}  = \frac{ (\,\vec{\,q}' + \,\vec{\,q}'')\times \vec{\omega} + (\vec{\,q}''\times\vec{\,q}')\times \vec{\omega}}{1-(1/4)\; \vec{\,q}'\cdot \vec{\,q}''}.
\end{equation}

\noindent The appearance of the {\bf common factor} $\vec{\omega}$ in \eqref{Delqq} is the major fruit of \eqref{elimompr} and \eqref{elimomprpr}.}

Now by comparison of this expression for $\Delta\vec{r}$ with the postulated expression for the same quantity in terms of $\Theta, \hat{T}$, we arrive at the following relation:
{\it The author wishes to compare \eqref{Delqq} with the third equation in \eqref{rotcomp}.  This comparison yields the equation
\begin{equation}\label{Drmatch}
\vec{Q}\times\vec{\omega} = \vec{s}\times\vec{\omega}
\end{equation}

\noindent where
\begin{equation}\label{sdef}
\vec{s}=\frac{\vec{\,q}'+ \vec{\,q}'' + (1/2)\vec{\,q}''\times\vec{\,q}'}{1-(1/4)\vec{\,q}'\cdot\vec{\,q}''}.
\end{equation}

\noindent One must now remember that $\vec{Q}$, as well as the variables entering into $\vec{s}$, are fixed for a particular displacement of the \emph{whole solid}, whereas 
$\vec{r}$, which enters into $\vec{\omega}$, can be \emph{any point} in the solid.  Therefore \eqref{Drmatch} holds for every possible $\vec{\omega}$.  This implies the author's \lq\lq following relation", namely}
\begin{equation}\label{Qs}
\vec {Q} = \vec{s},
\end{equation} 

\noindent from which one deduces, for the value of the resultant rotation,
\begin{equation}\label{cTH/2}
\cos(\Theta/2) = \cos(\theta'/2)\cos(\theta''/2) - \sin(\theta'/2) \sin(\theta''/2)\cos\nu,
\end{equation}

\noindent where $\nu$ is the angle between the two axes of rotation, thus $\cos\nu = \hat{q}'\cdot\hat{q}''$.

{\it The deduction may proceed as follows. First we replace the three $q$-vectors by their definitions:
generically, $\vec{q} = 2\hat{q} \tan(\theta/2)$.  Thus $\tan^2(\Theta/2)=Q^2/4$, and 
\begin{equation}\label{cossqTh}
\cos^2(\Theta/2) = (1+ (Q^2/4))^{-1}.
\end{equation}

\noindent Meanwhile we may write $\vec{s} = 2\vec{N}/D$, so that substituting \eqref{Qs} into \eqref{cossqTh} we have
\begin{equation}\label{ThND}
\cos^2 (\Theta/2) = (1+ (s^2/4))^{-1} = \frac{D^2}{N^2+D^2}
\end{equation}

\noindent where
\begin{equation}\label{Ndef}
\vec{N} = \vec{q}'/2 + \vec{\,q}''/2 + \vec{q}''\times\vec{q}'/4,
\end{equation}

\begin{equation}\label{Ddef}
D = 1 - (1/4)\cos\nu.
\end{equation}

From \eqref{ThND}  we see that if $N^2 + D^2$ can be exhibited as a perfect square we can obtain $\cos(\Theta/2)$ without radicals.
The third term of \eqref{Ndef} is orthogonal to both of the first two; therefore
\begin{eqnarray}\label{Nsq}
N^2 & = & [\vec{\,q}'/2 + \vec{\,q}''/2)]^2 + (\vec{\,q}''\times\vec{\,q}')^2] /16 \nonumber \\
& = & [\tan^2(\theta'/2)\! + \! \tan^2(\theta''/2) + 2\tan(\theta'/2)\tan(\theta''/2)\cos\nu]  \nonumber \\
& + & \tan^2(\theta'/2)\tan^2(\theta''/2)(1-\cos^2\nu);
\end{eqnarray}

\noindent  of course $(\hat{q}''\times\hat{q}')^2 = \sin^2\nu = 1 - \cos^2\nu$. At the same time \eqref{Ddef} gives us
\begin{equation}\label{Dsq}
D^2 = 1 - 2\tan(\theta'/2) \tan(\theta''/2) \cos\nu + \tan^2\theta'/2) \tan^2(\theta''/2) \cos^2\nu.
\end{equation}

\noindent The terms in $\cos\nu$ and $\cos^2\nu$ obligingly cancel between \eqref{Nsq} and \eqref{Dsq}, leaving
\begin{eqnarray}\label{NsqDsq}
N^2 + D^2 = 1 + \tan^2(\theta'/2) +  \tan^2(\theta''/2) +  \tan^2(\theta'/2)\tan^2(\theta''/2) \nonumber\\
 = (1+ \tan^2(\theta'/2))(1+ \tan^2(\theta''/2)) = 1/[\cos^2(\theta'/2) \cos^2(\theta''/2)].
\end{eqnarray}

\noindent This enables us to take the square root of \eqref{ThND}:
\begin{eqnarray}\label{coshalfTH}
\cos(\Theta/2) & = & \cos(\theta'/2)\cos(\theta''/2)(1 - \tan(\theta'/2) \tan(\theta''/2) \cos\nu) \nonumber \\
& = & \cos(\theta'/2)\cos(\theta''/2) - \sin(\theta'/2) \sin(\theta''/2)\cos\nu,
\end{eqnarray}

\noindent in agreement with \eqref{cTH/2}.  This startling relation was deduced from \eqref{Qs} by the author, who shows  no intermediate steps in the text. 
One can only wonder at his ability to navigate the maze of substitutions without the help of our vector relations.}

- and for the inclination of the resultant axis,
\begin{equation}\label{sTH/2}
\hat{Q}\sin(\Theta/2)=\hat{q}'\sin(\theta'/2)\cos(\theta''/2)+\hat{q}''\sin(\theta''/2) \cos(\theta'/2)+\hat{q}''\times\hat{q}''\sin(\theta'/2)\sin(\theta''/2).
\end{equation} 

{\it Recalling the definitions given after \eqref{rotcomp},
\begin{equation}\label{qgens}
\vec{Q} = 2\hat{Q}\tan(\Theta/2),\;\;\vec{q\,}' = 2\hat{q}''\tan(\theta'/2)\,\;\;\vec{q\,}'' = 2\hat{q}''\tan(\theta''/2),
\end{equation}

\noindent we can write \eqref{Qs} as
\begin{equation}\label{Qstan}
\hat{Q}\tan(\Theta/2) =  \frac{\hat{q}'\tan(\theta'/2)+ \hat{q}''\tan(\theta''/2) + \hat{q}''\times \hat{q}'\tan(\theta'/2)\tan(\theta''/2)}{1-\hat{q}''\cdot\hat{q}'\tan(\theta'/2)\tan(\theta''/2))}
\end{equation}

\noindent and \eqref{cTH/2} as
\begin{equation}
\cos(\Theta/2) = \cos(\theta'/2)\cos(\theta''/2)(1-\hat{q}''\cdot\hat{q}'\tan(\theta'/2)\tan(\theta''/2)).
\end{equation}

Multiplying these two equations, we obtain  \eqref{sTH/2}.}

%
%
%

In these formulas one notices immediately that the order of succession of the rotations $\theta', \theta''$ has no effect on the amplitude $\Theta$
of the resultant rotation {\it (see \eqref{coshalfTH})}, but that it does affect the direction of the axis of that rotation {\it (see the last term of \eqref{sTH/2})},
unless the rotations $\theta',\theta''$ are infinitely small.

Now, the expression for $\cos(\Theta/2)$ holds if $\Theta/2$ is an angle of a spherical triangle, of which the opposing side is $\nu$ and the two
other angles are  $\theta'/2,\theta''/2$. {\it This situation is dual to the one described by the usual spherical law of cosines, which gives a side in terms of the
opposite angle and the two other sides.  That is why the formula \eqref{coshalfTH} has a minus sign before the second term instead of a plus sign.}
To specify more precisely the position of the resultant axis in relation to the two given axes, let us suppose, which is always possible, that
\begin{equation}\label{spcoords}
\cos l' = \cos l'' = 0, \cos h' = 0, \cos g' = 1, \cos g'' = \cos\nu, \; \cos h'' = \sin\nu.
\end{equation} 

{\it Generally, the author has regarded the $x,y,z$ coordinate system as fixed and independent of what displacement is being applied to the
solid system.  Occasionally he has made slight departures: in the previous section he discussed linear orthogonal transformations of the coordinates,
and in the present inquiry he has already specified that the origin of coordinates is at the intersection of the two given (and hence also of the 
resultant) axes of rotation. Now he chooses an orientation of the coordinate system completely tailored to the problem at hand.  He has
the $x$-axis coinciding with the axis of the rotation $\theta'$, and the $z$-axis perpendicular to both rotation axes $\theta,\theta''$. Thus
the axis of rotation $\theta''$ lies in the $x$-$y$ plane making an angle $\nu$ with the $x$-axis. This will simplify his calculations, but at the cost of
making it impossible to apply modern vector notation to his formulas; they will have to be written out in components as in the original text.}

One now has
\begin{eqnarray}\label{sp1}
\sin(\Theta/2)\cos G & = & \sin(\theta'/2)\cos(\theta''/2) + \sin(\theta''/2)\cos(\theta'/2)\cos\nu, \nonumber \\
\sin(\Theta/2)\cos H & = & \sin(\theta''/2)\cos(\theta'/2)\sin\nu, \nonumber \\
\sin(\Theta/2)\cos L & = & \sin(\theta'/2)\sin(\theta''/2)\sin\nu. 
\end{eqnarray}

\noindent {\it These equations can be obtained by resolving \eqref{sTH/2} into its $x$-,$y$-, and $z$-components, and applying \eqref{spcoords}.}

If the order of rotations is reversed, the only change in \eqref{sp1} is that the sign of $\cos L$ becomes negative, from which it follows that the
new resultant axis is placed in a \emph{symmetric} position to the old one relative to the plane of the two given axes.

{\it This result can be obtained directly from \eqref{sTH/2} more easily than by struggling through \eqref{sp1}. The sum of the first
two terms on the right side of \eqref{sTH/2} is symmetric in $\hat{q}'$ and $\hat{q}''$, and both lie in the 
$\hat{q}'$-$ \hat{q}''$ plane; the last term is antisymmetric and perpendicular to this plane.  Therefore the interchange
of $\hat{q}'$ with $\hat{q}''$ causes $\hat{Q}$ to be reflected in the $\hat{q}'$-$\;\hat{q}''$ plane.}

If we denote by $H'$ the angle formed by the resultant axis with the axis of the rotation $\theta''$, we shall have
\begin{equation}\label{sp2}
\cos H' = \cos G\cos g'' + \cos H\cos h'' = \frac{\sin(\theta''/2)\cos(\theta'/2) + \sin(\theta'/2)\cos(\theta''/2)\cos\nu}{\sin(\Theta/2)}
\end{equation}

\noindent {\it (the first equality comes from the decomposition $\hat{t}'' = \hat{x}\cos g'' + \hat{y}\cos h''$, the second
by applying \eqref{spcoords} to \eqref{sp1})};  besides this we have
\begin{eqnarray}\label{sp3}
&&  \sin^2(\Theta/2)  = \sin^2(\theta'/2) \sin^2\nu + [\sin(\theta''/2)\cos(\theta'/2)+\sin(\theta'/2)\cos(\theta''/2)\cos\nu]^2 \nonumber \\
& = & \sin^2(\theta''/2) \sin^2\nu + [\sin(\theta'/2)\cos(\theta''/2)+\sin(\theta''/2)\cos(\theta'/2)\cos\nu]^2
\end{eqnarray}

\noindent and consequently
\begin{equation}\label{sp4}
\sin^2 G = \frac{\sin^2(\theta''/2)\sin^2\nu}{\sin^2(\Theta/2)}, \;\; \sin^2 H' = \frac{\sin^2(\theta'/2)\sin^2\nu}{\sin^2(\Theta/2)},
\end{equation}

\noindent equations that establish the proportionality of the sines of the half-rotations 
to those of the angles formed by the resultant axis with the given axes inversely corresponding 
{\it (that is, $G$ to $\theta''$ and $H'$ to $\theta'$)}, and which lead to the construction 
that we have indicated from the outset for the composition of rotations.  

{\it The first line of \eqref{sp3} is obtained by summing the squares of the three equations in \eqref{sp1}. The second line follows from
the first in view of \eqref{coshalfTH} which gives $\Theta$ as a symmetric function of $\theta'$ and $\theta''$.  In regard to \eqref{sp4},
note that $G$ bears the same relation to $\theta'$ as $H'$ to $\theta''$, in view of the condition $\cos g' = 1$ from \eqref{spcoords}.
Of course the author means to consider the square root of \eqref{sp4}.  The resulting formula is the one already described in
Section 6, corollary \# 2.}

If we were to follow an analogous procedure for the composition of rotations about 
an arbitrary number of intersecting axes, the resulting formulas would be rendered 
exceedingly complex by the terms of second order; hence we omit that subject.

\section{On the analytic composition of rotations about nonintersecting axes.}

As for the composition of rotations about nonintersecting axes, and generally of an arbitrary succession of displacements
of a solid system, given by individual displacements $\Delta',\Delta'',etc.$ whose analytic form is known, we shall have
\begin{eqnarray}\label{nonintcomp}
\Delta x = \Delta' x + \Delta'' x  + \Delta''' x  + ... = A + 2\tan(\Theta/2) (Y\cos L - Z\cos H), \nonumber \\
\Delta y = \Delta' y + \Delta'' y  + \Delta''' y  + ... = B + 2\tan(\Theta/2) (Z\cos G - X\cos L), \nonumber \\
\Delta z = \Delta' z + \Delta'' z  + \Delta''' z  + ... = C + 2\tan(\Theta/2) (X\cos H - Y\cos G).
\end{eqnarray}

The constants $A,B,C\; etc.$ {\it (pertaining to the resultant displacement)} are to be found in terms of analogous constants
$A', B', C' etc.$; ; $A'', B'', C'' etc.$; ... belonging to each of the consecutive displacements to be combined, and of the
other elements of these displacements.  {\it The expressions beginning $A$,$B$,$C$ are closely related to those
seen in \eqref{DelGam+omLqorig} and \eqref{qdeforig}.} But it is evident that the elements $\Theta,G,H,L$ of the
resultant rotation depend only on the \emph{rotational} elements of these displacements, just as we have seen earlier
from geometric considerations. 

{\it Having posed a problem of sweeping generality, the author proceeds to solve only the simplest case, that of composing
{\bf two} rotations about nonintersecting axes. Moreover, he will select a made-to-order orientation of the coordinate system
similar to the one he introduced to treat two intersecting axes.  Whereas in the previous section both axes passed
through the point $(0,0,0)$, and \eqref{spcoords} gave their direction cosines as $(1,0,0)$ and $(\cos\nu,\sin\nu,0)$,
in the nonintersecting case \eqref{spcoords} will still give the direction cosines but the second axis will pass not through
$(0,0,0)$ but through $(0,0,u)$ where $u$ is geometrically the closest distance of approach between the two axes.}

Let us take, for example, the composition of rotations about two fixed nonintersecting axes, one identified with the
$x$-axis and the other normal to the $z$-axis and cutting that axis at a distance $u$ from the origin. We shall, as before,
denote by $\nu$ the angle between these two rotation axes and by $\theta',\theta''$ the amplitudes of the respective
rotations. To begin with, we shall have for the amplitude of the resultant rotation and the direction of its axis, the same
formulas as obtained in the previous section. To fix the position of the central {\it  (i.e. composite)} axis, as well as to find the displacements
of the coordinates, one needs only to calculate the displacements $\alpha,\beta,\gamma$ from the 
origin of coordinates. Now, those displacements arising from the first
rotation about the $x$ axis are null; hence it suffices to calculate those arising from the rotation $\theta''$, for which
one has in general {\it (cf. \eqref{Deldelorig} near the end of Section 15)}
\begin{eqnarray}\label{nonintsec}
\Delta'' x = \alpha + 2\tan(\theta''/2) [(Y'' - (1/2)\beta)\cos l'' - (Z'' - (1/2)\gamma)\cos h''] \nonumber \\
\Delta'' y = \beta + 2\tan(\theta''/2) [(Z'' - (1/2)\gamma)\cos g'' - (X'' - (1/2)\alpha)\cos l''] \nonumber \\
\Delta'' z = \gamma + 2\tan(\theta''/2) [(X'' - (1/2)\alpha)\cos h'' - (Y'' - (1/2)\beta)\cos g'';
\end{eqnarray}

\noindent these displacements must vanish for all points on the axis of rotation $\theta''$, for which one has
\begin{equation}\label{spnonintcoords}
\cos l'' = 0,\; \cos h'' = \sin\nu,\; \cos g'' = \cos\nu, \;Y'' = X''\tan\nu, \;Z''= u.
\end{equation} 

\noindent From this there result the following values for $\alpha,\beta,\gamma$, which one could easily
derive from the construction itself,
\begin{equation}\label{alphbetgammnonint}
\alpha = u\sin\nu\sin\theta'',\;\;\beta = -u\cos\nu\sin\theta'',\;\; \gamma = 2u\sin^2(\theta''/2),
\end{equation}

\noindent in agreement with the theorems of Section 13.

Thus the resultant axis and the two given divergent axes become more nearly parallel to a single plane, the greater
the distance between the two latter axes in comparison with the absolute translation of the resultant displacement.
{\it That is, $\cos L$ approaches zero as $u>>T$.

Since the two rotation axes do not intersect, Euler's fixed point theorem (Section {\bf 3}) does not hold, and the resultant displacement
is not a pure rotation but contains an absolute translation $\vec{\,T} = T\hat{T}$.

In \eqref{nonintsec} the coordinates of the point $\xi,\eta,\zeta$ defined in \eqref{xietazetaorig}
have been replaced in \eqref{Deldelorig} by $X,Y,Z$ as given in \eqref{spnonintcoords}.  Note that in the case $u=0$
(two intersecting axes) $X,Y,Z$ would be strictly proportional to $\cos g'',\cos h'',\cos l''$ and consequently the 
terms in these quantities would vanish; this would leave only $(\alpha,\beta,\gamma)$ on the right side 
of \eqref{nonintsec} so that the left side is made to vanish by setting $(\alpha,\beta,\gamma) = 0$.  This is why 
$(\alpha,\beta,\gamma)$ do not appear in Section 19.

For $u\neq 0$, only the terms in $X,Y$ cancel out, and to solve for $\alpha,\beta,\gamma$ it suffices to make
$\Delta''(x,y,z)$ vanish when $(X,Y,Z) = (0,0,u)$. This yields (using \eqref{spnonintcoords} as well)
\begin{eqnarray}\label{nonintseceqs}
0 = \alpha + 2\tan(\theta''/2)[ - (u-(1/2)\gamma)\sin\nu] \nonumber \\
0 = \beta + 2\tan(\theta''/2)[(u-(1/2)\gamma)\cos\nu] \nonumber \\
0 = \gamma + 2\tan(\theta''/2)[-(1/2)\alpha\sin\nu - (-(1/2)\beta\cos\nu]
\end{eqnarray}

\noindent The first two lines give
\begin{equation}\label{alphbetnu}
\alpha = (2u-\gamma)\tan(\theta''/2)\sin\nu,\;\;\;\beta = - (2u-\gamma)\tan(\theta''/2)\cos\nu
\end{equation}

\noindent and the third gives
\begin{equation}\label{gamalphbet}
\gamma = (2u-\gamma)\tan^2(\theta''/2)(\sin^2\nu + \cos^2\nu) = 2u\tan^2(\theta/2) - \gamma\tan^2(\theta/2).
\end{equation} 

\noindent Transposing the last term to the left, and multiplying the equation by $\cos^2(\theta''/2$, we have
\begin{equation}\label{gamyes}
\gamma = 2u\tan^2(\theta''/2)\cos^2(\theta''/2) = 2u\sin^2(\theta''/2),
\end{equation}
 
 \noindent and substituting \eqref{gamyes} into \eqref{alphbetnu} gives
 \begin{equation}\label{alphbettrue}
 \alpha = u\sin\nu\sin\theta'',\;\;\; \beta = - u\cos\nu\sin\theta''
 \end{equation} 

\noindent as in $\eqref{alphbetgammnonint}$.  (The text has a mistake: $\sin\nu$ instead of $\cos\nu$ in $\beta$.)

The analogy between this calculation and the one discussed in the last paragraphs of Section 15 may be
obscured by the fact that then only a single rotation was under consideration, whereas now there are two. 
But there is no difference if one makes the appropriate correspondences.  The role played by the \emph{central} axis
in Section 15 is here played by the \emph{second} of the two axes, having direction $\hat{q}''$ and rotation angle $\theta''$.
Call this axis the \emph{active} axis.  In both cases the active axis does not pass through the origin of coordinates.  
In both cases a \emph{substitute} axis is introduced, parallel to the active axis but passing through the origin of coordinates.  
In both cases $(\alpha,\beta,\gamma) = \vec{\delta}$ is the displacement of a point starting at the origin.  In the present
case there is an additional axis, the first, which plays no part in determining $\vec{\delta}$ because it lies on the substitute
axis and its action on the point of origin is null.

In the present case, the values of $\alpha,\beta,\gamma$ can be easily obtained because the direction $\hat{q}'' = (\cos\nu,\sin\nu,0)$
is given.  Starting at the origin $(0,0,0)$, the trajectory of rotation about the active axis describes an arc $\theta''$ of the circle 
whose center is at $(0,0,u)$ and whose plane is normal to $\hat{q}''$. This gives immediately $\gamma = u(1-\cos\theta'')= 2u\sin^2(\theta''/2)$
as well as $\sqrt{\alpha^2+\beta^2} = u\sin\theta''$. From the direction of $\hat{q}''$ we infer that $\alpha:\beta::\sin\nu:(-\cos\nu)$
and so $\alpha = u\sin\theta''\sin\nu'',\beta = - u\sin\theta''\cos\nu$. In this way \eqref{gamyes} and \eqref{alphbettrue} are
derived from geometry alone.}

From \eqref{alphbetgammnonint} one derives, for the value
$T$ of the absolute translation resulting from two rotations
about two nonintersecting fixed axes:
\begin{eqnarray}\label{Tnonintsec}
T & = & \alpha\cos G + \beta\cos H +\gamma\cos L \nonumber \\
& = &\frac{2u\sin\nu\sin(\theta'/2)\sin(\theta''/2)}{\sin(\Theta/2)} \nonumber \\
& = & 2u\cos L.
\end{eqnarray}

\noindent {\it This array contains three equalities.  To prove the first equality, one notes that $T=\vec{\,T}\cdot\hat{T}$
where $\vec{\,T}$ is the displacement of the composite axis along its length due to the successive actions of the
two partial rotations.  We now treat the composite axis just as we treated the second partial axis in deriving
\eqref{gamyes} and \eqref{alphbettrue}:  we introduce a \emph{substitute composite axis} parallel to the true
composite axis but passing through the origin.  Referring rotations to this substitute axis, we deduce as before
that $(\alpha,\beta,\gamma) = \vec{\delta}$ is the displacement undergone by the point originally at $(0,0,0)$,
that is, $\vec{\delta}$ is the translation, common to all points, that must be added to rotation about the substitute
axis in order to simulate the original composite displacement $\vec{\,T}$.  But rotation about the substitute axis, applied to
an arbitrary point, produces a trajectory perpendicular to its direction $\hat{T}$, so that $\vec{\,T}\cdot\hat{T} = \vec{\delta}\cdot\hat{T}$.  
Hence $T = \delta\cdot\hat{T}$, which is the first equality.

The second equality is obtained by takng $(\alpha,\beta,\gamma)$ from \eqref{alphbetgammnonint} and 
$\hat{T}$ from \eqref{sp1}. This gives
\begin{eqnarray}\label{deldotThat}
\alpha \cos G \sin(\Theta/2) & = & u\sin\nu\sin\theta'' [\sin(\theta'/2)\cos(\theta''/2) + \cos(\theta'/2)\sin(\theta''/2)\cos\nu] \nonumber \\
\beta \cos H \sin(\Theta/2) & = & - u\cos\nu\sin\theta'' [\cos(\theta'/2)\sin(\theta''/2)\sin\nu] \nonumber \\
\gamma \cos L \sin(\Theta/2) & = & 2u\sin^2(\theta''/2) [\sin(\theta'/2)\sin(\theta''/2)\sin\nu] 
\end{eqnarray}

\noindent The second term in the first line cancels the second line so that
\begin{equation}\label{alphG+betH}
\alpha \cos G \sin(\Theta/2) + \beta \cos H \sin(\Theta/2) = u\sin\nu\sin\theta'' [\sin(\theta'/2)\cos(\theta''/2)] 
 = 2u\sin\nu\sin(\theta'/2)\sin(\theta''/2)\cos^2(\theta''/2)
\end{equation}

\noindent and adding the third line we obtain the second equality of \eqref{Tnonintsec}.  The third
equality is obtained by simply comparing the numerator of the second line of \eqref{Tnonintsec}
with the third line of \eqref{sp1}.}

The equations of the resultant axis are obtained by substituting for $\alpha,\beta,\gamma$ their values from 
\eqref{alphbetgammnonint} in the general equations
\begin{eqnarray}\label{centrnonint}
\frac{x - (1/2)\alpha - (1/2)\cot(\Theta/2)(\beta\cos L-\gamma\cos H)}{\cos G} \nonumber \\
= \frac{y - (1/2)\beta - (1/2)\cot(\Theta/2)(\gamma\cos G-\alpha\cos L)}{\cos H} \nonumber \\
= \frac{z - (1/2)\gamma - (1/2)\cot(\Theta/2)(\alpha\cos H-\beta\cos G)}{\cos L}.
\end{eqnarray}

\noindent {\it These equations are the same as \eqref{alphbetgamcentralaxorig} with $g,h,l,\theta$ capitalized.}

In the case of infinitely small rotations these formulas are considerably simplified, and it is found
that the {\it (resultant)} central axis is parallel to the plane of the two contributing axes and intersects
their shortest distance.  Effectively, in this case, with neglect of infinitesimals of second order, one has
\begin{eqnarray}\label{nonintsmall}
\cos G &=& \frac{\theta' + \theta''\cos\nu}{\Theta}, \;\cos H = \frac{\theta''\sin\nu}{\Theta},\;\cos L = \frac{\theta'\theta''\sin\nu}{2\Theta}, \nonumber \\
\Theta^2 &=& \theta'^2 + \theta''^2 + 2\theta'\theta''\cos\nu, \nonumber \\
\alpha &=& u\theta''\sin\nu,\;\; \beta = - u\sin\theta''\cos\nu, \;\;\gamma = 0, \;\;T = \frac{u\theta'\theta''\sin\nu}{\Theta},
\end{eqnarray}

\noindent {\it (I have corrected three errors in the French text of this array:  (i) top line, $\cos H$ numerator,
 $\sin\nu$ incorrectly given as $\sin^2\nu$; (ii) second line, left side of equation, $\Theta^2$ incorrectly given
as $\Theta$; (iii) bottom line, $T$ numerator incorrectly given as $u\theta\theta\theta\sin\nu$)}

\noindent and for the equations of the {\it (resultant)} central axis, 
\begin{equation}\label{nonintsmcent}
y = \frac{x\theta''\sin\nu}{\theta'+\theta''\cos\nu},\;\;\; z = \frac{u\theta''(\theta''+\theta'\cos\nu)}{\Theta^2}.
\end{equation}

\section{Composition of successive rotations about three perpendicular axes.}

We shall end this subject by giving the following formulas for the composition 
of three successive rotations $\theta,\theta',\theta''$ about the three coordinate 
axes $x,y,z$.  The elements $\Theta,G,H,L$ of the composite rotation are
expressed as follows:
\begin{equation}\label{xyzTh}
\cos(\Theta/2) = \cos(\theta/2)\cos(\theta'/2)\cos(\theta''/2) - \sin(\theta/2)\sin(\theta'/2)\sin(\theta''/2),
\end{equation}

\begin{eqnarray}\label{xyzGHL}
\sin^2 G & = & \frac{1-\cos\theta'\cos\theta''}{2\sin^2(\Theta/2)}, \nonumber \\
\sin^2 H & = & \frac{1-\cos\theta\cos\theta''+\sin(\theta/2)\sin(\theta'/2)\sin(\theta''/2)}{2\sin^2(\Theta/2)}, \nonumber \\
\sin^2 L & = & \frac{1-\cos\theta\cos\theta'}{2\sin^2(\Theta/2)}.
\end{eqnarray}

These formulas become symmetric with respect to each of the three successive rotations 
only when, the rotations being infinitely small, the term $\sin(\theta/2)\sin(\theta'/2)\sin(\theta''/2)$ 
vanishes in $\sin^2 H$.  This vanishing also makes the order of these rotations indifferent; one
then finds, in accordance with the law of composition of infinitesimally small rotations, 
\begin{equation}\label{xyzinftes}
\Theta^2 = \theta^2 + \theta'^2 + \theta''^2, \cos G = \theta/\Theta,\cos H = \theta'/\Theta,\cos L = \theta''/\Theta.
\end{equation} 

The problem inverse to the one we have just solved would have for its object the \emph{decomposition} of a 
finite rotation about a given axis into three rotations about the three coordinate axes. This amounts to solving
the above equations for $\theta,\theta',\theta''$ in terms of $\Theta,G,H,L$, which cannot be done for finite
rotations but is utterly simple when the rotations are infinitely small.

\section{On the composition of successive infinitesimal displacements of a solid system.}

We shall now consider, directly and with a special extension, the laws of composition and decomposition 
of successive \emph{infinitesimal} displacements, drawn from the analytic expression for the infinitely small coordinate
changes of a solid system.

These changes are generally expressed in the following way,
\begin{equation}\label{deltainftesorig}
\delta x = \alpha + py - nz,\;\; \delta y = \beta + mz - px,\;\; \delta z = \gamma + nx - my,
\end{equation}

{\it ($\vec{\Delta} = \vec{\delta} + \vec{r}\times\vec{q}$, where $\vec{\Delta}$ and $\vec{q}$ are infinitesimal)}

\noindent as \emph{linear} functions of the infinitesimal elements of the displacement, $\alpha,\beta, \gamma, m,n, p$ 
{\it (that is, $\vec{\delta}$ and $\vec{q}$)}.  It results that the changes arising from several successive infinitesimal
displacements combine by adding together the changes due separately to each of these successive displacements,
referred to the original situation of the system.  The {\bf elements} of the composite displacement are
the sums of the analogous elements of the partial displacements.

This is the {\it (usual)} way in which the complete differential of a function of several variables is formed by 
adding the partial differentials relative to each of those variables.  And since one neglects infinitesimals
of the second order, it makes no difference whether the differentials are expressed in terms of the 
starting values of the finite variables, or in terms of their values successively augmented by the infinitesimal
increments they undergo.

Hence if one denotes by $\vec{\delta}', \vec{q}'; \vec{\delta}'', \vec{q}''; \vec{\delta}''', \vec{q}'''$ etc.
the elements of the successive displacements to be composed, for each of which one has
\begin{equation}\label{delqsucc}
\vec{\Delta}' = \vec{\delta}' + \vec{q}'\times\vec{r}, \;\;\vec{\Delta}'' = \vec{\delta}'' + \vec{q}''\times\vec{r}, etc.,
\end{equation}

\noindent the elements of the composite displacement, represented by $A.B,C,M,N,P$ {\it (that is, by $\vec{\Gamma},\vec{Q}$)},
will be respectively the sums of the given partial elements; one will have
\begin{equation}\label{ABCinftesorig}
A=\alpha' + \alpha'' +... = \Sigma\alpha, B=\Sigma\beta, C=\Sigma\gamma, M=\Sigma m, N=\Sigma n, P=\Sigma p,
\end{equation}

\noindent and for the expressions of the composite {\it (infinitesimal)} variations of the coordinates $x,y,z$,
\begin{equation}\label{deltainftesxyz}
\vec{\Delta} = \Gamma + \vec{r}\times\vec{Q} = \Gamma + \Theta\vec{r}\times\hat{Q}
\end{equation} 

\noindent where we have introduced the  {\it (infinitesimal)} rotation $\Theta$ and the direction $\hat{Q}$.

{\it The author now points out that since the displacements are infinitesimal, one might equally
decompose the same total displacement (resulting from pure rotations about arbitrarily many axes
arbitrarily placed and directed) into exactly \emph{three} partial rotations taken successively
about the $x$,$y$,$z$ axes, with infinitesimal angles $m,n,p$, each accompanied by the appropriate
\lq\lq screw" translation $\vec{t}\cdot\hat{x}, \vec{t}\cdot\hat{y}, \vec{t}\cdot\hat{z}$.}

\section{Geometry and analytical mechanics.}

These successive rotations $m,n,p$ are known in mechanics by the name of {\bf elementary rotations}, and 
considered as {\bf simultaneous} in the passage from the geometric to the mechanical laws of
the displacement of bodies, notwithstanding that geometry cannot take account of them except by supposing them
to be {\bf successive}.  For it is evident that the system, in turning about the axis whose rotation is $\theta$ and 
whose angles with the $x,y,z$ axes are $g,h,l$, does not achieve {\bf at the same time} the three rotations
$m,n,p$ about those coordinate axes: this would require four axes of rotation instead of a single one
(M\'{e}canique Analytique, vol. 1, p. 52).\cite{Lagmech}

 {\it This reference is to the epoch-making two-volume treatise (1788-9)
by Joseph-Louis Lagrange, drawing on discoveries and insights from Euler, D'Alembert, the Bernoulli brothers,
and others, which for the first time showed that all mechanical properties of a system could be derived from a single
algebraic formula and embodied in a single differential equation, without (according to Lagrange himself) requiring
\lq\lq either geometrical or mechanical constructions or reasoning".}

We encounter here a fundamental point in the philosophy of mathematics, that which separates geometry from mechanics,
and the importance of which it is the object of this Memoir to establish in its totality.

The rotation $\theta$ results from the successive composition of the rotations $\theta\cos g,\theta\cos h,\theta\cos l$, because
the displacement due to this rotation $\theta$ is for each coordinate axis the sum of the displacements which would be due
separately to each of the elementary rotations.  In fact, if the system were to turn only about the $x$-axis with a rotation $\theta\cos g$, 
one would have for this displacement 
\begin{equation}\label{delprg}
\delta'x = 0,\;\;\delta'y = \theta z\cos g, \;\;\delta'z =  -\theta y\cos g;
\end{equation} 

\noindent if on the contrary the system were to turn only about the $y$-axis with a rotation $\theta\cos h$, 
one would have for this displacement 
\begin{equation}\label{delprprh}
\delta''x = -\theta z\cos h,\;\;\delta''y = 0, \;\;\delta''z =  \theta x\cos h;
\end{equation} 

\noindent {\it  (the expressions for $\delta''x$ and $\delta''z$ in the French text contain slight errors)}

\noindent and finally the rotation $\theta\cos l$ about the $z$-axis. considered alone, would give
\begin{equation}\label{delprprprl}
\delta'''x = \theta y\cos l,\;\;\delta'''y = -\theta x\cos l, \;\;\delta'''z = 0.
\end{equation} 

The sum of these composite displacements relative to each axis will therefore give for the 
whole displacement, effected definitively by the rotation $\theta$ about the axis $(g,h,l)$,
the following expression:
\begin{equation}\label{delghl}
\delta\vec{r} = \theta \vec{r}\times\hat{t}.
\end{equation} 

\noindent {\it This is the sum, in modern notation, of $\delta'\vec{r} + \delta''\vec{r} + \delta'''\vec{r}$
as given by the three previous equations, where $\hat{t} = (\cos g,\cos h,\cos l)$. It appears that the
author is continuing to make all displacements infinitesimal as in the previous section.}

\section{Successive infinitesimal displacements}

Let us return to the composition of arbitrary displacements, or  changes, given successively for
the same system. The formula for a single such displacement can be written, in conformity with
paragraph {\bf 15}, as
\begin{eqnarray}\label{delrinftes}
\delta x = \alpha + \theta(y\cos l - z\cos h) = t\cos g + \theta u\cos G \nonumber \\
\delta y = \beta + \theta(z\cos g - x\cos l) = t\cos h + \theta u\cos H \nonumber \\
\delta z = \gamma + \theta(x\cos h - y\cos g) = t\cos l + \theta u\cos L.
\end{eqnarray}

\noindent These equations, in which $u$ denotes the distance from the point $(x,y,z)$ to 
the central axis of the displacement, and $(G,H,L)$ are the angles formed with the 
coordinate axes by the direction of the infinitely small arc $u\theta$,  describe a
rotation about that central axis.

{\it In modern form the above array becomes
\begin{equation}\label{delrinftesmod}
\delta\vec{r} = \vec{\delta} + \theta\vec{r}\times\hat{t} = t\hat{t} - \theta \vec{u}\times\hat{t}.
\end{equation}

\noindent We are dealing again with a substitute axis having the same direction $\hat{t}$
as the central axis, but displaced so as to pass through the origin $\vec{r}=0$.  The rotation
about the central axis is equivalent to the same rotation about the substitute axis, augmented
by a translation $\vec{\delta}$ that is common to all points $\vec{r}$.  In the case $\vec{r} = 0$
the substitute rotation has no effect and $\vec{\delta}$ is the whole displacement.  For general
$\vec{r}$ one has also the substitute rotation $\theta\vec{r}\times\hat{t}$. 

The expression after the second = sign is obtained by decomposing $\vec{\delta}$ into a part
parallel and a part perpendicular to the central axis.  The parallel part 
$\vec{\delta}\cdot\hat{t}\hat{t} = t\hat{t}$ is the translation of the central axis along itself, which 
in general accompanies the rotation about the central axis in accordance with the \lq\lq screw"
principle.  The perpendicular part describes the rotation of the origin of coordinates about the
central axis; since $\vec{u}$ is the perpendicular from the point $\vec{r}$ to the central axis,
$\vec{r}+\vec{u}$ is the perpendicular from the origin to the central axis.  Therefore the perpendicular  
part of $\vec{\delta}$ is $-\theta(\vec{r}+\vec{u})\times\hat{t}$. Combining the two parts, we obtain
$\vec{\delta} = t\hat{t} - \theta\vec{r}\times\hat{t} - \theta\vec{u}\times\hat{t}$ or
\begin{equation}\label{deldecomp}
\vec{\delta} + \theta\vec{r}\times\hat{t} = t\hat{t} - \theta \vec{u}\times\hat{t}
\end{equation}

\noindent in keeping with \eqref{delrinftesmod}.}

In general,  the displacement relative to an {\bf arbitrary} direction $s$, 
making angles $a,b,c$ with the $x,y,z$ coordinates, is given by
\begin{equation}\label{dels} 
\delta s = t\cos(t,s) + \theta u\cos(tu,s)
\end{equation}

\noindent where $\cos(t,s)$ , $\cos(tu,s)$ are the cosines of the angles made by
this direction with the central axis and with the infinitely small arc of rotation $\theta$.
[[{\it If we put this equation in modern form, 
\begin{equation}\label{delsmod}
\hat{s}\cdot\delta\vec{r} = t\hat{s}\cdot\hat{t} + \theta u \hat{s}\cdot(-\hat{u}\times\hat{t}),
\end{equation}

\noindent we see that it results from taking the dot product of \eqref{delrinftesmod} with $\hat{s}$.
We can also recover the three equations \eqref{delrinftes} by replacing $s$ with $x$, $y$, or $z$.}

But if two lines are given in the space, one knows that the {\it (shortest)} distance from a {\it (particular)}
point of one line to the other line is reciprocal to the sine of the angle formed by the first line
with the plane containing the second line and the point {\it (of the first line)} under consideration.
Or equivalently, the product of this sine and this distance is constantly equal to the product of
the shortest distance between the two lines and the sine of their inclination {\it (that is, of the angle
between their directions)}.  If, therefore, we define $D$ as the distance between the central axis
and the line passing through the point of the system under consideration
{\it (that is, the point $\vec{r}$)} in the direction $s$, and $\nu$ as the angle between this line 
and the central axis, we shall have
\begin{equation}\label{Dseqs}
u\cos(tu,s) = D\sin\nu,\;\;\;\delta s = t\cos\nu + D\theta\sin\nu.
\end{equation}

\noindent {\it Here the author liberates himself from the given point $\vec{r}$ and refers in his
description to whatever point of the $s$-line through $\vec{r}$ is closest to the central
axis. This completes his analysis of a {\bf single}  displacement.}

If we now consider the successive {\it (infinitesimal)} displacements of the system about 
central axes whose elements are $t,\theta, g, h,l$;\;
$t',\theta', g', h',l'$;\; $t'',\theta'', g'', h'',l''$; etc., and denote by $\delta S$ the resultant
displacement of a point of the system relative to that same fixed direction $s$ and
by $T, \Theta, G, H, L$ the elements of the resultant central axis,  we shall have
\begin{equation}\label{delS}
\delta S = \Sigma t\cos\nu + \Sigma D\theta\sin\nu  = T\cos V + \Theta{\cal D}\sin V,
\end{equation}

\noindent ${\cal D}$ being the distance from the resultant axis to the line $s$ drawn
through the point $(x,y,z)$.

This equation, owing to the indeterminate parameters $a,b,c$ implicit in it and to
the fact that it must hold for all points of the system, is equivalent to the following
six equations which give the position of the resultant central axis and the resultant
translation and rotation:

\begin{eqnarray}\label{compinftes}
\Sigma\alpha - T\cos G + \Theta Y\cos L - \Theta Z\cos H= 0, \nonumber \\
\Sigma\beta - T\cos H + \Theta Z\cos L - \Theta X\cos L= 0, \nonumber \\
\Sigma\gamma - T\cos L + \Theta X\cos L - \Theta Y\cos G= 0, \nonumber \\
\Theta\cos G = \Sigma\theta\cos g,\Theta\cos H = \Sigma\theta\cos h,\Theta\cos L = \Sigma\theta\cos l,
\end{eqnarray} 

\noindent where $X,Y,Z$ are the coordinates of an arbitrary point on the resultant central axis.

Without loss of rigor, one could consider only pure rotations in this analysis, since the translations
$t,t',...$ can always be represented by couples of rotations; and simply put
\begin{equation}\label{delSrot}
\delta S= T\cos V + \Theta{\cal D}\sin V = \Sigma\theta D\sin\nu.
\end{equation}

{\it Here I omit some commentary that seems to me both tedious and repetitious.}

\section{Conditions for {\bf equilibrium} from many infinitely small successive displacements.}

We are now led to seek out what conditions need to be satisfied by the elements of the
successive displacements proposed for a solid system, in order that the system, passing
successively through various infinitesimally neighboring situations, should return to its
initial position; this would amount to a condition of {\bf equilibrium}, or neutralization,
on the totality of the successive displacements.  Now, it is evident that all the relevant
conditions are contained in a single equation, which can be decomposed into six others
on account of the indeterminate quantities implicit in it, to wit:
\begin{equation}\label{delS0}
\delta S = 0,
\end{equation}

\noindent since this equation expresses that each point of the system has returned
to its initial position.

The six equations hidden in $\delta S = 0$ are
\begin{eqnarray}\label{6S0}
\delta x_0 & = & \Sigma \alpha = 0, \nonumber \\
\delta y_0  & = & \Sigma \beta = 0, \nonumber \\
\delta z_0  & = & \Sigma \gamma = 0, \nonumber \\
\Sigma\theta\cos g = 0,\;\; \Sigma\theta\cos h & = & 0,\;\; \Sigma\theta\cos l = 0,
\end{eqnarray} 

\noindent where $\delta x_0,\delta y_0,\delta z_0$ stand for the resultant
variations of the coordinates of the origin.  The first three equations express
the {\bf immobility} of the origin of the coordinates; and the the three others,
that no resultant rotation has occurred in the displaced system {\it (i.e. no 
rotation about a resultant fixed axis through the origin, such as Euler's theorem
permits)}.  This double condition excludes the possibility of any resultant
displacement whatever. 

The double condition is an immediate consequence of the relations \eqref{delSrot} 
and\eqref{delS0}, which cannot be satisfied for every point in the system unless
\begin{equation}\label{TTh0}
T = 0, \Theta = 0.
\end{equation}

\noindent For $\Theta = 0$ implies the three last equations of\eqref{6S0}, while
the first three follow from \eqref{Delsq} of Section {\bf 15}.

These six equations of equilibrium are analytically contained in a single
equation which expresses the \emph{general} law of this equilibrium in
the simplest way {\it (in terms of the elements of the successive displacements)}, namely
\begin{equation}\label{SumthDnu0}
\Sigma\theta D\sin\nu = 0.
\end{equation}

\noindent But the equilibrium of these infinitesimal {\bf successive} displacements,
all else remaining the same, will continue to hold no matter how rapidly they succeed
one another.  Passing to the limit, one arrives at the {\bf identity} of the laws
of equilibrium due to successive infinitesimal displacements with those due to
{\bf simultaneous} infinitesimal displacements.  {\it This equivalence is needed
in order to justify the analogy between geometry and mechanics, to be presented
in the following section.}

\section{Analogy of these {\it (geometric)} laws of composition and equilibrium with those
of composition and equlibrium of {\bf forces} applied to an immovable system.}

The analogy between this general {\it (geometric)} law and that of the equilibrium of forces
applied to an immovable system is striking.  Let the applied forces follow the axes of
rotation and suppose them proportional to those rotations; then the {\bf moment} of 
a force on the system is exactly proportional to that of the corresponding rotation, and
each translation is replaced by a couple {\it (in the sense of earlier sections)} of applied forces
following the axes of the couple of rotations equivalent to the translation.  The analogy,
however, extends to the laws of composition and may be stated thus:

A system of successive displacements being given to be composed into a resultant displacement,
and at the same time  a system of forces proportional to the successive rotations given
for each displacement and applied following the same axes as the rotations, the translations 
of the successive displacements, supposing that they are not implicitly included in the rotations
by being represented by couples of rotations, being then represented in the system of
forces under consideration by couples of forces whose moments would be equal to
those of the translations, relatively to the three coordinate axes, the system of displacements
will amount to a resultant displacement composed of a rotation and an absolute translation
relative to the central axis of rotation; just as the system of forces will resolve itself by
the successive composition of its elements into a single force and a single couple of forces
situated in a plane normal to the resultant force.  This resultant force will be applied at the
central axis of the resultant displacement, which will be at the same time the central axis
of the static system; it will be proportional to the resultant rotation, and the moment of the couple
normal to this force will be proportional to the absolute translation of the system that operates
in a way parallel to the central axis.  Should the axes of the rotations to be combined be all
parallel and pass through determined points, the resultant central axis is also parallel to them
and passes through a certain point that corresponds to the center of the parallel forces,
the same point no matter what be the direction of the axes of rotation, and which is nothing 
other than the center of gravity of the points on the composing axes of rotation that are
determined when all the rotations are equal.

{\it This passage, as translated faithfully above, is somewhat in need of a retranslation or decipherment.
The basic key to the passage is that the phrase \lq\lq moment d'une force suivant un axe"
means, at least in modern French, the {\bf torque} $\vec{r}\times\vec{F}$ about that axis, 
supposing that the force $\vec{F}$ is applied to or through a point whose position vector 
relative to some point on the axis is $\vec{r}$. Thus, in general, the force is not applied along the axis, 
nor even toward some point on the axis; if it were, the \lq\lq moment" would be zero. Of course the 
words \lq\lq torque" and \lq\lq moment" were not as cleanly defined in the author's time as now,
but this identification does clarify the intention of the passage.

With this key in hand, we can appreciate what otherwise would be perplexing: that a {\bf single force}
is repeatedly associated with a geometric {\bf rotation}, while  a geometric {\bf translation} corresponds
to the moment of a {\bf pair or couple} of forces acting in a plane normal to the geometric translation.
As to the geometric meaning of the word \lq\lq moment", it may be helpful to refer to the discussion
following \eqref{Delsq} in {\bf 15}.  Of course the modern form of a cross-product using the right-hand
rule was not available; instead the author measures the \lq\lq moment" of a rotation by an area,
an idea introduced independently by H. Grassmann almost at the same time.}

In fact, the equations of the resultant central axis of the composition of the fixed axes of rotation
reduce in this case to
\begin{equation}\label{centaxcomp0} 
\vec{W}\times\hat{T} + \frac{\Sigma\vec{\delta}}{\Sigma\vec{\theta}} = 0,
\end{equation}

\noindent and we have
\begin{eqnarray}\label{TThetc0}
T & = & 0, \;\;\Theta = \Sigma\theta \nonumber \\
= \Sigma\alpha & = & \Sigma\theta(Z\cos h - Y\cos h), \nonumber \\
= \Sigma\beta & = & \Sigma\theta(X\cos h - Z\cos h), \nonumber \\
= \Sigma\gamma & = & \Sigma\theta(Y\cos h -X\cos h)
\end{eqnarray}

\noindent where $X,Y,Z$ denote the coordinates of the axis of rotation $\theta$.
Hence the resultant axis passes through the point whose coordinates are
\begin{equation}\label{xyzSig}
x = \frac{\Sigma\theta X}{\Sigma\theta},\;\;y = \frac{\Sigma\theta Y}{\Sigma\theta},\;\;z = \frac{\Sigma\theta Z}{\Sigma\theta}.
\end{equation}

\section{Determination of the changes in coordinates of a solid due to an arbitrary displacement, deduced 
analytically from the conditions of invariability of the system.}

We consider the displacement of the coordinate axes, given that they are rigidly attached to
the system as it is displaced.  This leads immediately to the algebraic expression for the
changes in the coordinates of any point, and there remains only to reduce to a minimum
the number of arbitrary constants that enter the calculation, as we shall now see.

{\it The author now proposes to derive the preceding results {\bf without} using the formulas 
derived from \eqref{Deldelmodq}, but by a different method.}

Let us denote by $a,b,c;\;a',b',c';\;a'',b'',c''$ the cosines of the angles formed by the displaced
coordinate axes with their original directions.  Then the {\bf new}  coordinates $x+\Delta x,
y+\Delta y, z+\Delta z$, of a point of the system after displacement relative to the \lq\lq old" axes, 
will be expressed as a function of the same displaced coordinates relative to the {\bf new}
axes. But this will be the same function, already well known, that gives the new {\bf axes} 
in terms of the old, that is
\begin{eqnarray}\label{chgaxes}
x+ \Delta x = \alpha + ax + by + cz, \nonumber \\
y+ \Delta y = \beta + a'x + b'y + c'z, \nonumber \\
z+ \Delta z = \gamma + a''x + b''y + c''z.
\end{eqnarray}

{\it The idea is that since the axes are rigidly attached to the system, they suffer the same
displacement from old to new as does the point under consideration. So if one compares
the final position of the point to its old position, both relative to the old axes, one finds the
same transformation formulas as in comparing the coordinates of the same {\bf new}
point, relative to the old axes, with its coordinates relative to the new ones.}

These formulas express $x+\Delta x, y+\Delta y, z+\Delta z$ as linear functions of the
initial coordinates $x,y,z$ for any displacement of the solid whatever.  They contain
12 arbitrary coefficients, but really only 6 of these are independent.  The first three,
$\alpha, \beta, \gamma$, determine the {\bf movement} of the origin, and the three
\lq \lq diagonal" coefficients $a, b', c''$ serve to define the {\bf directions} of the
three displaced axes relative to the old ones.  In terms of the latter, one may eliminate
the six off-diagonal elements by the formulas of Monge.

But the reduction of twelve constants to six can be achieved even more simply, without
using the formulas of Monge (complicated by radicals) and by a route that leads to
the simplest possible expressions for $\Delta x,\Delta y,\Delta z$ in terms of $x,y,z$.
We introduce into \eqref{chgaxes} the coordinates $\xi,\eta,\zeta$ of the midpoint of
the line joining the initial to the final position of the point in question. This gives us
\begin{equation}\label{rfromom}
x = \xi - (1/2)\Delta x, \;\; y = \eta - (1/2)\Delta y, \;\;z = \zeta - (1/2)\Delta z.
\end{equation}

Then we unite these three equations into one by multiplying them respectively 
by three indeterminate factors $\mu,\nu,\pi$ and adding together the results.  We now have 
the single equation
\begin{eqnarray}\label{xyzfrommunupi}
(1/2)\Delta x (\bar{a}+\mu) & + & (1/2)\Delta y (\bar{b}+\nu) + (1/2)\Delta z (\bar{c}+\pi) = (\alpha\mu + \beta\nu + \gamma\pi) \nonumber \\
& + & \xi(\bar{a}-\mu) + \eta(\bar{b}-\nu) + \zeta(\bar{c}-\pi),
\end{eqnarray}

\noindent where we define
\begin{equation}\label{barabcdef}
\bar{a} = a\mu + a'\nu + a''\pi, \;\bar{b} = b\mu + b'\nu + b''\pi, \;\bar{c} = c\mu + c'\nu + c''\pi.
\end{equation}

To determine {\it (for example)} $\Delta x$, we must assign to $\mu,\nu,\pi$ those values 
that cause the coefficients of $\Delta y$ and $\Delta z$ to vanish; that is, we must set
\begin{equation}\label{yzcoefs0}
\bar{b} + \nu = \bar{c}+ \pi = 0, 
\end{equation}

\noindent from which we find
\begin{equation}\label{munupix}
\mu = (1+b')(1+c'') -c'b'', \nu = b''c - b(1+c''),\pi = bc' - c(1+b').
\end{equation}

\noindent {\it (French text has incorrectly $a''$ instead of $b''c$ in $\nu$)}.

Now, the nine cosines $\vec{a}=(a,a',a''),\vec{b}=(b,b',b''),\vec{c}=(c,c',c'')$ satisfy, as 
is well known, the relations
\begin{eqnarray}\label{abccross}
\vec{a} = \vec{c}\times\vec{b}, \nonumber \\
\vec{b} = \vec{a}\times\vec{c}, \nonumber \\
\vec{c} = \vec{b}\times\vec{a}.
\end{eqnarray}

\noindent These can be deduced from six others,
\begin{eqnarray}\label{orthonorm}
|\vec{a}|^2  = |\vec{b}|^2 = |\vec{c}|^2 =1, \nonumber \\
\vec{a}\cdot\vec{b} = \vec{b}\cdot\vec{c} = \vec{c}\cdot\vec{a} = 0,
\end{eqnarray}

\noindent which say that the coordinate axes, both the old and the new,
form orthonormal systems.  The ambiguity of signs that enter into the
deduction is resolved by another condition altogether necessary in
considering the displacement of a solid system, namely that the system
of new axes arising from the displacement must always remain
in the condition of superposition with the old axes that is possible for
each respective axis and its correspondent.

{\it I have adhered to my usual convention of reversing cross-products
so as to maintain right-handedness (see {\bf 15}). I hope that by
doing so I satisfy the author's condition above.}

This being understood, it is clear that the above relations lead to
\begin{equation}\label{munupisimp}
\mu = 1 + a+ b' + c'', \nu = a' - b, \pi = a''-c,
\end{equation}

\noindent from which
\begin{equation}\label{baramu}
\bar{a} = 1 + a + b' + c'' = \mu
\end{equation}

\noindent and hence
\begin{equation}\label{Delxfin}
\Delta x - \alpha = \frac{2[(\eta-(1/2)\beta)(b-a') - (\zeta-(1/2)\gamma)(a'' - c)]}{1+a+b'+c''}. 
\end{equation}

\noindent {\it In eqs \eqref{yzcoefs0}, \eqref{munupix},\eqref{munupisimp},\eqref{baramu}, and \eqref{Delxfin} ,
the author has broken the threefold symmetry of \eqref{xyzfrommunupi} by singling out $\Delta x$ as the 
component to be determined.  He now retreats from that choice by considering in turn each of
the other components $\Delta y,\Delta z$.}

Similarly one can obtain
\begin{equation}\label{Delyfin}
\Delta y - \beta = \frac{2[(\zeta-(1/2)\gamma)(c'-b'') - (\xi-(1/2)\alpha)(b - a')]}{1+a+b'+c''}. 
\end{equation}

\begin{equation}\label{Delzfin}
\Delta z - \gamma = \frac{2[(\xi-(1/2)\alpha)(a''-c) - (\eta-(1/2)\beta)(c' - b'')]}{1+a+b'+c''}. 
\end{equation}

These formulas are identical to those of {\bf 15}], provided that we set
\begin{equation}\label{mnp15-19}
m = \frac{2(c'-b'')}{1+a+b'+c''},\;\;n = \frac{2(a''-c)}{1+a+b'+c''},\;\; p = \frac{2(b-a')}{1+a+b'+c''}.
\end{equation}

\section{Infinitesimal version of {\bf 27}.}

In the case of infinitesimal displacements, we neglect angular displacements of second order
and directly obtain $a = b'= c'' =1$, from which
\begin{equation}\label{delxyz}
\delta x = \alpha+by+cz,\;\; \delta y = \beta+c'z+a'x,\;\; \delta z = \gamma + a''x + b''y.
\end{equation}

\noindent Also, the distance from any point to the origin, when the latter is drawn along with
the displacement of the system, is invariable, so that {\it (vectorially)} $\vec{r}\cdot(\delta\vec{r}-\vec{\delta}) = 0$
for any arbitrary point $\vec{r}$.  This necessitates the following relations among the 
{\it\lq\lq off-diagonal"} first-order infinitesimal cosines:
\begin{equation}\label{antisym}
b+a' = c+ a'' = c' + b'' = 0.
\end{equation}

\noindent Hence by setting $m = c' = -b'', n = a'' = -c, p = b = - a'$, we have finally
\begin{equation}\label{delxyzfin}
\delta x= \alpha +  py - nz, \;\; \delta y = \beta +  mz - px, \;\;\  \delta z = \gamma +  nx - my.
\end{equation}

\section{Algebraic deduction of coordinate changes.}

{\it And now, yet another method.}

But it is interesting to arrive at these same formulas for finite or infinitesimal coordinate changes
by a purely algebraic route, independent of any geometric consideration {\it (other than the Pythagorean
formula for the distance between two points)} starting from the invariability of the distances 
between points of the solid.

So let
\begin{equation}\label{four points}
\vec{r}_0,\vec{r}_1,\vec{r}_2,\vec{r}
\end{equation}

\noindent be the positions of four points invariably linked together and belonging to the solid system.
The first three points and their displacements are to be considered as known, and the displacement
of the fourth point is to be calculated as a function of its starting position and the known quantities.

The distances between these four points will remain constant under an arbitrary displacement of the
system; this condition, when expressed algebraically, will give the following six equations:
\begin{eqnarray}\label{tetrig}
|\vec{r}_1+\Delta\vec{r}_1 -\vec{r}_0-\Delta\vec{r}_0 |^2  =  |\vec{r}_1 - \vec{r}_0| ^2 \nonumber \\
|\vec{r}_2+\Delta\vec{r}_2- \vec{r}_0-\Delta\vec{r}_0|^2  =  |\vec{r}_2 - \vec{r}_0|^2 \nonumber \\
|\vec{r}_2+\Delta\vec{r}_2 - \vec{r}_1-\Delta\vec{r}_1|^2  =  |\vec{r}_2 - \vec{r}_1| ^2 \nonumber \\
|\vec{r}+\Delta\vec{r} - \vec{r}_0-\Delta\vec{r}_0|^2  =  |\vec{r} - \vec{r}_0| ^2 \nonumber \\
|\vec{r}+\Delta\vec{r} - \vec{r}_1-\Delta\vec{r}_1|^2  =  |\vec{r} - \vec{r}_1| ^2 \nonumber \\
|\vec{r}+\Delta\vec{r} - \vec{r}_2-\Delta\vec{r}_2|^2  =  |\vec{r} - \vec{r}_2)| ^2. 
\end{eqnarray}

{\it I am rendering equations when possible in vector form for brevity and ease of reading.  The Pythagorean formula 
for distance is implicit in the vectors; thus $|\vec{r}-\vec{r}_0|^2= (x-x_0)^2+(y-y_0)^2+(z-z_0)^2$, etc.

We may think of a tetrahedron with vertex at $\vec{r}$ and base having the three points $\vec{r}_0,\vec{r}_1,\vec{r}_2.$
Then the first three equations of \eqref{tetrig} pertain to the edges of the base triangle, and the last three to
the edges that meet at the vertex.}

The six equations are quadratic in the  displacements $\Delta\vec{r_0}$,$\Delta\vec{r}_1$,$\Delta\vec{r}_2$,$\Delta\vec{r}$ 
but can be made linear by introducing the midpoint position $\vec{\omega} = \vec{r} + (1/2)\Delta\vec{r}$, and likewise for
the three base points. One then obtains
\begin{eqnarray}\label{tetriglin}
2(\vec{\omega}-\vec{\omega}_0)\cdot(\Delta\vec{r}-\Delta\vec{r}_0)= 0, \nonumber \\
2(\vec{\omega}_1-\vec{\omega}_0)\cdot(\Delta\vec{r}_1-\Delta\vec{r}_0) = 0, \nonumber \\
2(\vec{\omega}_2-\vec{\omega}_0)\cdot (\Delta\vec{r}_2-\Delta\vec{r}_0)= 0, \nonumber \\
(\vec{\omega}-\vec{\omega}_0)\cdot(\Delta\vec{r}_1-\Delta\vec{r}_0)+(\Delta\vec{r}-\Delta\vec{r}_0)\cdot(\vec{\omega}_1-\vec{\omega}_0) = 0, \nonumber \\
(\vec{\omega}-\vec{\omega}_0)\cdot(\Delta\vec{r}_2-\Delta\vec{r}_0)+(\Delta\vec{r}-\Delta\vec{r}_0)\cdot(\vec{\omega}_2-\vec{\omega}_0) = 0, \nonumber \\
(\vec{\omega}_1-\vec{\omega}_0)\cdot(\Delta\vec{r}_2-\Delta\vec{r}_0)+(\Delta\vec{r}_1-\Delta\vec{r}_0)\cdot(\vec{\omega}_2-\vec{\omega}_0) = 0.
\end{eqnarray}

{\it These equations may be derived as follows.  We note that each linear factor in \eqref{tetriglin} contains either $\vec{\omega}_0$ or $\Delta\vec{r}_0$.
Let $(\vec{R},\vec{\Omega})$ and $(\vec{R}',\vec{\Omega}')$ be any of the pairs 
$(\vec{r},\vec{\omega})$, $(\vec{r}_1,\vec{\omega}_1)$, $(\vec{r}_2,\vec{\omega}_2)$.  From the definition $\vec{\Omega} = \vec{R} + (1/2)\Delta\vec{R}$
we have 
\begin{equation}\label{RR+DfromOm}
\vec{R} = \vec{\Omega} - (1/2)\Delta\vec{R}, \;\;\vec{R} + \Delta\vec{R} = \vec{\Omega} + (1/2)\Delta\vec{R},
\end{equation}

\noindent and likewise with $\vec{R}',\vec{\Omega}'$ in place of $\vec{R},\vec{\Omega}$.  Therefore 
\begin{eqnarray}\label{tetlindiffsq}
& & [(\vec{R}+\Delta\vec{R})-(\vec{r}_0+\Delta\vec{r}_0)]\cdot[(\vec{R}'+\Delta\vec{R}')-(\vec{r}_0+\Delta\vec{r}_0)]  - (\vec{R}-\vec{r}_0)\cdot(\vec{R}'-\vec{r}_0) \nonumber \\
& = & [(\vec{\Omega}+(1/2)\Delta\vec{R}) - (\vec{\omega}_0+(1/2)\Delta\vec{r}_0)]\cdot[(\vec{\Omega}'+(1/2)\Delta\vec{R}') - (\vec{\omega}_0+(1/2)\Delta\vec{r}_0)]  \nonumber \\
& - &[(\vec{\Omega}-(1/2)\Delta\vec{R}) - (\vec{\omega}_0-(1/2)\Delta\vec{r}_0)]\cdot[(\vec{\Omega}'-(1/2)\Delta\vec{R}') - (\vec{\omega}_0-(1/2)\Delta\vec{r}_0)] \nonumber \\
& = & (\vec{\Omega} - \vec{\omega}_0)\cdot(\Delta\vec{R}' - \Delta\vec{r}_0) + (\Delta\vec{R} - \Delta\vec{r}_0)\cdot(\vec{\Omega}' - \vec{\omega}_0)
\end{eqnarray}

\noindent and since the top line of \eqref{tetlindiffsq} vanishes by \eqref{tetrig}, the bottom line vanishes, giving us the first three lines of \eqref{tetriglin} if $(\vec{R},\vec{\Omega})$
and $(\vec{R}',\vec{\Omega}')$ are the same, and the last three if they are different.

It will be observed that the author has selected the pair $(\vec{r}_0,\vec{\omega}_0)$ to play the starring r\^{o}le in \eqref{tetriglin}, inasmuch as it is the only pair that
is represented in every linear factor in every equation.  Other equations could be written down, but they can be deduced from the six displayed.  For example, the 
equation $(\vec{\omega}_0-\vec{\omega})\cdot(\Delta\vec{r}_1-\Delta\vec{r})+(\Delta\vec{r}_0-\Delta\vec{r})\cdot(\vec{\omega}_1-\vec{\omega}) = 0$, obtained by
interchanging $\vec{r}, \vec{\omega}$ with  $\vec{r}_0, \vec{\omega}_0$ in the fourth line of \eqref{tetriglin}, can be deduced by subtracting that line from the top line.}

We now multiply these six equations by the respective six factors $\mu^2, \nu^2, \pi^2, \mu\nu, \mu\pi, \nu\pi$ 
and add them together.  This results in the following single equation which {\it (when $\mu,\nu,\pi$ vary independently)} contains
the foregoing six:
\begin{equation}\label{tetriglinmunupi}
[\mu(\vec{\omega}-\vec{\omega}_0)+\nu(\vec{\omega_1}-\vec{\omega}_0)+\pi(\vec{\omega_2}-\vec{\omega}_0)]
\cdot[\mu(\Delta\vec{r}-\Delta\vec{r}_0)+\nu(\Delta\vec{r}_1-\Delta\vec{r}_0)+\pi(\Delta\vec{r}_2-\Delta\vec{r}_0)] = 0.
\end{equation}

{\it When a similar operation was performed in {\bf 27}, the coefficients $\mu,\nu,\pi$ were attached to the $x,y,z$-components
and the resulting manipulations could not be written vectorially.  Here, the three coefficients are attached to three different vectors.}

By a suitable choice of $\mu,\nu,\pi$, this equation can be made to contain only $\Delta\vec{r}$ [[and not $\Delta\vec{r_1}$ or $\Delta\vec{r_2}$]]
and the value of this displacement can be found in the simplest way. 

{\it Having thus far allowed us the luxury of our vector notation, the author is now going to forbid it by singling
out one component $x$ for study.  What makes this necessary is that he wishes to use the principle
that the vanishing of a product of two scalar quantities implies the vanishing of at least one of the factors,
and this is not true of the dot product of two vectors.}

To wit, suppose we set 
\begin{eqnarray}\label{yz0ofomom1om2}
\mu(\eta-\eta_0) + \nu(\eta_1-\eta_0) + \pi(\eta_2-\eta_0) = 0, \nonumber \\
\mu(\zeta-\zeta_0) + \nu(\zeta_1-\zeta_0) + \pi(\zeta_2-\zeta_0) = 0:
\end{eqnarray}

\noindent then \eqref{tetriglinmunupi} reduces to
\begin{equation}\label{xtetriglinmunupi}
[\mu(\xi-\xi_0)+\nu(\xi_1-\xi_0)+\pi(\xi_2-\xi_0)][\mu(\Delta x-\Delta x_0)+\nu(\Delta x_1-\Delta x_0)+\pi(\Delta x_2-\Delta x_0)] = 0.
\end{equation}

\noindent {\it He has arranged $\mu, \nu, \pi$ to make the first factor of the $y$ and $z$ parts of \eqref{tetriglinmunupi} vanish, and so the corresponding products 
vanish; therefore the product in the $x$ part must also vanish.  The principle of factorization then dictates that at least one of the $x$ factors vanish,
but he claims (deferring proof to the following section) that the first factor cannot vanish if the second does not, and so infers that the second factor
vanishes:}

\noindent Now, the first factor of \eqref{xtetriglinmunupi} cannot be zero unless the second is also, as we shall demonstrate hereafter; we therefore have
the simultaneous equations
\begin{eqnarray}\label{Delxetazeta}
\mu(\Delta x - \Delta x_0) + \nu(\Delta x_1- \Delta x_0) + \pi(\Delta x_2 -\Delta x_0) = 0, \nonumber \\
\mu(\eta - \eta_0) + \nu(\eta_1 - \eta_0) + \pi(\eta_2 - \eta_0) = 0, \nonumber \\
\mu(\zeta - \zeta_0) + \nu(\zeta_1 - \zeta_0) + \pi(\zeta_2 - \zeta_0) = 0.
\end{eqnarray}

\noindent These equations evidently imply the following ones:
\begin{eqnarray}\label{Delxnp}
\Delta x_1 - \Delta x_0 = p(\eta_1 - \eta_0) - n(\zeta_1 - \zeta_0), \nonumber \\
\Delta x_2 - \Delta x_0 = p(\eta_2 - \eta_0) - n(\zeta_2 - \zeta_0), \nonumber \\
\Delta x - \Delta x_0 = p(\eta - \eta_0) - n(\zeta - \zeta_0),
\end{eqnarray}  

\noindent  {\it (French text has mistakenly $p(\xi_1-\xi_0)$
in third equation)} where $n$ and $p$ are two constants bound to the displacements of the
three first points {\it ($\vec{r}_1,\vec{r}_2,\vec{r}_0$, the base of the tetrahedron)}
by the first two of these three equations. 

{\it The three simultaneous equations \eqref{Delxetazeta} are homogeneous in $\mu,\nu, 
\pi$; therefore the determinant must vanish.  Interchanging rows and columns, we deduce
that the three homogeneous equations
\begin{eqnarray}\label{mbarnbarpbar}
(\Delta x - \Delta x_0)\bar{m} + (\eta - \eta_0)\bar{n} + (\zeta - \zeta_0)\bar{p} = 0, \nonumber \\
(\Delta x_1 - \Delta x_0)\bar{m} + (\eta_1 - \eta_0)\bar{n} + (\zeta - \zeta_0)\bar{p} = 0, \nonumber \\
(\Delta x_2 - \Delta x_0)\bar{m} + (\eta_2 - \eta_0)\bar{n} + (\zeta_2 - \zeta_0)\bar{p} = 0
\end{eqnarray}

\noindent have a solution $\bar{m},\bar{n},\bar{p}$.  Defining $n=\bar{p}/\bar{m}$, $p = -\bar{n}/\bar{m}$, 
we obtain \eqref{Delxnp}.}

The same analysis {\it (if $y$ or $z$ instead of $x$ had been singled out)} would give
\begin{eqnarray}\label{Delympznm}
\Delta y - \Delta y_0 = m'(\zeta - \zeta_0) - p'(\xi - \xi_0), \nonumber \\
\Delta z - \Delta z_0 = n''(\xi - \xi_0) - m''(\eta - \eta_0).
\end{eqnarray}

\noindent [[French text has $n',m'$ in third line]] But we have also {\it (from the top line of \eqref{tetriglin})}
\begin{equation}\label{tetriglintop}
(\xi-\xi_0)(\Delta x-\Delta x_0) + (\eta-\eta_0)(\Delta y-\Delta y_0) + (\zeta-\zeta_0)(\Delta z-\Delta z_0) = 0;
\end{equation}

\noindent it follows that $m=m'=m''$, $n=n'=n''$, $p=p'=p''$, and so we have at last
\begin{eqnarray}\label{Delxyz}
\Delta x-\Delta x_0 = A + p\eta - n\zeta, \nonumber \\
\Delta y-\Delta y_0 = B + m\zeta - p\xi, \nonumber \\
\Delta z-\Delta z_0 = C + n\xi -m\eta,
\end{eqnarray}

\noindent where the six constants $A,B,C,m,n,p$ are functions of the displacements $\Delta\vec{r}_0,\Delta\vec{r}_1,\Delta\vec{r}_2$.

{\it Substituting \eqref{Delympznm} and the third line of \eqref{Delxnp} into \eqref{tetriglintop}, we find for example that the product
$(\xi-\xi_0)(\eta-\eta_0)$ appears with coefficient $p-p'$ so that $p' = p$.  This is one of six equations that together justify
dropping all the primes in \eqref{Delxyz}.}
\section{Proof of prior claim.}

But it remains to prove what we have claimed, that in  \eqref{xtetriglinmunupi} the first factor
$[\mu(\xi-\xi_0)+\nu(\xi_1-\xi_0)+\pi(\xi_2-\xi_0)]$ cannot vanish unless the second factor
$[\mu(\Delta x-\Delta x_0)+\nu(\Delta x_1-\Delta x_0)+\pi(\Delta x_2-\Delta x_0)]$ does so as well.

The three simultaneous equations {\it (they are easily expressed as one, in vector notation)}
\begin{equation}\label{tetmunupi}
\mu(\vec{\omega} - \vec{\omega}_0) + \nu(\vec{\omega}_1 - \vec{\omega}_0) + \pi(\vec{\omega}_2 - \vec{\omega}_0) = 0
\end{equation}

\noindent express that the midpoints of the lines traversed by the four points we are considering
lie in the same plane. But this unusual condition can be fulfilled only if either, on the one hand,
the four points themselves also lie in one plane, or on the other hand, the pyramid formed by
these points after displacement is not genuinely superposable on the one formed initially, but
is only symmetric to it.  The second hypothesis is not admissible in our problem, but the first
must be examined.

{\it Consider the phrase \lq\lq midpoint plane" as referring to the plane containing the four midpoints
$\vec{\omega}$, etc. If we examine the four points of the tetrahedron in relation to the 
midpoint plane, we see that the only way to make the actual midpoints lie in this plane is
to make the displacement take each vertex of the tetrahedron to the opposite side of the plane.  But the resulting pyramid
is not \lq\lq superposable" on the original one (see {\bf1}, second paragraph) unless all four
points already lie in the midpoint plane, so that the tetrahedron is identical to its mirror image.
This exceptional case is precisely the \lq\lq first hypothesis" advanced by the author.}

We shall now prove that if the four points of the tetrahedron, as well as those to which they are
displaced, all lie in one plane, then both factors of \eqref{xtetriglinmunupi} vanish.  

Let $(x,y,z)$, $(x_1,y_1,z_1)$, $(x_2,y_2,z_2)$, $(x_0,y_0,z_0)$ denote the initial positions of the four 
points, and $(x',y',z')$, etc. their positions after displacement.  If the four points are initially in 
the same plane, there must be coefficients $g,h,l$ such that
\begin{eqnarray}\label{ghl4pts}
g(x-x_0)+h(x_1-x_0)+l(x_2-x_0)=0,\nonumber\\
g(y-y_0)+h(y_1-y_0)+l(y_2-y_0)=0,\nonumber\\
g(z-z_0)+h(z_1-z_0)+l(z_2-z_0)=0.
\end{eqnarray}

\noindent Likewise, if the four points after displacement are in the same plane, there must be
$g',h',l'$ such that

\begin{eqnarray}\label{ghlpr4pts}
g'(x'-x'_0)+h'(x'_1-x'_0)+l'(x'_2-x'_0)=0,\nonumber\\
g'(y'-y'_0)+h'(y'_1-y'_0)+l'(y'_2-y'_0)=0,\nonumber\\
g'(z'-z'_0)+h'(z'_1-z'_0)+l'(z'_2-z'_0)=0.
\end{eqnarray}

\noindent Now, from the invariability of distances we must have
\begin{eqnarray}\label{rr0sq}
|\vec{r}-\vec{r}_0|^2 = |\vec{\,r}'-\vec{\,r}'_0|^2, \nonumber  \\
|\vec{r}_1-\vec{r}_0|^2 = |,\vec{\,r}'_1-\vec{\,r}'_0|^2, \nonumber  \\
|\vec{r}_2-\vec{r}_0|^2 = |\vec{\,r}'_2-\vec{\,r}'_0|^2
\end{eqnarray}

\noindent as well as
\begin{eqnarray}\label{rr0sq}
(\vec{r}-\vec{r}_0)\cdot(\vec{r}_1-\vec{r}_0)= (\vec{\,r}'-\vec{r}'_0)\cdot(\vec{\,r}'_1-\vec{\,r}'_0), \nonumber  \\
(\vec{r}-\vec{r}_0)\cdot(\vec{r}_2-\vec{r}_0) = (\vec{\,r}'-\vec{r}'_0)\cdot(\vec{\,r}'_2-\vec{\,r}'_0), \nonumber  \\
(\vec{r}_1-\vec{r}_0)\cdot(\vec{r}_2-\vec{r}_0) = (\vec{\,r}_1'-\vec{r}'_0)\cdot(\vec{\,r}'_2-\vec{\,r}'_0).
\end{eqnarray}

\noindent But the coefficients $g',h',l'$ satisfy the same relations \eqref{ghlpr4pts} to 
$\vec{\,r}'-\vec{\,r}'_0,\vec{\,r}'_1-\vec{r}'_0, \vec{\,r}'_2-\vec{r}'_0$ as do the coefficients
$g,h,l$ (see \eqref{ghl4pts}) to $\vec{\,r}-\vec{\,r}_0,\vec{\,r}_1-\vec{r}_0, \vec{\,r}_2-\vec{r}_0$.
Consequently $g',h',l'$ are proportional to $g,h,l$ and we may as well set the former equal to
the latter.

{\it In vector notation, \eqref{ghl4pts} and \eqref{ghlpr4pts} become
\begin{eqnarray}\label{ghl4ptsvec}
g(\vec{r} - \vec{r}_0) +  h(\vec{r} - \vec{r}_0) +  l(\vec{r} - \vec{r}_0) = 0 \nonumber \\
g(\vec{\,r}' - \vec{\,r}'_0) +  h(\vec{\,r}'_1- \vec{\,r}'_0) +  l(\vec{\,r}'_2 - \vec{\,r}'_0) = 0.]]
\end{eqnarray}

Putting $(\vec{r},\vec{\,r}')$, etc. in terms of $(\vec{\omega}, \Delta\vec{r})$, etc. we have finally
\begin{eqnarray}\label{ghl4ptsvecomDel}
g(\vec{\omega} - \vec{\omega}_0) + h(\vec{\omega}_1 - \vec{\omega}_0) + l(\vec{\omega}_2 - \vec{\omega}_0) = 0, \nonumber \\
g(\Delta\vec{r} - \Delta\vec{r}_0) + h(\Delta\vec{r}_1 - \Delta\vec{r}_0) + l(\Delta\vec{r}_2 - \Delta\vec{r}_0) = 0.
\end{eqnarray}

\noindent Here the upper line, compared with \eqref{tetmunupi}, shows that $g,h,l$ are proportional to $\mu,\nu,\pi$. Therefore
the former can be replaced by the latter in the lower line, yielding
\begin{equation}
\mu(\Delta\vec{r}-\Delta\vec{r}_0)+\nu(\Delta\vec{r}_1-\Delta\vec{r}_0)+\pi(\Delta\vec{r}_2-\Delta\vec{r}_0),
\end{equation} 

\noindent Q. E. D.}

\section{Coordinate changes from infinitesimal analysis.}

{\it And a third way!}

The two analytic methods that we have just presented for determining the formulas for the coordinate changes of solid system
are based on purely algebraic procedures.  From infinitesimal analysis we can derive an even simpler proof of these formulas,
resembling that given by Lagrange in his \emph{Analytic Mechanics}, and describing in the same analysis both the expression
for finite changes and that for infinitesimal ones.  Here is that demonstration.

Given a point $\vec{r}$ in the undisplaced system, consider a neighboring point $\vec{r} + d\vec{r}$ also before displacement, 
$d\vec{r}$ being infinitesimal.  Let the changes due to displacment be denoted by a prefix $\Delta$ or $\delta$ according as
the displacement is finite or infinitesimal.  We wish to determine the general expression for $\Delta\vec{r}$ by means of the
equation $\Delta |d\vec{r}|^2 = 0$.  

Now, letting $\vec{\omega} = \vec{r} + (1/2) \Delta\vec{r}$, and using the commutation of the signs $d$ and $\Delta$,  the foregoing
equation becomes 
\begin{equation}\label{domdr}
d\vec{\omega}\cdot d\Delta\vec{r}= 0,
\end{equation}

\noindent which is to be satisfied in the most general way {\it (that is, for arbitrary choices of $d\vec{r}$)}. 

{\it The proof of \eqref{domdr} is made easier by defining $\vec{r\,}' = \vec{r} + \Delta\vec{r}$.  Then
\begin{equation}\label{domdrpf}
\Delta |d\vec{r}|^2  = d\vec{r\,}'\cdot d\vec{r'} - d\vec{r}\cdot d\vec{r} = (d\vec{r\,}'+ d\vec{r})\cdot(d\vec{r\,}' - d\vec{r}) = 2d\vec{\omega}\cdot\Delta\vec{r},
\end{equation}

\noindent which establishes \eqref{domdr}.}

To this end, let us consider $\Delta\vec{r}$ as a function of $\vec{\omega}$;
then \eqref{domdr} can be written, taking $d\vec{\omega}$ as constant {\it (that is, as unchanged by the displacement)}, as 
$\Delta(d\vec{\omega}\cdot d\vec{r}) = 0$, from which follows
\begin{equation}\label{omdotDelconst}
\frac{d\vec{\omega}\cdot\Delta\vec{r}}{|d\vec{\omega}|} = constant
\end{equation}

\noindent {\it (in the same sense of \lq\lq constant")}.

This equation \eqref{omdotDelconst} represents algebraically the property of a quadrilateral two of whose sides are equal, that these 
two sides, and the two other sides as well, project {\bf equally} {\it (that is, each opposed pair makes a pair of equal projections)}
on the line joining the midpoints of the second pair of sides. 

 {\it This geometrical theorem needs some interpretation as well as a proof.  The two sides that are equal (in length) are the initial and final vectors
$d\vec{r}$ and $d(\vec{r}+\Delta\vec{r})$, since the two neighboring points are rigidly connected during the displacement. The other 
two sides are the lines traversed by the two neighboring points, which may be finite; these two lines need not be of equal lengths (although
the difference $d|\Delta\vec{r}|$ is infinitesimal) since different points can move differently under a given displacement.  There is no
requirement that the quadrilateral be all in one plane - this makes visualization even more challenging. I have been unable to find
a geometric proof simpler than the algebraic one already given leading to \eqref{omdotDelconst}. 

However, it may help to restate the theorem without infinitesimals, but with vector notation.  Let the four vertices of the quadrilateral be 
called $\vec{r}_1, (\vec{r}_1)', \vec{r}_2, (\vec{r}_2)'$, and let the opposite lengths $|\vec{r}_1 - \vec{r}_2|$ and $|(\vec{r}_1)' - (\vec{r}_2)'|$
be equal.  Let the second pair of opposite sides be called $\vec{\Delta}_1 = (\vec{r}_1)' - \vec{r}_1$ and $\vec{\Delta}_2 = (\vec{r}_2)' - \vec{r}_2$.
Let $\vec{\omega}_1,\vec{\omega}_2$ be the midpoints of these two sides.  Then the four vertices can be renamed
$\vec{\omega}_1 - (1/2)\vec{\Delta}_1, \vec{\omega}_1 + (1/2)\vec{\Delta}_1, \vec{\omega}_2 - (1/2)\vec{\Delta}_2, \vec{\omega}_2 + (1/2)\vec{\Delta}_2$.
The two equal lengths are now $|(\vec{\omega}_1 - \vec{\omega}_2) - (1/2)(\vec{\Delta}_1 - \vec{\Delta}_2)|$ and 
$|(\vec{\omega}_1 - \vec{\omega}_2) + (1/2)(\vec{\Delta}_1 - \vec{\Delta}_2)|$, from which it follows trivially that $(\vec{\omega}_1 - \vec{\omega}_2)$ is
orthogonal to $(\vec{\Delta}_1 - \vec{\Delta}_2)|$ as well as to $\vec{r}_1 - \vec{r}_2$, as stated by the theorem.}


Consider now the equation 
\begin{equation}\label{domdrxyz}
d\xi d\Delta x + d\eta d\Delta y + d\zeta d\Delta z = 0.
\end{equation}

\noindent {\it (This is eq \eqref{domdr} written out in components.)} We shall expand the complete differentials $d\Delta x$, etc. in terms
of partials $\partial$.  
{\it (That is, $d\Delta x = d\xi \frac{\partial \Delta x}{\partial\xi} + d\eta\frac{\partial \Delta x}{\partial\eta} + d\zeta \frac{\partial \Delta x}{\partial\zeta}$., etc.)}
\eqref{domdrxyz} thus becomes
\begin{eqnarray}\label{domdrxyzpart}
& & (d\xi)^2\frac{\partial \Delta x}{\partial\xi} + (d\eta)^2\frac{\partial \Delta y}{\partial\eta} + (d\zeta)^2\frac{\partial \Delta z}{\partial\zeta} 
+ d\xi d\eta (\frac{\partial \Delta x}{\partial\eta}  + \frac{\partial \Delta y}{\partial\xi})            \nonumber \\
& & + d\xi d\zeta (\frac{\partial \Delta x}{\partial\zeta} + \frac{\partial \Delta z}{\partial\xi})
+ d\eta d\zeta (\frac{\partial \Delta y}{\partial\zeta} + \frac{\partial \Delta z}{\partial\eta}) = 0.
\end{eqnarray}

\noindent Inasmuch as the differentials $d\xi,d\eta,d\zeta$ are independent, \eqref{domdrxyzpart} implies
\begin{eqnarray}\label{part0}
& & \frac{\partial \Delta x}{\partial\xi} = \frac{\partial \Delta y}{\partial\eta} = \frac{\partial \Delta z}{\partial\zeta} =0, \nonumber \\
& & \frac{\partial \Delta x}{\partial\eta} + \frac{\partial \Delta y}{\partial\xi} = \frac{\partial \Delta x}{\partial\zeta} + \frac{\partial \Delta z}{\partial\xi}
= \frac{\partial \Delta y}{\partial\zeta} + \frac{\partial \Delta z}{\partial\eta} = 0.
\end{eqnarray}

{\it By setting $d\eta = d\zeta = 0$, one establishes the vanishing of $\frac{\partial \Delta x}{\partial\xi}$, similarly 
$\frac{\partial \Delta y}{\partial\eta}$ and $\frac{\partial \Delta z}{\partial\zeta}$.  Then by constraining only $d\zeta$ to vanish, 
one isolates $\frac{\partial \Delta x}{\partial\eta} + \frac{\partial \Delta y}{\partial\xi}$, and so forth.}


This system of six equations can be easily integrated.  The first shows that $\Delta x$ is independent of $\xi$, but then
its derivatives $\frac{\partial\Delta x}{\partial\xi}$, $\frac{\partial\Delta x}{\partial\eta}$, $\frac{\partial\Delta x}{\partial\zeta}$
are also independent of $\xi$, and likewise starting with the second or third term, so that each of the nine partial derivatives
is independent of each of the independent variables $\xi.\eta,\zeta$; {\it (i.e., each is constant.)} Moreover, the matrix of derivatives
$\frac{\partial(\Delta x,y,z)}{\partial(\xi,\eta,\zeta)}$ is an antisymmetric matrix of constants, and this leads finally to the expressions
at which we have previously arrived,
\begin{equation}\label{npmdejavu}
\Delta x = A + p\eta - n\zeta,\;\;\Delta y = B + m\zeta - p\xi,\;\;\Delta z = C+ n\xi - m\eta.
\end{equation}

{\it Here the sense of \lq\lq constant" is that the quantity so called vanishes under the operator $d$, unlike its sense
in \eqref{omdotDelconst}. The point is that the equations \eqref{domdr} and \eqref{domdrxyz} remain true with no
change in form, if the \lq\lq neighboring point" $\vec{r} + d\vec{r}$ is replaced by a different neighboring point.
Therefore all the consequences can be differentiated again and again if one wishes. In
particular, \eqref{part0} can be differentiated so as to show that all 27 second partial derivatives are zero, and
consequently all nine first partial derivatives $\frac{\partial(\Delta x,y,z)}{\partial(\xi,\eta,\zeta)}$ are constant.
Then appealing again to \eqref{part0}, one sees that this constant matrix is antisymmetric.}

We shall not return to the transformations undergone by these formulas in reestablishing the variables $x,y,z$; it suffices
to have shown how the method of variations applies to the study of these formulas and gives in a single algebraic form 
both finite and infinitesimal coordinate changes in displaced points, these points being replaced in the case of finite
displacements by the midpoints of the interval traversed in a straight line by them.

\section{Recapitulation}

To wind up this work it remains to deduce quickly, from the expression for these changes, the geometric laws for the 
displacement of solid bodies that we developed synthetically in the first place, and took as point of departure for our
first analysis.

The formula
\begin{equation}\label{npmdejavuvec}
\vec{\Delta} = \vec{\Gamma}+ \vec{q}\times\vec{\omega}
\end{equation}

\noindent {\it (this is the vectorial form of \eqref{npmdejavu})} immediately gives the following fundamental relation:
\begin{equation}\label{qproj}
\vec{q}\cdot\vec{\Delta}= \vec{q}\cdot\vec{\Gamma}, 
\end{equation} 

\noindent from which one sees that the {\bf lines} actually traversed by all the points of the system,
in passing from one situation to the other, all have equal projections onto a particular direction $\hat{q}$.
{\it Recall that the components of $\vec{\Gamma}$ and those of $\vec{q}$ are those six constants
$A,B,\Gamma, m, n, p$ which characterize the displacement as a whole. Therefore, \eqref{qproj}
tells us that although different points 
$\vec{r}$ traverse different lines $\vec{\Delta}$, these lines all project equally on the fixed direction $\hat{q}$.}
Denoting this projection by $t$, one will have for all the points of the displaced system
\begin{equation}\label{tproj}
\hat{q}\cdot\vec{\Delta} = t
\end{equation}

\noindent and for two different points,
\begin{equation}\label{tprojsame} 
\hat{q}\cdot(\vec{\Delta}_1 - \vec{\Delta}_2) = 0.
\end{equation}

\noindent Here $\vec{\Delta}_1 - \vec{\Delta}_2$ is the chord of the arc that would be described by the
first point about an axis of rotation drawn through the second point and parallel to the direction $\hat{q}$;
the two points can thereafter be brought to their final positions by translating them both by $\vec{\Delta}_2$.

If $\theta$ represents the angle of that rotation, and $u$ the distance from that chord to that axis of rotation,
one has evidently
\begin{equation}\label{4tansqdejavu}
4u^2\tan^2(\theta/2) = |\vec{\Delta}_1 - \vec{\Delta}_2|^2 
= \vec{\,q}^2[|\vec{\omega}_1 - \vec{\omega}_2|^2 - ((\vec{\omega}_1 - \vec{\omega}_2))\cdot\hat{q})^2],
\end{equation}

\begin{equation}\label{usqdejavu}
u^2 = |\vec{\omega}_1 - \vec{\omega}_2|^2 - [(\vec{\omega}_1 - \vec{\omega}_2))\cdot\hat{q})]^2
\end{equation}

and therefore, regardless of what two points are being considered,
\begin{equation}\label{qsqdejavu}
4 \tan^2(\theta/2) = \vec{\,q}^2.
\end{equation} 

\noindent {\it The author has worked his way back to what followed immediately from the definition 
$\vec{q} = (m,n,p)$ given in Section 15.  It should be noted that $\hat{q}$ is the same unit vector as $\hat{t}$.}

\noindent Thus, the displacement given to a solid, from one situation to another, can always be resolved into
two consecutive displacements, one of rotation and one of translation, just as has been explained at the 
start of this treatise.

Furthermore, let $\nu$ represent the amplitude of the angular displacement of a \emph{line} within the solid,
and $\phi$ the angle formed by the axis of rotation with this line, then if $\vec{r}_1,\vec{r}_2$ are two points
on the line we have
\begin{equation}\label{cosphi}
\cos\phi = \frac{(\vec{r}_1- \vec{r}_2)\cdot\hat{q}}{|\vec{r}_1-\vec{r}_2|},
\end{equation}

\noindent and, in view of the invariablity of the distances between points of the solid,
\begin{equation}\label{Delcosphi}
\Delta\cos\phi = \frac{(\vec{\Delta}_1 - \vec{\Delta}_2)\cdot\hat{q}}{|\vec{r}_1-\vec{r}_2|} = 0.
\end{equation}

\noindent The angle $\phi$ is the same before and after the displacement.  As for the angle $\nu$,
we have
\begin{equation}\label{cosnu}
\cos\nu = \frac{(\vec{r}_1 - \vec{r}_2)\cdot(\vec{r}_1 - \vec{r}_2+ \vec{\Delta}_1- \vec{\Delta}_2)}{|\vec{r}_1-\vec{r}_2|^2},
\end{equation}

\noindent whence finally the remarkable relation 
\begin{equation}\label{sinhalfnudejavu}
\sin(\nu/2) = \sin\phi \sin(\theta/2) 
\end{equation}

\noindent expressing the theorem stated in {\bf 5}. {\it (Actually the relation is given in {\bf 6}.)}

{\it It is not trivial to derive \eqref{sinhalfnudejavu} from the preceding equations.  An essential
preliminary step is to derive the identity
\begin{equation}\label{om12Del12}
(\vec{\omega}_1 - \vec{\omega}_2)\cdot(\vec{\Delta}_1 - \vec{\Delta}_2) = 0.
\end{equation}

\noindent This is easily done by noting that for any point $\vec{r}$ we have $\vec{r} = \vec{\omega} - (1/2)\vec{\Delta}$
and $\vec{r} + \vec{\Delta} = \vec{\omega} + (1/2)\vec{\Delta}$, so that the equivalence of the two lengths $|\vec{r}_1 -  \vec{r}_2|$
and $|(\vec{r}_1 + \vec{\Delta}_1) -  (\vec{r}_2 + \vec{\Delta}_2)|$ can be written as
\begin{equation}\label{ompmDel}
|(\vec{\omega}_1 - (1/2)\vec{\Delta}_1) - (\vec{\omega}_2 - (1/2)\vec{\Delta}_2)|^2
 =      |(\vec{\omega}_1 + (1/2)\vec{\Delta}_1) - (\vec{\omega}_2 + (1/2)\vec{\Delta}_2)|^2,
\end{equation}

\noindent which is equivalent to \eqref{om12Del12}.

From \eqref{om12Del12} there follow the two useful relations
\begin{equation}
(\vec{\omega}_1 - \vec{\omega}_2)\cdot\vec{\Delta}_1 = (\vec{\omega}_1 - \vec{\omega}_2)\cdot\vec{\Delta}_2
\end{equation}

\noindent and
\begin{equation}\label{r12sqom12sq}
|\vec{r}_1 - \vec{r}_2|^2 = |\vec{\omega}_1 - \vec{\omega}_2|^2 + (1/4)|\vec{\Delta}_1 - \vec{\Delta}_2|^2.
\end{equation}

Now consider \eqref{cosnu}. The numerator can be written as 
\begin{equation}\label{cosnunum}
[(\vec{\omega}_1 - (1/2)\vec{\Delta}_1) - (\vec{\omega}_2 - (1/2)\vec{\Delta}_2)]
\cdot [(\vec{\omega}_1 + (1/2)\vec{\Delta}_1) - (\vec{\omega}_2 + (1/2)\vec{\Delta}_2)]
=  |\vec{\omega}_1 - \vec{\omega}_2|^2 - (1/4)|\vec{\Delta}_1 - \vec{\Delta}_2|^2
\end{equation}

\noindent in view of \eqref{om12Del12}, and the denominator as
\begin{equation}\label{cosnuden} 
|\vec{\omega}_1 - \vec{\omega}_2|^2 + (1/4)|\vec{\Delta}_1 - \vec{\Delta}_2|^2
\end{equation}

\noindent by \eqref{r12sqom12sq}.  Then
\begin{equation}\label{sinsqhalfnu}
\sin^2(\nu/2) = (1/2)(1-\cos\nu) = (1/4) \frac{|\vec{\Delta}_1 - \vec{\Delta}_2|^2}{|\vec{\omega}_1 - \vec{\omega}_2|^2 + (1/4)|\vec{\Delta}_1 - \vec{\Delta}_2|^2}
= (1/4) \frac{|\vec{\Delta}_1 - \vec{\Delta}_2|^2}{|\vec{r}_1 - \vec{r}_2|^2}.
\end{equation}.

In dealing with $\phi$, we recall that it is the angle between the line containing $\vec{r}_1$ and  $\vec{r}_2$ and the fixed direction $\hat{q}$, also called the
axis of rotation.  If we think of $\hat{q}$ as vertical,  $\vec{r}_1 - \vec{r}_2$ may be resolved into a vertical part $(\vec{r}_1 - \vec{r}_2)\cdot\hat{q}$ and
a horizontal part which we may call $(\vec{r}_1 - \vec{r}_2)_\perp$. Then $\cos\phi$ is given by \eqref{cosphi}, and for $\sin\phi$ we have
\begin{equation}\label{sinphi}
\sin\phi = \frac{|\vec{r}_1 - \vec{r}_2|_\perp}{|\vec{r}_1 - \vec{r}_2|}.
\end{equation}

Combining \eqref{sinphi} with \eqref{sinsqhalfnu}, we find
\begin{equation}\label{nuphi}
\frac{\sin^2(\nu/2)}{\sin^2\phi} = (1/4)\frac{|\vec{\Delta}_1 - \vec{\Delta}_2|^2}{|\vec{r}_1 - \vec{r}_2|_\perp^2}.
\end{equation}

To obtain $\sin^2(\theta/2)$, it is not sufficient to proceed directly from \eqref{qsqdejavu}, as this will involve a factor
$\vec{\,q}^2 = m^2+n^2+p^2$ which does not enter into \eqref{nuphi}.  Instead, consider the arc swept out by rotating
$\vec{r}_1$ an angle $\theta$ about the \lq\lq vertical" axis passing through $\vec{r}_2$. The radius of this arc is
$|\vec{r}_1 - \vec{r}_2|_\perp$, and the chord subtended by $\theta$ has length $|\vec{\Delta}_1 - \vec{\Delta}_2|$.
Therefore 
\begin{equation}\label{sinsqhalfth}
\sin^2(\theta/2) = (1/4)\frac{|\vec{\Delta}_1 - \vec{\Delta}_2|^2}{|\vec{r}_1 - \vec{r}_2|_\perp^2}.
\end{equation}

Combining \eqref{sinsqhalfth} with \eqref{nuphi}, we obtain \eqref{sinhalfnudejavu}.

%
%
%
%
%

|

The lines parallel to the direction of the axis of rotation are therefore transported parallel to themselves.
Among all these lines there is one that simply glides upon itself; for this line the change $\Delta\vec{r}$ is
evidently in the direction $\hat{q}$.  The equation of this line is therefore
\begin{equation}\label{tcent}
\Delta\vec{r} = \hat{q}t = \vec{t};
\end{equation}

\noindent and since
\begin{equation}\label{omegdejavu}
\vec{\omega} = \vec{r} + (1/2)\Delta\vec{r},
\end{equation}

\noindent we obtain the same equation already given in {\bf 16} for the \emph{central} axis of the displacement:
\begin{equation}\label{centdejavu}
\Gamma + \vec{q}\times\vec{r} = t\hat{q} = \vec{t}.
\end{equation}}

\section{Conclusion - General law of Statics.}

Geometry considers the displacements, finite or infinitesimal, of \emph{solid} bodies, brought about 
by the \emph{successive} action of causes or of forces capable of producing them.

Mechanics considers consecutive displacements of solid bodies, and more generally of 
\emph{arbitrary} systems of points, brought about by the \emph{simultaneous} and prolonged action of
causes or of forces capable of producing them.

Statics is that most elementary part of Mechanics in which one considers only the possibility
of infinitesimal or \emph{virtual} displacements of these systems, resulting from the
simultaneous and discontinuous action of those same causes.

{\it The idea seems to be that the distinction between Geometry and Mechanics disappears
when one considers only infinitesimal displacements in Geometry, or instantaneous 
applications of force in Mechanics.}

Geometry teaches that the displacement of a solid body reduces to a turning about
one or two fixed axes.  {\it (Two, if one wishes to eliminate translations.)}

It follows that if the forces that act {\bf simultaneously} on a solid system cannot impress
on it any rotation about any fixed axis whatever, these forces equilibrate or 
neutralize one another, and the body remains at rest.

These forces, considered {\bf separately}, can act only in two ways, either by tending to 
turn this solid body about a fixed axis, or by tending to displace a certain point
of the system, or more expressively to change the coordinates of that point.
That is the most general way to consider and to examine, in Mechanics, the
action of forces.

At any rate the law of equilibrium is identical in these two modes of thought,
as we shall see. {\it (Probably the two modes are that of Geometry and that
of Mechanics.)}

If any forces or causes of displacement tend successively, or else simultaneously
in passing to the limit (that is to say in passing from Geometry to Mechanics),
to impress on a solid {\bf elementary} or virtual rotations $\theta,\theta',\theta'',...$,
about given fixed axes, the law of equilibrium of these forces is that the
sum of the moments of these rotations should vanish relatively to any axis whatever,
This law is rendered algebraically by the equation
\begin{equation}\label{SigthD} 
\Sigma\theta D \sin{\nu} = 0,
\end{equation}

\noindent implying, on account of the indeterminacy of that arbitrary axis, six special
equations, which reduce to three when the solid system reduces to a point.

Let us now examine what happens when the forces acting simultaneously on the solid
are applied individually to various points of the system.  As any displacement refers
{\bf virtually} to a fixed axis of rotation, it will suffice to consider for each point the change
that can result, from the action of the forces affecting that point, in the coordinate of that
point {\bf orthogonal} to that fixed axis {\it (and also orthogonal to the perpendicular dropped
from the point to the axis)}, whose resistance is opposed to any change in its {\it (the point's)}
other rectangular coordinates.  {\it The fixed axis is understood to be {\bf mechanically} 
fixed: it cannot move either along its own length or perpendicularly to it.  The only way
the point can move without the axis moving is on the tangent to the circle it would describe
if the system rotates about the axis.}

Now, it is evident that {\it (wordy redundant passage omitted)}.

From another point of view, on account of the rigidity of the system, it is evident that
two forces {\bf equal} {\it (in magnitude)} will, in their action on any given point, 
be in equilibrium 
(1) if they are applied through the point in opposite directions;
(2) if they are applied in opposite sense through the extremities of a fixed line segment;
(3) if they tend to turn in opposite sense a circumference whose center is fixed, in the plane of
which they are applied tangentially;
(4) if, more generally, they tend to turn oppositely a right cylinder whose axis is fixed, at whose
surface they are applied tangentially and orthogonal to its axis.

It results from these propositions that if one considers all the forces that tend
to displace the individual points of a solid system which contains a fixed axis, and
acting in given directions, and to impress on them individually given {\bf virtual} 
translations, there will be equilibrium among all these forces if, supposing that
they are all applied to points equidistant from the fixed axis - which is always possible -
the sum of the {\bf moments} of the virtual translations that measure the effect of
these forces is zero relative to that axis. {\it For \lq\lq moment", read \lq\lq torque: the
total torque about the fixed axis should vanish, in order to produce equilibrium.  I don't
understand the language about points equidistant from the fixed axis.} 

Passing from a fixed axis to an arbitrary axis, it will follow necessarily that the general
equation \eqref{SigthD} is the algebraic expression of the equilibrium of a set of forces
capable of producing virtual or infinitesimal translations proportional to the rotations 
$\theta,\theta',\theta'',...$, those forces being applied about the axes of these rotations
positively or negatively according to the signs of the rotations.

This explains the remarkable analogy between the laws of equilibrium (and consequently 
of the composition) of infinitesimal rotations and the laws of equilibrium and composition
of forces, considered in Statics.  {\it This is the analogy discussed in {\bf 26}.}

If one denotes these forces by their  finite magnitudes $P,P',...,$,the equation of their equilibrium will be
\begin{equation}\label{SigPDnu}
\Sigma PD\sin\nu = 0.
\end{equation}

\noindent Here each term $PD\sin\nu$ expresses the static moment  [[torque]] of the
force $P$ about the  fixed axis; it is equal to the product of the distance from the point
of application of that force to the fixed axis with the component of that force normal
to that axis. {\it The author here combines $P$ and $\nu$ into a single quantity, a component
of the force. He should say also, normal to the distance $D$.

Repetitive paragraph omitted.}

The conditions of equilibrium of forces applied to a solid system, which the secondary
but admirably ingenious consideration of {\bf couples} reduces in finite terms to
two conditions, are therefore comprised in a single law, similarly expressed in finite terms,
that the sum of the moments of the forces be null about an arbitrary axis.
This law is general and applies, as does that of the principle of virtual speeds - equivalent
to an infinitesimal transformation -  to the equilibrium of {\bf any system rigid or not}, provided
that the conditions of the connections among points of the system be replaced by the introduction
of forces that make it possible to regard these points as entirely free.  {\it Here the author
definitely is speaking of Mechanics; his idea is that even if the system is not solid but consists
of independently mobile points, it can be made to act like a solid if the forces acting on it
are such as to leave invariant all the distances $|\vec{r}_i - \vec{r}_j|$.} These forces are
determined by analysis and eliminated in accordance with the equations that give the equilibrium
of each point.

By this means the law of equilibrium of a point immediately implies that of a solid system.
We shall not linger over this. We merely remark that in the particular case of the
forces introduced being equal and opposite in pairs, the sum of the moments of
all the forces applied to all the points, which must be null for equilibrium, contains
only the sum of the given forces, and thus expresses, in identical form, the law
of equilibrium of one or many points entirely free and that of a rigid system.

\subsection{On the equation of virtual speeds.}

Consider an infinitesimal displacement of the system, producing changes characterized
by the symbol $\delta$ in the coordinates. If the force $P$ tends to change the coordinate $p$,
{\it (we claim that)} the product $P\delta p$ will be equal to $PD\theta\sin\nu$, where $D\theta\sin\nu$ is 
the {\bf moment} of the virtual rotation $\theta$.

In fact, since this infinitesimal displacement must reduce to either a single or to two successive rotations,
it need only be considered in its unique rotation or in the first of the two.  And so the infinitesimal arc 
described by the point at which $P$ is applied, projected along the direction of that force, will be
equal to the infinitesimal change in the coordinate $p$ on which the force acts. This {\it (indeed)} gives
$\delta p = \theta D\sin\nu$, or
\begin{equation}\label{PDdp}
P\delta p = \theta PD\sin\nu.
\end{equation}

\noindent Hence the general equation \eqref{SigPDnu} of equilibrium of forces

\noindent transforms into
\begin{equation}\label{SigPdp}
\Sigma P\delta p = 0,
\end{equation}

\noindent which says that {\bf if a solid system is in equilibrium under a set of forces,
and by any cause this system is infinitesimally dislodged from its present position, 
the sum of the products of each force with the infinitesimal distance traversed by
the point {\it (of application of the force)} of the system along the direction of that force
must be zero, and conversely; this is the principle of virtual speeds.}{\it As we would
say nowadays, the solid is in equilbrium under a set of forces applied at certain points
in a certain manner if and only if the {\bf virtual work} that would be done by these forces
in {\bf any} hypothetical infinitesimal displacement {\bf obeying the geometrical constraints
of the solid} is zero.}

Equation \eqref{SigPdp}, although it is certainly superior, algebraically speaking, 
to \eqref{SigPDnu}, is no more general at bottom; but it expresses in the simplest
possible way the law of equilibrium of any system in which the conditions binding
the parts together can be transformed into a set of linear equations among the
changes of coordinates of the various points of the system.

{\it The remainder of this section, and hence of the whole treatise, is devoted to
elucidating the above remark.}

In fact, {\it (suppose that)} these conditions are expressed by equations such as
$\delta L= 0, \delta L'= 0, \delta L''= 0,$ etc. or more generally by a single 
equation such as 
\begin{equation}\label{SiglamdelL}
\lambda\delta L + \lambda'\delta L' + \lambda''\delta L'' + ... =0,
\end{equation}
in which $\lambda,\lambda',\lambda'',...$ are arbitrary multipliers {\it (what we
nowadays call Lagrange multipliers)}. Let us denote
by $\bar{x},\bar{y},\bar{z}$ the coefficients of the changes $\delta x,\delta y, \delta z$ {\it (hidden within)}
this equation; by $r$ a linear distance in the direction of $(\bar{x},\bar{y},\bar{z})$; and by $R$
a force equal to $\sqrt{\bar{x}^2+\bar{y}^2+\bar{z}^2}$, applied at the point $(x,y,z)$ 
in the direction of the line that it tends to change, Then we have
\begin{equation}
\delta r = \frac{\bar{x}\delta x + \bar{y}\delta y + \bar{z}\delta z}{\sqrt{\bar{x}^2+\bar{y}^2+\bar{z}^2}}
\end{equation}

\noindent and similarly for $\delta r'$,$\delta r''$, etc.  Equation \eqref{SiglamdelL} will
now take the form $R\delta r + R'\delta r' + R''\delta r'' + ... = 0$, and the equations \eqref{SigPdp}
and 
\begin{equation}\label{SigPp+Rr}
\Sigma P\delta p + R\delta r + R'\delta r' + R''\delta r'' + ... = 0
\end{equation}

\noindent will have equal generality, the changes $\delta p$ being \emph{limited} in
\eqref{SigPdp} by the equation of conditions \eqref{SiglamdelL}, and in \eqref{SigPp+Rr}
being completely {\bf independent}.  Now, in the second case, \eqref{SigPp+Rr} expresses 
the condition of equilibrium of all the points of the system, independent of all binding, but
implied by the external forces $P,P'P'',...$ and by other forces $R,R',R''...$ which \emph{statically}
have replaced the assumed conditions of binding. 

{\it The relation between \eqref{SiglamdelL} and \eqref{SigPp+Rr} is the familiar one that arises
in the Lagrange multiplier method: one relaxes the rigidity conditions requiring $\delta L$, etc.
to vanish, at the cost of introducing the extra terms $R\delta r$, etc. into the equation to be
solved.  It is interesting that the whole procedure, nowadays presented as a purely mathematical
manipulation in the spirit of \cite{Lagmech}, is here explicated in a totally
physical way.}

By eliminating these forces $R,R',R''...$ from \eqref{SigPp+Rr} one can obtain the definitive 
equations of equilibrium of the external forces subject to the binding of the system.  And conversely, 
if these equations hold, there will be equilibrium, since the external forces are then determined,
and the equations \eqref{SigPp+Rr} on independent values of $\delta r$, etc. establish the immobility of 
all the points of the system resulting from the action of the external forces plus that of the others
that are {\bf statically} equivalent to the given binding among the various points of the system.

When the system under consideration is {\bf continuous}, the equations of binding contain
{\bf definite integrals} which represent in some way an infinite number of linear conditions
among the changes of coordinates of the system.  The arbitrary multipliers may be moved
inside the integral sign, and it then remains to solve for the changes by a method entirely
analytic and independent of any static consideration.  One easily arrives at the following
general formula:
\begin{equation}\label{nonlin}
\delta S^n Udx_1dx_2dx_3...dx_n = S^n dx_1dx_2dx_3...dx_n [\delta U + 
U(\frac{\partial\delta x_1}{\partial x_1}+\frac{\partial\delta x_2}{\partial x_2}+...+\frac{\partial\delta x_n}{\partial x_n})],
\end{equation}

\noindent where $U$ is an arbitrary function of the independent variables $x_1,x_2,...,x_n$, and the symbol
$S^n$ denotes a multiple definite integral of order $n$.

{\it In today's notation \eqref{nonlin} could be written as
\begin{equation}\label{nonlinmod}
\delta\int dx_1...\int dx_n U(x_1,...x_n) = \int dx_1...\int dx_n [\delta U(x_1,...x_n) + U(x_1,...x_n)(\frac{\partial\delta x_1}{\partial x_1}
 +... +\frac{\partial\delta x_n}{\partial x_n})].
 \end{equation}
 
\noindent The point of interest is that $U$ is an {\bf arbitrary} function, not necessarily linear in the $x_i$.  This is why the term $\delta U$ enters
the right side of the equation.}

\emph{December 5, 1840.}

\section{Acknowledgment}

This translation has been truly a labor of love.  I could never have dared to undertake it, much less carry it through to completion, without the unflagging interest
and support of my associate Dr. Johannes Familton, with whom I discussed every part of the translation as it took shape.  Dr. Familton also gave me indispensable
help in designing and implementing the figures.

Figures

1) Section 3 (Euler)

2) Section 5 (Central axis)

3) Section 6 (spherical geodesic)

4) Section 10 (parallel axes)

\begin{figure}
   \includegraphics[width=0.88\textwidth]{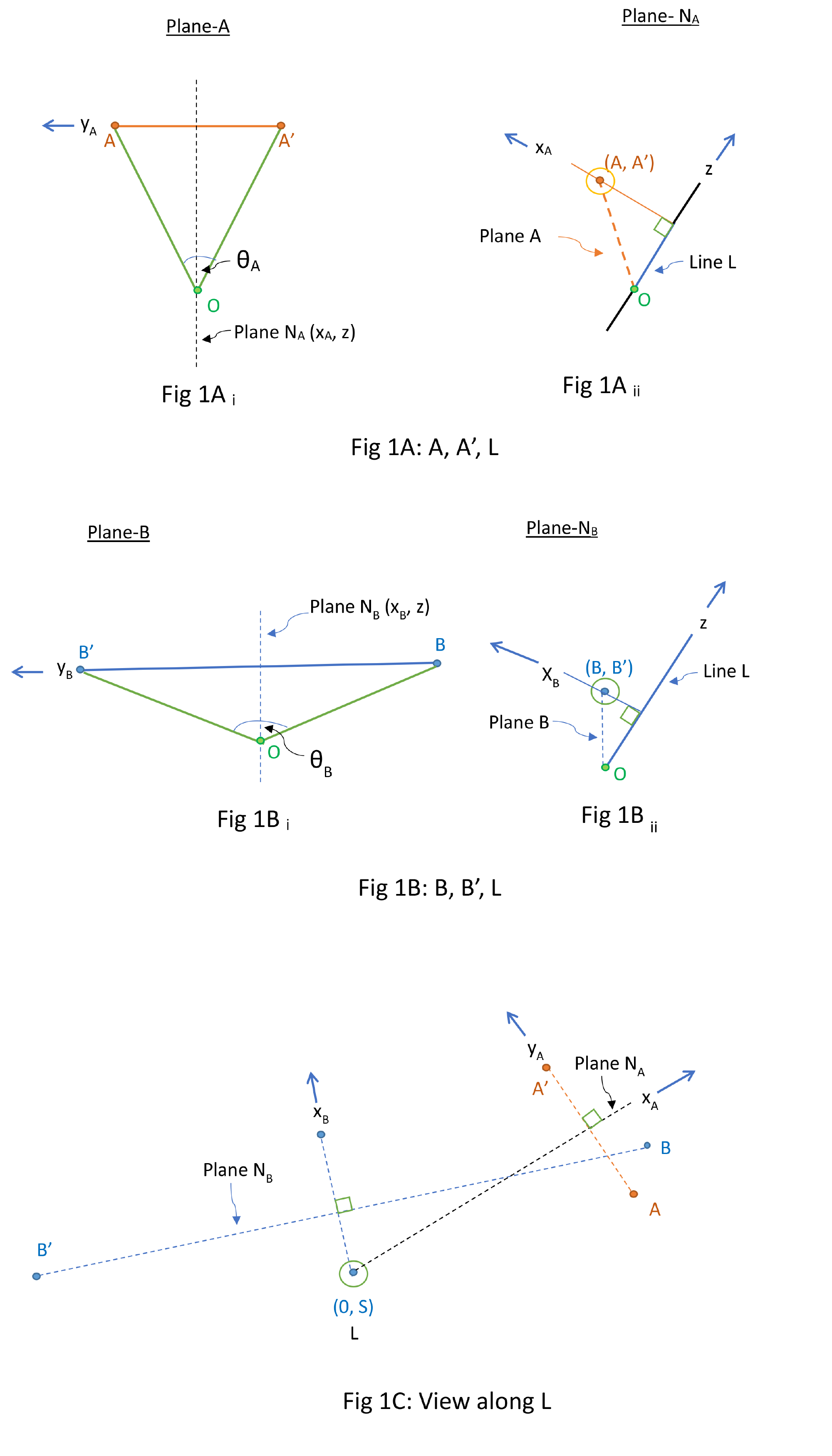}
\end{figure}

\begin{figure}
   \includegraphics[width=0.88\textwidth]{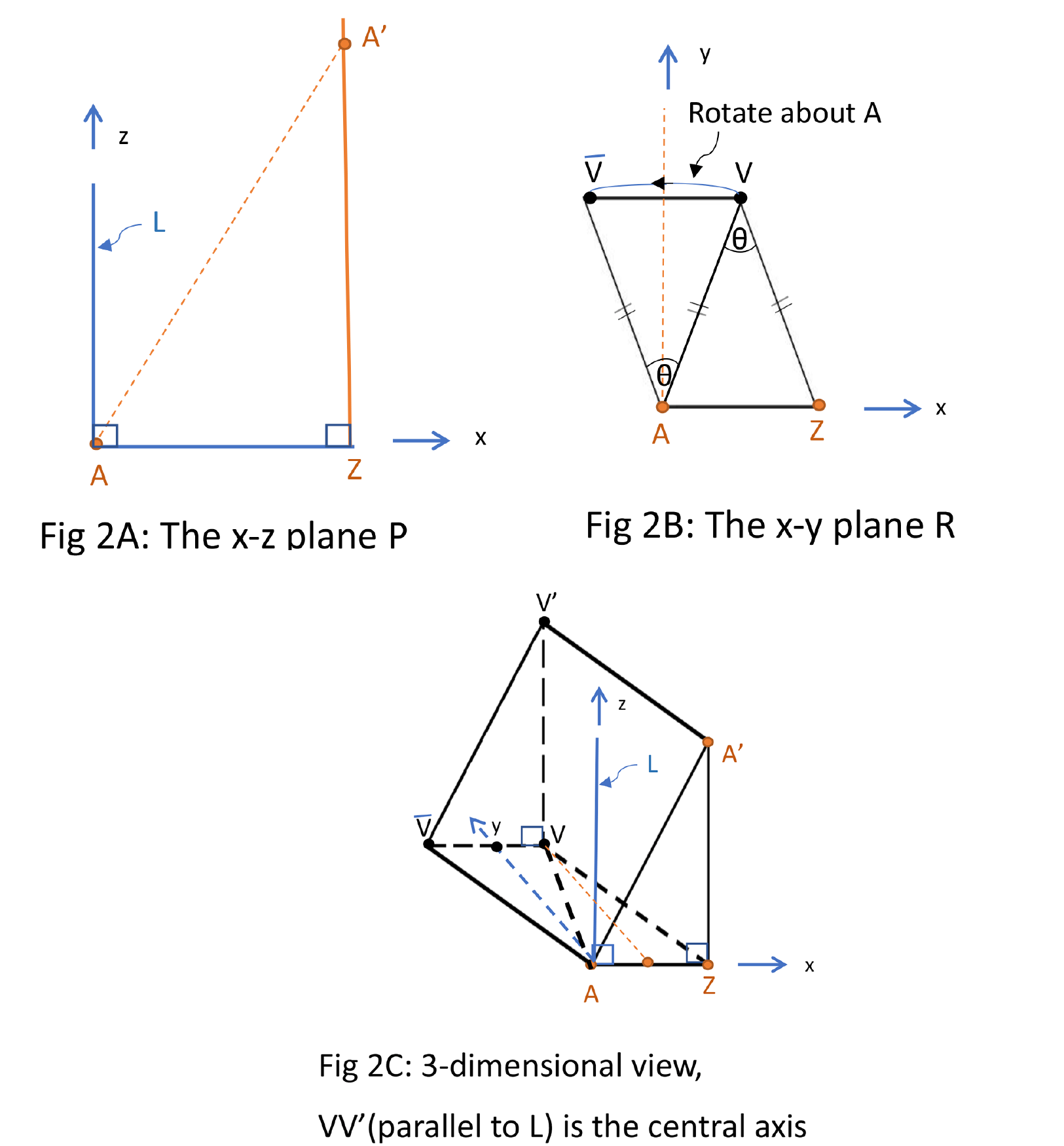}
\end{figure}

\begin{figure}
   \includegraphics[width=0.88\textwidth]{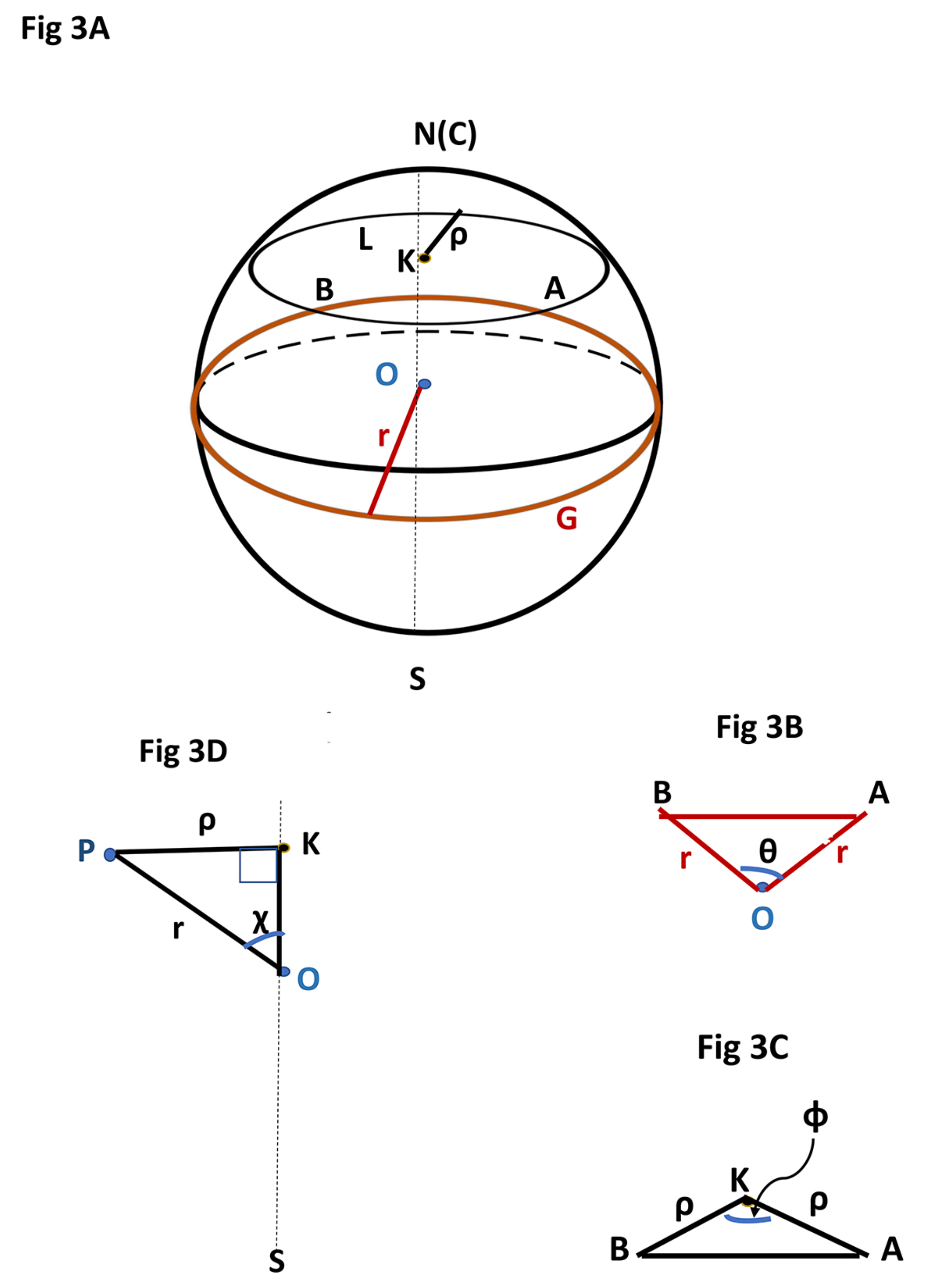}
\end{figure}

\begin{figure}
   \includegraphics[width=0.88\textwidth]{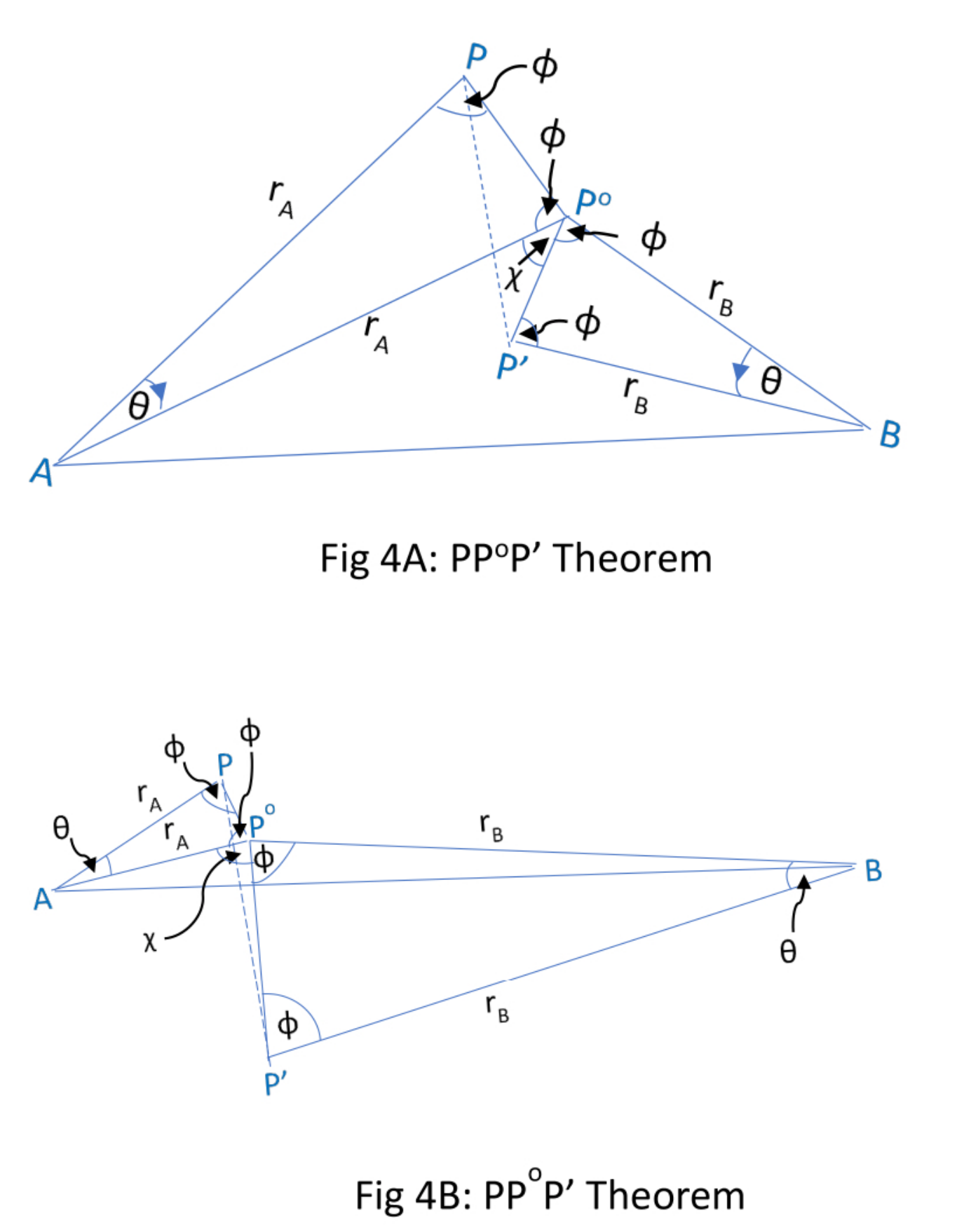}
\end{figure}

\end{document}